\documentclass[letterpaper, 12pt, oneside]{book}

\usepackage[left=1.5in, right=1in, top=1in, bottom=1in, dvips,pdftex]{geometry}
\usepackage{color,setspace, graphicx, verbatim,mathrsfs, amsmath, amsthm, amssymb, amsfonts,ntabbing,multirow,moreverb}
\newtheorem{defn}{Definition}[section]

\newtheorem{prop}[defn]{Proposition}
\newtheorem{lemma}[defn]{Lemma}
\newtheorem{thm}[defn]{Theorem}
\newtheorem{cor}[defn]{Corollary}

\theoremstyle{remark}

\numberwithin{equation}{section}
\newcommand{\com}[2]{\left[\,#1\,,#2\,\right]}
\newcommand{\inprod}[2]{\left(\,#1\,,#2\,\right)}

\DeclareMathOperator{\re}{Re}

\DeclareMathOperator{\sgn}{sgn}
\DeclareMathOperator{\pfaffian}{Pf} \DeclareMathOperator{\tr}{Tr}
\DeclareMathOperator{\prob}{P} \DeclareMathOperator{\airy}{Ai}
\newcommand{\rchi}{\raisebox{.4ex}{$\chi$}}
\def\ra{\rightarrow}
\def\iy{\infty}

\def\be{\begin{equation}}
\def\ee{\end{equation}}

\begin{document}

\singlespacing


\pagenumbering{roman}
\pagestyle{plain}


\begin{center}

\begin{Large}
  Distribution Functions for Edge Eigenvalues in Orthogonal and Symplectic Ensembles: Painlev\'e Representations
\end{Large}\\
\bigskip
By\\
\bigskip
MOMAR DIENG\\
B.A. (Macalester College, St Paul) 2000\\
M.A. (University of California, Davis) 2001\\
\bigskip
DISSERTATION\\
\bigskip
Submitted in partial satisfaction of the requirements for the degree of\\
\bigskip
DOCTOR OF PHILOSOPHY\\
\bigskip
in\\
\bigskip
MATHEMATICS\\
\bigskip
in the\\
\bigskip
OFFICE OF GRADUATE STUDIES\\
\bigskip
of the\\
\bigskip
UNIVERSITY OF CALIFORNIA,\\
\bigskip
DAVIS\\
\vspace{.75in}
Approved:\\
\bigskip
Craig A. Tracy
\par
\smallskip
\rule{4in}{1pt}\\
\bigskip
Bruno Nachtergaele
\par
\smallskip
\rule{4in}{1pt}\\
\bigskip
Alexander Soshnikov
\par
\smallskip
\rule{4in}{1pt}\\
\bigskip
\bigskip
Committee in Charge\\
\bigskip
\copyright\, Momar Dieng, MMV. All rights reserved.\\

\end{center}


\renewcommand{\baselinestretch}{1.6}\small\normalsize


\newpage

\tableofcontents


\newpage


\chapter*{}

\singlespacing
\begin{flushright}
Momar Dieng
\par
June 2005
\par
Mathematics
\end{flushright}
\par
\bigskip
\centerline{\textbf{\underline{Abstract}}}
\par
\bigskip
\doublespacing
We derive Painlev\'e--type expressions for the distribution of the $m^{th}$ largest eigenvalue in the Gaussian Orthogonal and Symplectic Ensembles in the edge scaling limit. This work generalizes to general $m$ the $m=1$ results of Tracy and Widom \cite{Trac2}. The results of Johnstone and Soshnikov (see \cite{John1}, \cite{Sosh2}) imply the immediate relevance of our formulas for the $m^{th}$ largest eigenvalue of the appropriate Wishart distribution.


\newpage
\chapter*{Acknowledgments}
\enlargethispage{\baselineskip}
My parents have been my first and best teachers. None of this could have existed without their love, the priceless upbringing they gave me, the sense of values, respect, and self they have instilled in me, and continue to nurture. It is my pleasure to dedicate this thesis to them in the hope that it will help justify my prolonged absence, and count as a very modest reward toward the countless sacrifices they made and are making.
\par\noindent
I could not have wished for a better thesis advisor than Professor Craig Tracy. This thesis is in great part the fruit of his patient support, guidance and expertise. Thank you Craig, not only for introducing me to Random Matrix Theory and so much beautiful Mathematics, but also for being a model of kindness, generosity and class.
\par\noindent
I would also like to thank Professor Alexander Soshnikov and Professor Bruno Nachtergaele for always pushing me to be a better mathematician.
\par\noindent
The staff of the Mathematics department at UC Davis has been instrumental in making this thesis possible, especially Celia Davis who has been a second mother to me. Thank you for your patience and professionalism.
\par\noindent
Last but not least, I would be remiss if I did not at least mention the many friends along the way who deserve more credit than I could possibly give here ...
\par\noindent
 This work was supported in part by the National Science Foundation under grant DMS-0304414.


\newpage
\pagestyle{headings}
\pagenumbering{arabic}

\chapter{Introduction}

\section{Motivation}

\begin{figure}[htbp] \label{TWdensities}
  \includegraphics{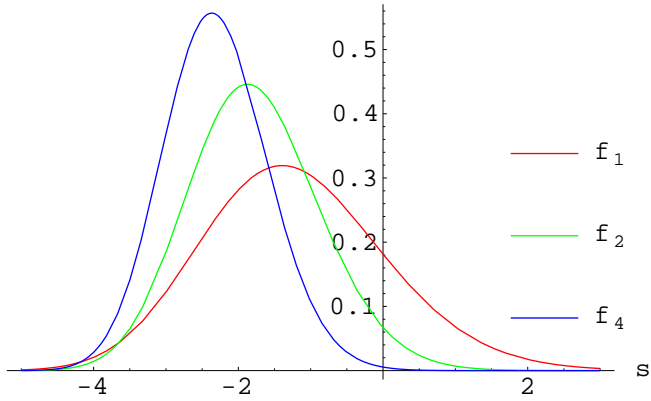}
  \caption{Tracy--Widom Density Functions}
\end{figure}

The Gaussian $\beta$--ensembles are probability spaces on $N$-tuples of random variables $\{\ell_{1},\ldots, \ell_{N}\}$, with joint density functions $P_{\beta}$ given by \footnote{In many places in the Random Matrix Theory literature, the parameter $\beta$ (times $1/2$) appears in front of the summation inside the exponential factor, in addition to being the power of the Vandermonde determinant. That convention originated in \cite{Meht1}, and was justified by the alternative physical and very useful interpretation of \eqref{jointdensity} as a one--dimensional Coulomb gas model. In that language the potential $W=\frac{1}{2}\sum_{i} l_{i}^{2}-\sum_{i<j}\ln|l_{i}-l_{j}|$ and $P_{\beta}^{(N)}(\vec{\ell}\,\,)=C\exp(-W/kT)=C\exp(-\beta\,W)$, so that $\beta=(k\,T)^{-1}$ plays the role of inverse temperature. However, by an appropriate choice of specialization in Selberg's integral (see Section~\ref{sec:selbergs-integral}), it is possible to remove the $\beta$ in the exponential weight, at the cost of redefining the normalization constant $C_{\beta}^{(N)}$. We choose the latter convention in this work since we will not need the Coulomb gas analogy. Moreover, with computer simulations and statistical applications in mind, this will in our opinion make later choices of standard deviations, renormalizations, and scalings a lot more transparent. It also allows us to dispose of the pesky square root of $2$ that is often present in the definition of the Tracy--Widom $\beta=4$ distribution in the literature. More on that later.}
\begin{equation}\label{jointdensity}
  P_{\beta}(\ell_{1},\ldots,\ell_{N})= P_{\beta}^{(N)}(\vec{\ell}\,\,) = C_{\beta}^{(N)}\,\exp\left[-\sum_{j=1}^{N}\ell_{j}^{\,2}\right]\prod_{j<k}|\ell_{j}-\ell_{k}|^{\beta}.
\end{equation}
The $C_{\beta}^{(N)}$ are normalization constants, given by
\begin{equation}
  \label{eq:34}
C_{\beta}^{(N)} = \pi^{-N/2}\,2^{-N-\beta\,N(N-1)/4}\cdot\prod_{j=1}^{N}\frac{\Gamma(1 + \gamma)\,\Gamma(1 + \frac{\beta}{2})}{\Gamma(1 + \frac{\beta}{2}\,j)}
\end{equation}
By setting $\beta=1, 2, 4$ we recover the (finite $N$) \emph{Gaussian Orthogonal Ensemble} ($\textrm{GOE}_{N}$), \emph{Gaussian Unitary Ensemble} ($\textrm{GUE}_{N}$), and \emph{Gaussian Symplectic Ensemble} ($\textrm{GSE}_{N}$), respectively. We restrict ourselves to those three cases in this thesis, and refer the reader to \cite{Dumi1} for recent results on the general $\beta$ case. Originally the $\ell_{j}$ are eigenvalues of randomly chosen matrices from corresponding matrix ensembles (see Section~\ref{cha:matrix-ensembles}), so we will henceforth refer to them as eigenvalues. With the eigenvalues ordered so that $\ell_{j}\geq\ell_{j+1}$, define
\begin{equation}
  \hat{\ell}_{m}^{(N)}=\frac{\ell_{m}-\sqrt{2\,N}}{2^{-1/2}\,N^{-1/6}},
\end{equation}
to be the rescaled $m^{th}$ eigenvalue measured from edge of spectrum. A standard result of Random Matrix Theory about the distribution of the largest eigenvalue in the $\beta$--ensembles (proved only in the $\beta=1,2,4$ cases) is that
\begin{equation}
  \hat{\ell}_{1}^{\,(N)}\xrightarrow{\mathscr{D}}\hat{\ell}_{1},
\end{equation}
whose law is given by the Tracy--Widom distributions.
\begin{thm}[Tracy, Widom \cite{Trac3},\cite{Trac2}]
  \begin{equation}\label{guemax}
    F_{2}(s):=\prob_{_{\textrm{GUE}}}(\hat{\ell}_{1}\leq s)=\exp\left[-\int_{s}^{\infty}(x-s)\,q^{2}(x)d\,x\right],
  \end{equation}
  \begin{equation}\label{goemax}
    F_{1}^{2}(s):=\left[\prob_{_{\textrm{GOE}}}(\hat{\ell}_{1}\leq s)\right]^{2}=F_{2}\cdot\exp\left[-\int_{s}^{\infty}q(x)d\,x\right],
  \end{equation}
  \begin{equation}\label{gsemax}
    F_{4}^{2}(\frac{s}{\sqrt{2}}):=\left[\prob_{_{\textrm{GSE}}}(\hat{\ell}_{1}\leq s)\right]^{2}=F_{2}\cdot\cosh^{2}\left[-\frac{1}{2}\int_{s}^{\infty}q(x)d\,x\right].
  \end{equation}
\end{thm}
The function $q$ is the unique (see \cite{Hasti1},\cite{Clar1}) solution to the Painlev\'e II equation
\begin{equation}\label{pII}
  q'' = x\,q + 2\,q^{3},
\end{equation}
such that $q(x)\sim \airy(x)$ as $x\to\infty$, where $\airy(x)$ is the solution to the Airy equation which decays like $\frac{1}{2}\,\pi^{-1/2}\,x^{-1/4}\,\exp\left(-\frac{2}{3}\,x^{3/2}\right)$ at $+\infty$. The density functions $f_{\beta}$ corresponding to the $F_\beta$ are graphed in Figure \ref{TWdensities}.\footnote{The square root of 2 in the argument of $F_{4}$ reflects a normalization chosen in \eqref{jointdensity} to agree with Mehta's original one (see \cite{Meht1}. See first footnote in this Introduction, as well as comments following Equation~\eqref{eq:26}, and Section~\ref{sec:stand-devi-matt}).}
Let $F_{2}(s,m)$ denote the distribution for the $m^{th}$ largest eigenvalue in GUE. Tracy and Widom showed in \cite{Trac3} that if we define $F_{2}(s,0)\equiv 0$, then
\begin{equation}
  F_{2}(s,m+1) -  F_{2}(s,m) = \frac{(-1)^{m}}{m\,!}\frac{d^{m}}{d\,\lambda^{m}}\,D_{2}(s,\lambda)\big{\vert}_{\lambda=1}\,,\quad m\geq 0, \label{lderiv}
\end{equation}
where
\begin{equation}\label{D2}
  D_{2}(s,\lambda)=\exp\left[-\int_{s}^{\infty}(x-s)\,q^{2}(x,\lambda)d\,x \right],
\end{equation}
and $q(x,\lambda)$ is the solution to \eqref{pII} such that $q(x,\lambda)\sim~\sqrt{\lambda}\,\airy(x)$ as $x\to\infty$. An intermediate step leading to \eqref{D2} is to first show that $D_{2}(s,\lambda)$ can be expressed as a Fredholm determinant
\begin{equation}
  D_{2}(s,\lambda)=\det(I-\lambda\,K_{\airy}),\label{fred}
\end{equation}
where $K_{\airy}$ is the integral operator on $L^{2}(s,\infty)$ with kernel
\begin{equation}\label{kairy}
  K_{\airy}(x,y)=\frac{\airy(x)\airy^{\prime}(y)-\airy^{\prime}(x)\airy(y)}{x-y}\,.
\end{equation}
In the $\beta=1,4$ cases a result similar to \eqref{fred} holds with the difference that the operators in $D_{\beta}(s,\lambda)$ have matrix--valued kernels (see e.g. \cite{Trac5}). In fact, the same combinatorial argument used to obtain the recurrence \eqref{lderiv} in the $\beta=2$ case also works for the $\beta=1,4$ cases, leading to
\begin{equation}\label{betaderiv}
  F_{\beta}(s,m+1) -  F_{\beta}(s,m) = \frac{(-1)^{m}}{m\,!}\frac{d^{m}}{d\,\lambda^{m}}\,D_{\beta}^{1/2}(s,\lambda)\big{\vert}_{\lambda=1}\,,\quad m\geq 0,\, \beta=1,4,
\end{equation}
where $F_{\beta}(s,0)\equiv 0$. Given the similarity in the arguments up to this point and comparing \eqref{D2} to \eqref{guemax}, it is natural to conjecture that $D_{\beta}(s,\lambda), \beta=1,4,$ can be obtained simply by replacing $q(x)$ by $q(x,\lambda)$ in \eqref{goemax} and \eqref{gsemax}. However the following fact, of which Corollary~\eqref{interlacingcor} gives a new operator theoretic proof, hints that this cannot be the case.
\begin{thm}[Baik, Rains \cite{Baik2}]\label{baikthm}
  In the appropriate  scaling limit, the distribution of the largest eigenvalue in GSE corresponds to that of the second largest in GOE. More generally, the joint distribution of every second eigenvalue in the GOE coincides with the joint distribution of all the eigenvalues in the GSE, with an appropriate number of eigenvalues.
\end{thm}
This so--called ``interlacing property'' between GOE and GSE had long been in the literature, and had in fact been noticed by Mehta and Dyson (see \cite{Meht3}). In this context, the remarkable work of Forrester and Rains in \cite{Forr1} classified all weight functions for which alternate eigenvalues taken from an Orthogonal Ensemble form a corresponding Symplectic Ensemble, and similarly those for which alternate eigenvalues taken from a union of two Orthogonal Ensembles form a Unitary Ensemble. Theorem~\eqref{baikthm} does not agree with the formulae we postulated for $D_{\beta}(s,\lambda), \beta=1,4$. Indeed, combining the two leads to incorrect relationships between partial derivatives of $q(x,\lambda)$ with respect to $\lambda$ evaluated at $\lambda=1$ and $q(x,1)$. To be precise, the conjecture is true for $D_{4}(s,\lambda)$ but it is false for $D_{1}(s,\lambda)$. The correct forms for both $D_{\beta}(s,\lambda), \beta=1,4$ are given below in Theorem~\eqref{mainthm}. \par This work also extends that of Johnstone in \cite{John1} (see also \cite{Elka1}), since $F_{1}(s,m)$ gives the asymptotic behavior of the $m^{th}$ largest eigenvalue of a $p$ variate  Wishart distribution on $n$ degrees of freedom with identity covariance. This holds under very  general conditions on the underlying distribution of matrix entries by Soshnikov's universality theorem (see \cite{Sosh2} for a precise statement). In Table~\ref{table}, we compare our distributions to finite $n$ and $p$ empirical Wishart distributions as in \cite{John1}.

\section{Statement of the Main Results}

\begin{thm}\label{mainthm}
  In the edge scaling limit, the distributions for the $m^{th}$ largest eigenvalues in the \textrm{GOE} and \textrm{GSE} satisfy the recurrence \eqref{betaderiv} with\,\footnote{See first footnote in this Introduction, as well as comments following Equation~\eqref{eq:26}, and Section~\ref{sec:stand-devi-matt} for an explanation of why we write $D_{4}(s,\lambda)$ instead of $D_{4}(s/\sqrt{2},\lambda)$ as is customary in the RMT literature.}
  \begin{equation}\label{goedet}
    D_{1}(s,\lambda)=D_{2}(s,\tilde{\lambda})\,\frac{\lambda - 1 - \cosh{\mu(s,\tilde{\lambda})} + \sqrt{\tilde{\lambda}}\,\sinh{\mu(s,\tilde{\lambda})}}{\lambda - 2},
  \end{equation}
  \begin{equation}\label{gsedet}
    D_{4}(s,\lambda)=D_{2}(s,\lambda)\,\cosh^{2}\left(\frac{\mu(s,\lambda)}{2}\right),
  \end{equation}
  where
  \begin{equation}
    \mu(s,\lambda):=\int_{s}^{\infty}q(x,\lambda)d\,x, \qquad \tilde{\lambda}:=2\,\lambda-\lambda^{2},
  \end{equation}
  and $q(x,\lambda)$ is the solution to \eqref{pII} such that $q(x,\lambda)\sim~\sqrt{\lambda}\,\airy(x)$ as $x\to\infty$.
\end{thm}
\begin{cor}[Interlacing property]
  \begin{equation}\label{interlacingcor}
    F_{4}(s,m)= F_{1}(s,2\,m),\quad m\geq 1.
  \end{equation}
\end{cor}
In the next section we give the proof of these theorems. In the last, we present an efficient numerical scheme to compute $F_{\beta}(s,m)$ and the associated density functions $f_{\beta}(s,m)$. We implemented this scheme using MATLAB\texttrademark\,,\footnote{MATLAB\texttrademark\, is a registered trademark of The MathWorks, Inc., 3 Apple Hill Drive, Natick, MA 01760-2098; Phone: 508-647-7000; Fax: 508-647-7001.} and compared the results to simulated Wishart distributions.
\begin{center}
  \begin{figure}
    \includegraphics[height=65mm]{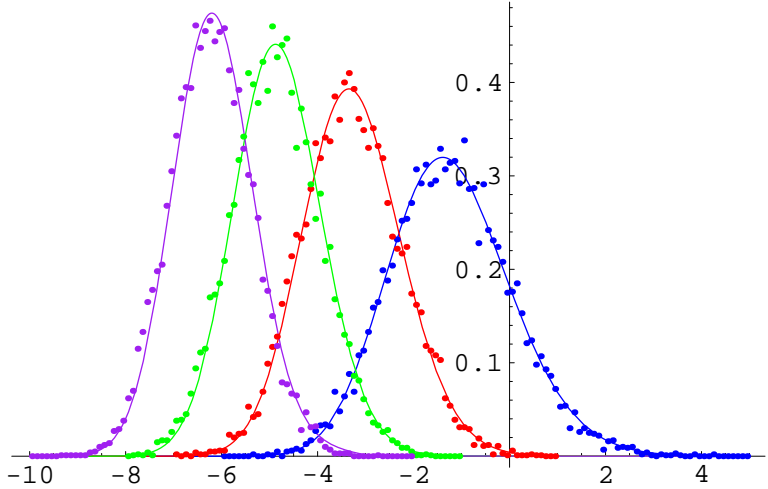}
    \caption{$10^{4}$ realizations of $10^{3}\times 10^{3}$ GOE matrices; the solid curves are, from right to left, the theoretical limiting densities for the first through fourth largest eigenvalue.}
  \end{figure}
\end{center}

\chapter{Preliminaries}

\section{Determinant matters}

We gather in this short chapter more or less classical results for further reference. Since these are well-known results, their proofs will be omitted. The exception to this is the first one which is a simple identity relating a determinant to the fourth power of a Vandermonde, whose proof is given because it is a little less common.
\begin{thm}\label{vandthm}
  \begin{equation*}
    \prod_{0\leq j<k\leq N}(x_{j}-x_{k})^{4}=\det\left(x_{k}^{j}\quad j\,x_{k}^{j-1}\right)_{\substack{j=0,\dots,2\,N-1 \\ k=1,\ldots,N}}
  \end{equation*}
\end{thm}
\begin{proof}
  First note that
  \begin{equation*}
    \det\left(x_{k}^{j}\quad j\,x_{k}^{j-1}\right)_{\substack{j=0,\dots,2\,N-1 \\ k=1,\ldots,N}}=\left. \frac{\partial}{\partial y_{1}\ldots\partial y_{N}}\,\det\left(x_{k}^{j}\quad y_{k}^{j}\right)_{\substack{j=0,\dots,2\,N-1 \\ k=1,\ldots,N}}\right|_{y_{k}=x_{k}}.
  \end{equation*}
  Now
    \begin{align*}
      \det \left( x_{k}^{j}\quad y_{k}^{j}\right)&_{\substack{j=0,\dots,2\,N-1 \\
          k=1,\ldots,N}} =   \left(
        \begin{array}{ccccc} 1 & 1 & 1 & 1 & \cdots \\
          x_{1} & y_{1} & x_{2} & y_{2} & \cdots \\
          x_{1}^{2} & y_{1}^{2} & x_{2}^{2} & y_{2}^{2} & \cdots \\
          \vdots & \vdots & \vdots & \vdots & \\
        \end{array}
      \right) \\
      &  = \prod_{j>k}\,(x_{j} - x_{k})\,(y_{j} - y_{k})\,(y_{j} - x_{k})\,(x_{j} - y_{k})\,\prod\,(y_{j} - x_{j}).
    \end{align*}
  Differentiate this product with respect to each $y_{i}$ and applying the product rule, we obtain
  \begin{eqnarray*}
    \begin{aligned}
      & \left[\frac{\partial}{\partial y_{1}\ldots\partial y_{N}} \prod(y_{j} - y_{k})\right]\,\prod_{j>k}\,(x_{j} - x_{k})\,(y_{j} - x_{k})\,(x_{j} - y_{k})\,\prod\,(y_{j} - x_{j})\\
      & \qquad + \left[\frac{\partial}{\partial y_{1}\ldots\partial y_{N}}\prod(y_{j} - x_{k})\right]\,\prod_{j>k}\,(x_{j} - x_{k})\,(y_{j} - y_{k})\,(x_{j} - y_{k})\,\prod\,(y_{j} - x_{j})\\
      & \qquad  + \left[\frac{\partial}{\partial y_{1}\ldots\partial y_{N}}\prod(x_{j} - y_{k})\right]\,\prod_{j>k}\,(x_{j} - x_{k})\,(y_{j} - y_{k})\,(y_{j} - x_{k})\,\prod\,(y_{j} - x_{j})\\
      & \qquad + \left[\frac{\partial}{\partial y_{1}\ldots\partial y_{N}}\prod(y_{j} - x_{k})\right]\,\prod_{j>k}\,(x_{j} - x_{k})\,(y_{j} - y_{k})\,(y_{j} - x_{k})\,(x_{j} - y_{k}).
    \end{aligned}
  \end{eqnarray*}
  We set $y_{k}=x_{k}$ after differentiating. In each of the first three terms the last factor will be $0$. Therefore, only the last term will remain with $y_{k}=x_{k}$.  Now
  \begin{equation*}
    \left[\frac{\partial}{\partial y_{1}\ldots\partial y_{N}}\prod(y_{j} - x_{k})\right]=1.
  \end{equation*}
  So we are left with
  \begin{equation*}
    \left. \prod_{j>k}\,(x_{j} - x_{k})\,(y_{j} - y_{k})\,(y_{j} - x_{k})\,(x_{j} - y_{k})\right|_{y_{k}=x_{k}}=\prod_{0\leq j<k\leq N}(x_{j}-x_{k})^{4}.
  \end{equation*}
\end{proof}
\begin{thm}\label{detthm}
  If $A, B$ are Hilbert--Schmidt operators on a general \footnote{See \cite{Gohb1} for proof.} Hilbert space $\mathcal H$, then
  \begin{equation*}
    \det(I + A\,B)=\det(I + B\,A).
  \end{equation*}
\end{thm}
\begin{thm}[de Bruijn, 1955]
  \begin{eqnarray}\label{debruijn1}
    \begin{aligned}
      \int\cdots\int\det(\varphi_{j}(x_{k}))_{_{1\leq j,k\leq N}}\cdot &\det(\psi_{j}(x_{k}))_{_{1\leq j,k\leq N}}\,d\,\mu(x_{1})\cdots d\,\mu(x_{N}) \\
      & = N!\,\det\left(\int\varphi_{j}(x)\,\psi_{k}(x)\,d\,\mu(x)\right)_{_{1\leq j,k\leq N}},
    \end{aligned}
  \end{eqnarray}
  \begin{eqnarray}\label{debruijn2}
    \begin{aligned}
      \underset{x_{1}\leq\ldots\leq x_{N}}{\int\cdots\int} & \det(\varphi_{j}(x_{k}))_{_{1\leq j,k\leq N}}  d\,\mu(x_{1})\cdots d\,\mu(x_{N}) \\
      & = \pfaffian\left(\int\int\sgn(x-y)\varphi_{j}(x)\,\varphi_{k}(x)\,d\,\mu(x)\,d\,\mu(y)\right)_{_{1\leq j,k\leq N}},
    \end{aligned}
  \end{eqnarray}
  \begin{eqnarray}\label{debruijn3}
    \begin{aligned}
      \int\cdots\int & \det(\varphi_{j}(x_{k})\quad \psi_{j}(x_{k}))_{_{\substack{1\leq j\leq 2\,N\\
            1\leq  k \leq N}}} d\,\mu(x_{1})\cdots d\,\mu(x_{N}) \\
      & = (2\,N)!\,\pfaffian\left(\int \varphi_{j}(x)\psi_{k}(x)-\varphi_{k}(x)\psi_{j}(x)\,d\,\mu(x)\right)_{_{1\leq j,k\leq 2\,N}},
    \end{aligned}
  \end{eqnarray}
\end{thm}
where $\pfaffian$ denotes the Pfaffian. The last two integral identities were discovered by de Bruijn \cite{Debr1} in an attempt to generalize the first one. The first and last are valid in general measure spaces. In the second identity, the space needs to be ordered. In the last identity, the left hand side determinant is a $2\,N\times 2\,N$ determinant whose columns are alternating columns of the $\varphi_{j}$ and $\psi_{j}$ (i.e. the first four columns are $\{\varphi_{j}(x_{1})\}$, $\{\psi_{j}(x_{1})\}$, $\{\varphi_{j}(x_{2})\}$, $\{\psi_{j}(x_{2})\}$, respectively for $j=1,\ldots,2\,N$), hence the notation, and asymmetry in indexing.

\section{Recursion formula for the eigenvalue distributions}

With the joint density function defined as in \eqref{jointdensity}, let $J$ denote the interval $(t,\infty)$, and $\rchi=\rchi_{_{J}}(x)$ its characteristic function.\footnote{Much of what is said here is still valid if $J$ is taken to be a finite union of open intervals in $\mathbb R$ (see \cite{Trac7}). However, since we will only be interested in edge eigenvalues we restrict ourselves to $(t,\infty)$ from here on.} We denote by $\tilde{\rchi}=1-\rchi$ the characteristic function of the complement of $J$, and define $\tilde{\rchi}_{_{\lambda}}=1-\lambda\,\rchi$. Furthermore, let $E_{\beta,N}(t,m)$ equal the probability that exactly the $m$ largest eigenvalues of a matrix chosen at random from a (finite $N$) $\beta$--ensemble lie in $J$. We also define
\begin{equation}
\label{eq:23}
  G_{\beta,N}(t,\lambda)= \underset{x_{i}\in \mathbb R}{\int\cdots\int} \tilde{\rchi}_{_{\lambda}}(x_{1})\cdots \tilde{\rchi}_{_{\lambda}}(x_{N})\,P_{\beta}(x_{1},\ldots,x_{N})\,d\,x_{1}\cdots d\,x_{N}.
\end{equation}
For $\lambda=1$ this is just $E_{\beta,N}(t,0)$, the probability that no eigenvalues lie in $(t,\infty)$, or equivalently the probability that the largest eigenvalue is less than $t$. In fact we will see in the following propositions that $G_{\beta,N}(t,\lambda)$ is in some sense a generating function for $E_{\beta,N}(t,m)$.
\begin{prop}
  \begin{equation}
    G_{\beta,N}(t,\lambda)=\sum_{k=0}^{N}(-\lambda)^{k}\binom{N}{k}\underset{x_{i}\in J}{\int\cdots\int} P_{\beta}(x_{1},\ldots,x_{N})\,d\,x_{1}\cdots d\,x_{N}.
  \end{equation}
\end{prop}
\begin{proof}
  Using the definition of the $\tilde{\rchi}_{_{\lambda}}(x_{1})$ and multiplying out the integrand of \eqref{eq:23} gives
  \begin{equation*}
    G_{\beta,N}(t,\lambda) =  \sum_{k=0}^{N}(-\lambda)^{k}\underset{x_{i}\in \mathbb R}{\int\cdots\int} e_{k}(\rchi(x_{1}),\ldots, \rchi(x_{N}))P_{\beta}(x_{1},\ldots,x_{N})\,d\,x_{1}\cdots d\,x_{N},
  \end{equation*}
  where, in the notation of \cite{Stan2}, $e_{k}=m_{1^{k}}$ is the $k^{th}$ elementary symmetric function. Indeed each term in the summation arises from picking $k$ of the $\lambda\,\rchi$-terms, each of which comes with a negative sign, and $N-k$ of the $1$'s. This explains the coefficient $(-\lambda)^{k}$. Moreover, it follows that $e_{k}$ contains $\binom{N}{k}$ terms. Now the integrand is symmetric under permutations of the $x_{i}$. Also if $x_{i}\not\in J$, all corresponding terms in the symmetric function are $0$, and they are $1$ otherwise. Therefore we can restrict the integration to $x_{i}\in J$, remove the characteristic functions (hence the symmetric function), and introduce the binomial coefficient to account for the identical terms up to permutation.
\end{proof}
\begin{prop}
  \begin{equation}
    \left. E_{\beta,N}(t,m)=\frac{(-1)^{m}}{m\,!}\,\frac{d^{m}}{d\,\lambda^{m}}\,G_{\beta,N}(t,\lambda)\right|_{\lambda=1}\,,\quad m\geq 0.
  \end{equation}
\end{prop}
\begin{proof}
  This is proved by induction. As noted above, $E_{\beta,N}(t,0)=G_{\beta,N}(t,1)$ so it holds for the degenerate case $m=0$. When $m=1$ we have
  \begin{equation*}
    \begin{aligned}
      \left.-\frac{d}{d\,\lambda}\,G_{\beta,N}(t,\lambda)\right|_{\lambda=1} &=\left.-\frac{d}{d\,\lambda}\, \int\cdots\int \tilde{\rchi}_{_{\lambda}}(x_{1})\cdots \tilde{\rchi}_{_{\lambda}}(x_{n})\,P_{\beta}^{(N)}(\vec{x})\,d\,x_{1}\cdots d\,x_{n}\right|_{\lambda=1} \\
      &= -\sum_{j=1}^{N}- \int\cdots\int \tilde{\rchi}(x_{1})\cdots\tilde{\rchi}(x_{j-1})\,\rchi(x_{j})\,\tilde{\rchi}(x_{j+1})\cdots\\
      &   \qquad\qquad \cdots\tilde{\rchi}(x_{N})\, P_{\beta}^{(N)}(\vec{x})\,d\,x_{1}\cdots d\,x_{N}.
    \end{aligned}
  \end{equation*}
  The integrand is symmetric under permutations so we can make all terms look the same. There are $N=\binom{N}{1}$ of them so we get
  \begin{equation*}
    \begin{aligned}
      \left.-\frac{d}{d\,\lambda}\,G_{\beta,N}(t,\lambda)\right|_{\lambda=1} &=\left.\binom{N}{1}\,\int\cdots\int \rchi(x_{1})\,\tilde{\rchi}(x_{2})\cdots \tilde{\rchi}(x_{N})\,P_{\beta}^{(N)}(\vec{x})\,d\,x_{1}\cdots d\,x_{N}\right|_{\lambda=1} \\
      &=\binom{N}{1}\,\int\cdots\int \rchi(x_{1})\,\rchi(x_{2})\cdots \rchi(x_{N})\,P_{\beta}^{(N)}(\vec{x})\,d\,x_{1}\cdots d\,x_{N}\\
      &= E_{\beta,N}(t,1).
    \end{aligned}
  \end{equation*}
  When $m=2$ then
  \begin{equation*}
    \begin{aligned}
      \frac{1}{2}&\left(-\frac{d}{d\,\lambda}\right)^{2} \, \left. G_{\beta,N}(t,\lambda)\right|_{\lambda=1} \\
      & = \frac{N}{2}\,\sum_{j=2}^{N} \int\cdots\int \rchi(x_{1})\tilde{\rchi}(x_{2})\cdots\tilde{\rchi}(x_{j-1})\,\rchi(x_{j})\,\tilde{\rchi}(x_{j+1})\cdots \\
      &  \qquad\qquad \cdots\left.\tilde{\rchi}(x_{N})\, P_{\beta}^{(N)}(\vec{x})\,d\,x_{1}\cdots d\,x_{N}\right|_{\lambda=1}  \\
      &=\left.\frac{N\,(N-1)}{2}\,\int\cdots\int \rchi(x_{1})\,\rchi(x_{2})\,\tilde{\rchi}(x_{3})\cdots \tilde{\rchi}(x_{N})\,P_{\beta}^{(N)}(\vec{x})\,d\,x_{1}\cdots d\,x_{N}\right|_{\lambda=1} \\
      & =\left. \binom{N}{2}\int\cdots\int \rchi(x_{1})\,\rchi(x_{2})\,\tilde{\rchi}(x_{3})\cdots \tilde{\rchi}(x_{N})\,P_{\beta}^{(N)}(\vec{x})\,d\,x_{1}\cdots d\,x_{N}\right|_{\lambda=1} \\
      & =\binom{N}{2}\int\cdots\int \rchi(x_{1})\,\rchi(x_{2})\,\rchi(x_{3})\cdots \rchi(x_{N})\,P_{\beta}^{(N)}(\vec{x})\,d\,x_{1}\cdots d\,x_{N} \\
      & = E_{\beta,N}(t,2),
    \end{aligned}
  \end{equation*}
  where we used the previous case to get the first equality, and again the invariance of the integrand under symmetry to get the second equality. By induction then,
  \begin{eqnarray*}
    \lefteqn{\frac{1}{m\,!}\left(-\frac{d}{d\,\lambda}\right)^{m}\, \left. G_{\beta,N}(t,\lambda)\right|_{\lambda=1} =}\\
    & & = \frac{N\,(N-1)\cdots(N-m+2)}{m\,!}\,\sum_{j=m}^{N}\int\cdots\int \rchi(x_{1})\,\tilde{\rchi}(x_{2})\cdots \\
    & & \qquad \cdots\tilde{\rchi}(x_{j-1})\,\rchi(x_{j})\,\tilde{\rchi}(x_{j+1})\cdots\tilde{\rchi}(x_{N})\,P_{\beta}^{(N)}(\vec{x})\,d\,x_{1}\cdots d\,x_{N}  \\
    & & =\frac{N\,(N-1)\cdots(N-m+1)}{m\,!}\,\int\cdots\int\rchi(x_{1})\cdots \\
    & & \qquad \cdots\rchi(x_{m})\,\tilde{\rchi}(x_{m+1})\cdots\tilde{\rchi}(x_{N})\,P_{\beta}^{(N)}(\vec{x})\,d\,x_{1}\cdots d\,x_{N} \\
    & & = \binom{N}{m}\,\int\cdots\int\rchi(x_{1})\cdots\rchi(x_{m})\,\tilde{\rchi}(x_{m+1})\cdots\tilde{\rchi}(x_{N})\,P_{\beta}^{(N)}(\vec{x})\,d\,x_{1}\cdots d\,x_{N} \\
    & & = E_{\beta,N}(t,m). \end{eqnarray*}
\end{proof}
If we define $F_{\beta,N}(t,m)$ to be the distribution of the $m^{th}$ largest eigenvalue in the (finite $N$) $\beta$--ensemble, then the following probabilistic result is immediate from our definition of $E_{\beta,N}(t,m)$.
\begin{cor}
\begin{equation}F_{\beta,N}(t,m+1) -  F_{\beta,N}(t,m) = E_{\beta,N}(t,m)\end{equation}
\end{cor}

\chapter{The distribution of the $m^{th}$ largest eigenvalue in the GUE}

\chaptermark{The $m^{th}$ largest eigenvalue in the GUE}

\section{The distribution function as a Fredholm determinant}

We closely follow \cite{Trac1} for the derivations which follow. The GUE case corresponds to the specialization $\beta=2$ in \eqref{jointdensity} so that
\begin{equation}
  \label{eq:1}
  G_{2,N}(t,\lambda)=C_{2}^{(N)}\underset{x_{i}\in\mathbb{R}}{\int\cdots\int} \prod_{j<k}\left(x_{j}-x_{k}\right)^{2}\,\prod_{j}^{N}w(x_{j})\,\prod_{j}^{N}\left(1+f(x_{j})\right)\,d\,x_{1}\cdots d\,x_{N}
\end{equation}
where $w(x)=\exp\left(-x^{2}\right)$, $f(x)=-\lambda\,\rchi_{J}(x)$, and $C_{2}^{(N)}$ depends only on $N$. In the steps that follow, additional constants depending solely on $N$ (such as $N!$) which appear will be lumped into $C_{2}^{(N)}$. A probability argument will show that the resulting constant at the end of all calculations simply equals $1$. Expressing the Vandermonde as a determinant
\begin{equation}
  \prod_{1\leq j<k\leq N}(x_{j}-x_{k})=\det\left(x_{k}^{j}\right)_{\substack{j=0,\ldots,N\\k=1,\ldots,N}}
\end{equation}
and using \eqref{debruijn1} with $\varphi_{j}(x)=\psi_{j}(x)=x^{j}$ and $d\,\mu(x)=w(x)(1+f(x))$ yields
\begin{equation}
  G_{2,N}(t,\lambda)=C_{2}^{(N)}\det\left(\int_{\mathbb{R}} x^{j+k}\,w(x)\left(1+f(x)\right)\,d\,x)\right)_{j,k=0,\ldots,N-1}.
\end{equation}
Let $\left\{\varphi_{j}(x)\right\}$ be the sequence obtained by orthonormalizing the sequence $\left\{x^{j}\,w^{1/2}(x)\right\}$. It follows that
\begin{eqnarray}
  G_{2,N}(t,\lambda) & = & C_{2}^{(N)}\det\left(\int_{\mathbb{R}} \varphi_{j}(x)\,\varphi_{k}(x)\,\left(1+f(x)\right)\,d\,x)\right)_{j,k=0,\ldots,N-1}\\
  & = & C_{2}^{(N)}\det\left(\delta_{j,k}+\int_{\mathbb{R}} \varphi_{j}(x)\,\varphi_{k}(x)\,f(x)\,d\,x)\right)_{j,k=0,\ldots,N-1}.
\end{eqnarray}
The last expression is of the form $\det(I+AB)$ for $A:L^{2}(\mathbb{R})\to \mathbb{C}^{N}$ with kernel $A(j,x)=~\varphi_{j}(x)f(x)$ whereas $B:\mathbb{C}^{N}\to L^{2}(\mathbb{R})$ with kernel $B(x,j)=\varphi_{j}(x)$. Note that $AB:\mathbb{C}^{N}\to \mathbb{C}^{N}$ has kernel
\begin{equation}
AB(j,k)=\int_{\mathbb{R}} \varphi_{j}(x)\,\varphi_{k}(x)\,f(x)\,d\,x
\end{equation}
whereas $BA:L^{2}(\mathbb{R})\to L^{2}(\mathbb{R})$ has kernel
\begin{equation}
BA(x,y)=\sum_{k=0}^{N-1}\varphi_{k}(x)\,\varphi_{k}(y) :=K_{2,N}(x,y).
\end{equation}
From \ref{detthm} it follows that
\begin{equation}
  G_{2,N}(t,\lambda)=C_{2}^{(N)}\det\left(I-K_{2,N}\,f\right),
\end{equation}
where $K_{2,N}$ has kernel $K_{2,N}(x,y)$ and $K_{2,N}\,f$ acts on a function by first multiplying it by $f$ and acting on the product with $K_{2,N}$. From \eqref{eq:1} we see that setting $f=0$ in the last identity yields $C_{2}^{(N)}=1$. Thus the above simplifies to
\begin{equation}
  \label{eq:22}
  G_{2,N}(t,\lambda)=\det\left(I-K_{2,N}\,f\right).
\end{equation}

\section{Edge scaling and differential equations}

Following the derivation in \cite{Trac4}, we specialize $w(x)=\exp\left(-x^{2}\right)$, $f(x)=-\lambda\,\rchi_{J}(x)$, so that the $\left\{\varphi_{j}(x)\right\}$ are in fact the Hermite polynomials times the square root of the weight. Using the asymptotics of Hermite polynomials, it follows that in the so-called \emph{edge scaling limit},
\begin{equation}
  \label{eq:edgescaling}
\lim_{N\to\infty}\frac{1}{2^{1/2}N^{1/6}}\,K_{N,2}\left(\sqrt{2N}+\frac{x}{2^{1/2}N^{1/6}},\sqrt{2N}+\frac{y}{2^{1/2}N^{1/6}}\right)\,\rchi_{J}\left(\sqrt{2N}+\frac{
y}{2^{1/2}N^{1/6}}\right)
\end{equation}
is $K_{\airy}(x,y)$ as defined in \eqref{kairy}. The fact that this limit holds in trace class norm is crucial  (see \cite{Bowi1}, \cite{Forr2} for proofs, and \cite{Trac5} for corresponding results in the $\beta=1,4$ cases). It is what allows us to take the limit inside the determinant in \eqref{eq:22}. For notational convenience, we denote the corresponding operator $K_{\airy}$ by $K$ in the rest of this subsection. We also think of $K$ as the integral operator with kernel
\begin{equation}
  \label{eq:15}
  K(x,y)=\frac{\varphi(x)\psi(y)-\psi(x)\varphi(y)}{x-y}\,\rchi_{J}(y),
\end{equation}
where $\varphi(x)=\sqrt{\lambda}\airy(x)$, $\psi(x)=\sqrt{\lambda}\airy^{\prime}(x)$ and $J$ is $\left(s,\infty\right)$ with
\begin{equation}
  t=\sqrt{2N}+ \frac{s}{\sqrt{2}N^{1/6}}.
\end{equation}
Note that although $K(x,y)$, $\varphi$ and $\psi$ are functions of $\lambda$ as well, this dependence will not affect our calculations in what follows. Thus we omit it to avoid cumbersome notation. The Airy equation implies that $\varphi$ and $\psi$ satisfy the relations
\begin{eqnarray}
  \frac{d}{d\,x}\,\varphi &= &\psi,\nonumber\\
  \frac{d}{d\,x}\,\psi &=& x\,\varphi.
\end{eqnarray}
We define $D_{2,N}(s,\lambda)$ to be the Fredholm determinant $\det\left(I-K\right)$. Thus in the edge scaling limit
\begin{equation*}
  E_{2}(s,\lambda):=E_{2,\infty}(s,\lambda)=D_{2}(s,\lambda).
\end{equation*}
We define the operator
\begin{equation}
  \label{eq:rdef}
  R=(I-K)^{-1}K,
\end{equation}
whose kernel we denote $R(x,y)$. Incidentally, we shall use the notation $\doteq$ in reference to an operator to mean ``has kernel''. For example $R\doteq R(x,y)$. We also let $M$ stand for the operator whose action is multiplication by $x$. It is well known that
\begin{equation}
  \label{rderiv}
  \frac{d}{ds}\log\det\left(I-K\right)=-R(s,s).
\end{equation}
For functions $f$ and $g$, we write $f\otimes g$ to denote the operator specified by
\begin{equation}
  f\otimes g \doteq f(x)g(y),
\end{equation}
and define
\begin{eqnarray}
  Q(x,s)&=& Q(x) = \left(\left(I-K\right)^{-1}\varphi\right)(x),\\
  P(x,s) &=& P(x) = \left(\left(I-K\right)^{-1}\psi\right)(x).
\end{eqnarray}
Then straightforward computation yields the following facts
\begin{eqnarray}
  \com{M}{K} &=& \varphi\otimes\psi - \psi\otimes\varphi,\nonumber\\
  \com{M}{\left(I-K\right)^{-1}} &=& \left(I-K\right)^{-1}\com{M}{K}\left(I-K\right)^{-1} \nonumber\\
  & = & Q\otimes P-P\otimes Q.
\end{eqnarray}
On the other hand if $\left(I-K\right)^{-1}\doteq\rho(x,y)$, then
\begin{equation}
  \rho(x,y)=\delta(x-y)+R(x,y),
\end{equation}
and it follows that
\begin{equation}
\com{M}{\left(I-K\right)^{-1}} \doteq (x-y)\rho(x,y)=(x-y)R(x,y).
\end{equation}
Equating the two representation for the kernel of $\com{M}{\left(I-K\right)^{-1}}$ yields
\begin{equation}
  R(x,y)=\frac{Q(x)P(y)-P(x)Q(y)}{x-y}.
\end{equation}
Taking the limit $y\to x$ and defining $q(s)=Q(s,s)$, $p(s)=P(s,s)$, we obtain
\begin{equation}
  \label{RDiag}
  R(s,s)=Q^{\prime}(s,s)\,p(s)-P^{\prime}(s,s)\,q(s).
\end{equation}
Let us now derive expressions for $Q^{\prime}(x)$ and $P^\prime(x)$. If we let the operator $D$ stand for differentiation with respect to $x$,
\begin{eqnarray}
  Q^\prime(x,s)&=& D \left(I-K\right)^{-1} \varphi \nonumber\\
  &=& \left(I-K\right)^{-1} D\varphi +
  \left[D,\left(I-K\right)^{-1}\right]\varphi\nonumber\\
  &=& \left(I-K\right)^{-1} \psi +
  \left[D,\left(I-K\right)^{-1}\right]\varphi\nonumber\\
  &=& P(x) + \left[D,\left(I-K\right)^{-1}\right]\varphi. \label{Qderiv1}
\end{eqnarray}
We need the commutator
\begin{equation}
  \left[D,\left(I-K\right)^{-1}\right]=\left(I-K\right)^{-1} \left[D,K\right] \left(I-K\right)^{-1}.
\end{equation}
Integration by parts shows
\begin{equation}
  \left[D,K\right] \doteq \left( \frac{\partial K}{\partial x} + \frac{\partial K}{\partial y}\right) + K(x,s) \delta(y-s).
\end{equation}
The $\delta$ function comes from differentiating the characteristic function $\chi$. Moreover,
\begin{equation}
  \left( \frac{\partial K}{\partial x} + \frac{\partial K}{\partial y}\right) = \varphi(x) \varphi(y).
\end{equation}
Thus
\begin{equation}
 \label{DComm}
\left[D,\left(I-K\right)^{-1}\right]\doteq - Q(x) Q(y) + R(x,s) \rho(s,y).
\end{equation}
(Recall $(I-K)^{-1}\doteq \rho(x,y)$.)  We now use this in (\ref{Qderiv1}) to obtain
\begin{eqnarray*}
  Q^\prime(x,s)&=&P(x) - Q(x) \left(Q,\varphi\right) + R(x,s) q(s) \\
  &=& P(x) - Q(x) u(s) + R(x,s) q(s),
\end{eqnarray*}
where the inner product $\left(Q,\varphi\right)$ is denoted by $u(s)$. Evaluating  at $x=s$  gives
\begin{equation}
  \label{Qderiv2}
  Q^\prime(s,s) = p(s) - q(s) u(s) +R(s,s) q(s).
\end{equation}
We now apply the same procedure to compute $P^\prime$.
\begin{eqnarray*}
  P^\prime(x,s)&=& D \left(I-K\right)^{-1} \psi \\
  &=& \left(I-K\right)^{-1} D\psi + \left[D,\left(I-K\right)^{-1}\right]\psi\\
  &=& M \left(I-K\right)^{-1} \varphi +  \left[\left(I-K\right)^{-1},M\right]\varphi+   \left[D,\left(I-K\right)^{-1}\right]\psi\\
  &=& x Q(x) +\left(P\otimes Q-Q\otimes P\right)\varphi +(-Q\otimes Q)\psi +  R(x,s) p(s)\\
  &=& x Q(x) + P(x)\left(Q,\varphi\right) -  Q(x) \left(P,\varphi\right)  - Q(x) \left(Q,\psi\right)+R(x,s)p(s)\\
  &=& x Q(x) - 2 Q(x) v(s) + P(x) u(s) + R(x,s) p(s).
\end{eqnarray*}
Here $v=\left(P,\varphi\right)=\left(\psi,Q\right)$. Setting $x=s$ we obtain
\begin{equation}
  \label{Pderiv}
  P^{\prime}(s,s) = s q(s) + 2 q(s) v(s) +p(s) u(s) +R(s,s) p(s).
\end{equation}
Using this and the expression for $Q^\prime(s,s)$ in (\ref{RDiag}) gives
\begin{equation}
  \label{RDiag2}
  R(s,s)= p^2-s q^2 + 2 q^2 v - 2 p q u.
\end{equation}
Using the chain rule, we have
\begin{equation}
  \label{qDeriv}
  \frac{d\,q}{d\,s} = \left( \frac{\partial}{\partial x} + \frac{\partial}{\partial s}\right) Q(x,s)\left\vert_{x=s}. \right.
\end{equation}
The first term is known. The partial with respect to $s$ is
\begin{eqnarray*}
  \frac{\partial Q(x,s)}{\partial s}&=& \left(I-K\right)^{-1} \frac{\partial K}{\partial s} \left(I-K\right)^{-1} \varphi\\
  &=& - R(x,s) q(s),
\end{eqnarray*}
where we used the fact that
\begin{equation}
  \frac{\partial K}{\partial s}\doteq -K(x,s)\delta(y-s).
\end{equation}
Adding the two partial derivatives  and evaluating at $x=s$ gives
\begin{equation}
  \label{qEqn}
  \frac{d\,q}{d\,s} = p - q u.
\end{equation}
A similar calculation gives
\begin{equation}
  \label{pEqn}
  \frac{d\,p}{d\,s}= s q - 2 q v + p u.
\end{equation}
We derive first order differential equations for  $u$ and $v$  by differentiating the inner products. Recall that
\begin{equation*}
  u(s) = \int_s^\infty \varphi(x) Q(x,s)\, d\,x.
\end{equation*}
Thus
\begin{eqnarray*}
  \frac{d\,u}{ds}&=& -\varphi(s) q(s) + \int_s^\infty \varphi(x) \frac{\partial Q(x,s)}{\partial s}\, d\,x \\
  &=& -\left(\varphi(s)+\int_s^\infty R(s,x) \varphi(x)\,d\,x\right) q(s)\\
  &=& -\left(I-K\right)^{-1} \varphi(s) \, q(s)\\
  &=& - q^2.
\end{eqnarray*}
Similarly,
\begin{equation}
  \frac{d\,v}{d\,s} = - p q.
\end{equation}
From the first order differential equations for $q$, $u$ and $v$ it follows immediately  that the derivative of   $ u^2-2v-q^2 $ is zero.  Examining the behavior near $s=\infty$ to check that the constant of integration is zero then gives
\begin{equation}
  u^2-2v=q^2.
\end{equation}
We now differentiate (\ref{qEqn}) with respect to $s$, use the first order differential equations for $p$ and $u$, and then the  first integral to deduce that $q$ satisfies the \textit{Painlev\'e II equation}
\begin{equation}
  \label{P2}
  q^{\prime\prime}=s q + 2 q^3.
\end{equation}
Checking the asymptotics of the Fredholm determinant $\det(I-K)$ for large $s$ shows we want the solution to the Painlev\'e II equation with boundary condition
\begin{equation}
  \label{bc}
  q(s,\lambda)\sim \sqrt{\lambda}\airy(s) \qquad \textrm{as} \qquad s\rightarrow\infty.
\end{equation}
That a solution $q$ exists and is unique follows from the representation of the Fredholm determinant in terms of it.  Independent proofs of this, as well as the asymptotics as $s\ra\-\iy$ were given by \cite{Hasti1}, \cite{Clar1}, \cite{Deif2}. Since $\com{D}{(I-K)^{-1}}\doteq (\partial/\partial x~+~\partial/\partial y)~R(x,y)$, (\ref{DComm}) says
\begin{equation}
  \left(\frac{\partial}{\partial x}+\frac{\partial}{\partial y}\right)R(x,y)=-Q(x)Q(y)+R(x,s)\rho(s,y).
\end{equation}
In computing $\partial Q(x,s)/\partial s$ we showed that
\begin{equation}
  \frac{\partial}{\partial s} \left(I-K\right)^{-1}\doteq \frac{\partial}{\partial s}R(x,y) = -R(x,s)\rho(s,y).
\end{equation}
Adding these two expressions,
\begin{equation}
  \left(\frac{\partial}{\partial x} + \frac{\partial}{\partial y}+ \frac{\partial}{\partial s}\right)R(x,y) = -Q(x)\,Q(y),
\end{equation}
and then evaluating at $x=y=s$ gives
\begin{equation}
  \frac{d}{d\,s}R(s,s)=-q^2. \label{Rderiv}
\end{equation}
Integration  (and recalling (\ref{rderiv})) gives,
\begin{equation}
  \frac{d}{d\,s}\log\det\left(I-K\right)=-\int_s^\infty q^2(x,\lambda) \, d\,x;
\end{equation}
and hence,
\begin{equation}
  \log\det\left(I-K\right)=-\int_s^\infty\left(\int_y^\infty q^2(x,\lambda)\,d\,x\right)\, d\,y =-\int_s^\infty (x-s) q^2(x,\lambda)\, d\,x.
\end{equation}
To summarize,  we have shown that $D_{2}(s,\lambda)$ has the Painlev\'e representation
\begin{equation}
  \label{F2}
  D_{2}(s,\lambda) = \exp\left(-\int_s^\infty (x-s) q^2(x,\lambda)\, d\,x\right)
\end{equation}
where $q$  satisfies the Painlev\'e II equation (\ref{P2}) subject to the boundary condition (\ref{bc}).

\chapter{The distribution of the $m^{th}$ largest eigenvalue in the GSE}

\chaptermark{The $m^{th}$ largest eigenvalue in the GSE}

\section{The distribution function as a Fredholm determinant}

The GSE corresponds case corresponds to the specialization $\beta=4$ in \eqref{jointdensity} so that
\begin{equation}
  G_{4,N}(t,\lambda)=C_{4}^{(N)}\underset{x_{i}\in\mathbb{R}}{\int\cdots\int} \prod_{j<k}\left(x_{j}-x_{k}\right)^{4}\,\prod_{j}^{N}w(x_{j})\,\prod_{j}^{N}\left(1+f(x_{j})\right)\,d\,x_{1}\cdots d\,x_{N}
\end{equation}
where $w(x)=\exp\left(-x^{2}\right)$, $f(x)=-\lambda\,\rchi_{J}(x)$, and $C_{4}^{(N)}$ depends only on $N$. As in the GUE case, we will lump into $C_{4}^{(N)}$ any constants depending only on $N$ that appear in the derivation. A simple argument  at the end will show that the final constant is $1$. These calculations more or less faithfully follow and expand on \cite{Trac1}. By \eqref{vandthm}, $G_{4,N}(t,\lambda)$ is given by the integral
\begin{equation*}
  C_{4}^{(N)}\underset{x_{i}\in \mathbb R}{\int\cdots\int} \det\left(x_{k}^{j}\quad j\,x_{k}^{j-1}\right)_{\substack{j=0,\dots,2\,N-1 \\ k=1,\ldots,N}}\,\prod_{i=1}^{N}w(x_{i})\prod_{i=1}^{N}(1+f(x_{i}))\,d\,x_{1}\cdots d\,x_{N}
\end{equation*}
which, if we define $\varphi_{j}(x)=x^{j-1}\,w(x)\,(1+f(x))$ and $\psi_{j}(x)=(j-1)\,x^{j-2}$ and use the linearity of the determinant, becomes
\begin{equation*}
  G_{4,N}(t,\lambda)= C_{4}^{(N)}\underset{x_{i}\in \mathbb R}{\int\cdots\int} \det\left(\varphi_{j}(x_{k})\quad \psi(x_{k})\right)_{\substack{1\leq j \leq 2\,N \\
      1\leq k \leq N}}\,d\,x_{1}\cdots d\,x_{N}.
\end{equation*}
Now using \eqref{debruijn3}, we obtain
\begin{eqnarray*}
  G_{4,N}(t,\lambda) & = & C_{4}^{(N)}\,\pfaffian\left(\int \varphi_{j}(x)\psi_{k}{x}-\varphi_{k}(x)\psi_{j}(x)\,d\,x\right)_{_{1\leq j,k\leq 2\,N}}\\
  &=& C_{4}^{(N)}\,\pfaffian\left(\int (k-j)\,x^{j+k-3}\,w(x)\,(1+f(x))\,d\,x\right)_{_{1\leq j,k\leq 2\,N}} \\
  &=&  C_{4}^{(N)}\pfaffian\left(\int (k-j)\,x^{j+k-1}\,w(x)\,(1+f(x))\,d\,x\right)_{_{0\leq j,k\leq 2\,N-1}},
\end{eqnarray*}
where we let $k\to k+1$ and $j\to j+1$ in the last line. Remembering that the square of a Pfaffian is a determinant, we
obtain
\begin{equation*}
  G_{4,N}^{2}(t,\lambda)=C_{4}^{(N)}\, \det\left(\int (k-j)\,x^{j+k-1}\,w(x)\,(1+f(x))\,d\,x\right)_{_{0\leq j,k\leq 2\,N-1}}.
\end{equation*}
Row operations on the matrix do not change the determinant, so we can replace $\{x^{j}\}$ by an arbitrary sequence $\{p_{j}(x)\}$ of polynomials of degree $j$ obtained by adding rows to each other. Note that the general $(j,k)$ element in the matrix can be written as
\begin{equation*}
  \left[\left(\frac{d}{d\,x}\,x^{k}\right)\,x^{j} - \left(\frac{d}{d\,x}\,x^{j}\right)\,x^{k} \right]\,w(x) \,\left(1+f(x)\right).
\end{equation*}
Thus when we add rows to each other the polynomials we obtain will have the same general form (the derivatives factor). Therefore we can assume without loss of generality that $G_{4,N}^{2}(t,\lambda)$ equals
\begin{equation*}
  C_{4}^{(N)}\, \det\left(\int \left[ p_{j}(x)\,p_{k}^{\prime}(x) - \,p_{j}^{\prime}(x)\,p_{k}(x)\right]\,w(x)\,\left(1+f(x)\right)\,d\,x\right)_{_{0\leq j,k\leq 2\,N-1}},
\end{equation*}
where the sequence $\{p_{j}(x)\}$ of polynomials of degree $j$ is arbitrary. Let $\psi_{j}=p_{j}\,w^{1/2}$ so that $p_{j}=\psi_{j}\,w^{-1/2}$. Substituting this into the above formula and simplifying, we obtain
\begin{eqnarray*}
  G_{4,N}^{2}(t,\lambda)&=&C_{4}^{(N)}\, \det\left[ \int\left( \left( \psi_{j}(x)\,\psi_{k}^{\prime}(x) - \,\psi_{k}(x)\,\psi_{j}^{\prime}(x)\right)\left(1+f(x)\right)\right)\,d\,x \right]_{_{0\leq j,k\leq 2\,N-1}}\\
  &=&C_{4}^{(N)}\, \det\left[M + L\right] =C_{4}^{(N)}\, \det[M]\,\det[I+ M^{-1}\cdot L],
\end{eqnarray*}
where $M,L$ are matrices given by
\begin{equation*}
  M=\left(\int \left( \psi_{j}(x)\,\psi_{k}^{\prime}(x) - \,\psi_{k}(x)\,\psi_{j}^{\prime}(x)\right)\,d\,x\right)_{_{0\leq j,k\leq 2\,N-1}},
\end{equation*}
\begin{equation*}
  L=\left(\int \left( \psi_{j}(x)\,\psi_{k}^{\prime}(x) - \,\psi_{k}(x)\,\psi_{j}^{\prime}(x)\right)f(x)\,d\,x\right)_{_{0\leq j,k\leq 2\,N-1}}.
\end{equation*}
Note that $\det[M]$ is a constant which depends only on $N$ so we can absorb it into $C_{4}^{(N)}$. Also if we denote
\begin{equation*}
  M^{-1}=\left\{\mu_{jk}\right\}_{_{0\leq j,k\leq 2\,N-1}}\quad , \quad \eta_{j}=\sum_{k=0}^{2N-1}\mu_{jk}\psi_{k}(x),
\end{equation*}
it follows that
\begin{equation*}
  M^{-1}\cdot N=\left\{\int \left( \eta_{j}(x)\,\psi_{k}^{\prime}(x) - \,\eta_{j}^{\prime}(x)\,\psi_{k}(x)\right)f(x)\,d\,x\right\}_{_{0\leq j,k\leq 2\,N-1}}.
\end{equation*}
Let $A:L^{2}(\mathbb R)\times L^{2}(\mathbb R)\to \mathbb C^{2N}$ be the operator defined by the $2N\times 2$ matrix
\begin{equation*}
  A(x)=\left(
    \begin{array}{cc}
      \eta_{0}(x) & -\eta_{0}^{\prime}(x) \\
      \eta_{1}(x) & -\eta_{1}^{\prime}(x) \\
      \vdots & \vdots  \\
    \end{array}
  \right).
\end{equation*}
Thus if
\begin{equation*}
  g=
  \left(
    \begin{array}{c}
      g_{0}(x)  \\
      g_{1}(x)
    \end{array}
  \right)
  \in L^{2}(\mathbb R)\times L^{2}(\mathbb R),
\end{equation*}
we have
\begin{equation*}
  A\,g=A(x)\,g=
  \left(
    \begin{array}{c}
      \int\left(\eta_{0}g_{0}-\eta_{0}^{\prime}g_{1}\right)\,d\,x \\
      \int\left(\eta_{1}g_{0}-\eta_{1}^{\prime}g_{1}\right)\,d\,x \\
      \vdots  \\
    \end{array}
  \right)
  \in \mathbb C^{2N}.
\end{equation*}
Similarly we define $B:\mathbb C^{2n}\to L^{2}(\mathbb R)\times L^{2}(\mathbb R)$ given by the $2\times 2n$ matrix
\begin{equation*}
  B(x) = f\cdot
  \left(
    \begin{array}{ccc}
      \psi_{0}^{\prime}(x) & \psi_{1}^{\prime}(x) & \cdots \\
      \psi_{0}(x) & \psi_{1}(x) & \cdots
    \end{array}
  \right).
\end{equation*}
Explicitly if
\begin{equation*}
  \alpha=
  \left(
    \begin{array}{c}
      \alpha_{0}  \\
      \alpha_{1} \\
      \vdots\end{array}
  \right)
  \in \mathbb C^{2N},
\end{equation*}
then
\begin{equation*}
  B\alpha= B(x)\cdot\alpha=
  \left(
    \begin{array}{c}
      f\displaystyle{\sum_{i=0}^{2N-1}\alpha_{i}\psi_{i}^{\prime}}\\
      f\displaystyle{\sum_{i=0}^{2N-1}\alpha_{i}\psi_{i}}
    \end{array}
  \right)\in L^{2}\times L^{2}.
\end{equation*}
Observe that $M^{-1}\cdot L=A B: \mathbb C^{2n}\to \mathbb C^{2n}$. Indeed
\begin{eqnarray*}
  A B \alpha &=&
  \left(
    \begin{array}{c}
      \displaystyle{\sum_{i=0}^{2N-1}\left[\int\left( \eta_{0}\psi_{i}^{\prime}-\eta_{0}^{\prime}\psi_{i}\right)f\, d\,x\right] \alpha_{i}} \\
      \displaystyle{\sum_{i=0}^{2N-1}\left[\int\left( \eta_{1}\psi_{i}^{\prime}-\eta_{1}^{\prime}\psi_{i}\right)f\, d\,x\right] \alpha_{i}} \\
      \vdots
    \end{array}
  \right) \\
  & = & \left\{\int\left(\eta_{j}\psi_{k}^{\prime}-\eta_{j}^{\prime}\psi_{k}\right)f\, d\,x \right\}\cdot
  \left(
    \begin{array}{c}
      \alpha_{0}  \\
      \alpha_{1} \\
      \vdots
    \end{array}
  \right)
  =\left( M^{-1}\cdot L\right) \alpha.
\end{eqnarray*}
Therefore, by \eqref{detthm}
\begin{equation*}G_{4,N}^{2}(t,\lambda) = C_{4}^{(N)}\det(I+M^{-1}\cdot L) =C_{4}^{(N)}\det(I+AB)=C_{4}^{(N)}\det(I+BA)\end{equation*}
where $BA:L^{2}\left(\mathbb R\right)\to L^{2}\left(\mathbb R\right)$. From our definition of $A$ and $B$ it follows that
\begin{eqnarray*}
  B\,A\,g & = &
  \left(
    \begin{array}{c}
      f\displaystyle{\sum_{i=0}^{2n-1}\psi_{i}^{\prime} \left( \int\left(\eta_{i}g_{0}-\eta_{i}^{\prime}g_{1}\right)\,d\,x \right)} \\
      f \displaystyle{\sum_{i=0}^{2N-1}\psi_{i}^{\prime} \left(\int\left(\eta_{i}g_{0}-\eta_{i}^{\prime}g_{1}\right)\,d\,x \right)}
    \end{array}\right) \\
  & =&  f
  \left(
    \begin{array}{c}
      \displaystyle{\int \sum_{i=0}^{2N-1}\psi_{i}^{\prime}(x)\eta_{i}(y)g_{0}(y)\,d\,y - \int \sum_{i=0}^{2N-1}\psi_{i}^{\prime}(x)\eta_{i}^{\prime}(y)g_{1}(y)\,d\,y } \\
      \displaystyle{\int\sum_{i=0}^{2N-1}\psi_{i}(x)\eta_{i}(y)g_{0}(y)\,d\,y - \int\sum_{i=0}^{2N-1}\psi_{i}(x)\eta_{i}^{\prime}(y)g_{1}(y)\,d\,y}
    \end{array}
  \right)
  \\
  & = &f\, K_{4,N}\,g,
\end{eqnarray*}
where $K_{4,N}$ is the integral operator with matrix kernel
\begin{equation*}
  K_{4,N}(x,y)=
  \left(
    \begin{array}{cc}
      \displaystyle{ \sum_{i=0}^{2N-1}\psi_{i}^{\prime}(x)\eta_{i}(y)} & \displaystyle{-\sum_{i=0}^{2N-1}\psi_{i}^{\prime}(x)\eta_{i}^{\prime}(y)} \\
      \displaystyle{ \sum_{i=0}^{2N-1}\psi_{i}(x)\eta_{i}(y)} & \displaystyle{- \sum_{i=0}^{2N-1}\psi_{i}(x)\eta_{i}^{\prime}(y)}
    \end{array}
  \right).
\end{equation*}
Recall that
$\displaystyle{\eta_{j}(x)=\sum_{k=0}^{2N-1}\mu_{jk}\psi_{k}(x)}$ so that
\begin{equation*}
  K_{4,N}(x,y)=
  \left(
    \begin{array}{cc}
      \displaystyle{ \sum_{j,k=0}^{2N-1}\psi_{j}^{\prime}(x)\mu_{jk}\psi_{k}(y)} & \displaystyle{-\sum_{j,k=0}^{2N-1}\psi_{j}^{\prime}(x)\mu_{jk}\psi_{k}^{\prime}(y)} \\
      \displaystyle{ \sum_{j,k=0}^{2N-1}\psi_{j}(x)\mu_{jk}\psi_{k}(y)}& \displaystyle{-\sum_{j,k=0}^{2N-1}\psi_{j}(x)\mu_{jk}\psi_{k}^{\prime}(y)}
    \end{array}
  \right).
\end{equation*}
Define $\epsilon$ to be the following integral operator
\begin{equation}
  \label{epsilonop}
\epsilon\doteq  \epsilon(x-y)=
  \begin{cases}
    \frac{1}{2} &\textrm{if $x>y$}, \\
    -\frac{1}{2} &\textrm{if $x<y$}.
  \end{cases}
\end{equation}
As before, let $D$ denote the operator that acts by differentiation with respect to $x$. The fundamental theorem of calculus implies that $D\,\epsilon~=~\epsilon\,D~=~I$. We also define
\begin{equation*}
  S_{N}(x,y)=\sum_{j,k=0}^{2N-1}\psi_{j}^{\prime}(x)\mu_{jk}\psi_{k}(y).
\end{equation*}
Since $M$ is antisymmetric,
\begin{equation*}
  S_{N}(y,x)=\sum_{j,k=0}^{2N-1}\psi_{j}^{\prime}(y)\mu_{jk}\psi_{k}(x)=-\sum_{j,k=0}^{2N-1}\psi_{j}^{\prime}(y)\mu_{kj}\psi_{k}(x)= -\sum_{j,k=0}^{2N-1}\psi_{j}(y)\mu_{kj}\psi_{k}^{\prime}(x),
\end{equation*}
after re-indexing. Note that
\begin{equation*}
  \epsilon\,S_{N}(x,y)=\sum_{j,k=0}^{2N-1}\epsilon\,D\,\psi_{j}(x)\mu_{jk}\psi_{k}(y) =\sum_{j,k=0}^{2N-1}\psi_{j}(x)\mu_{jk}\psi_{k}(y),
\end{equation*}
and
\begin{equation*}
  -\frac{d}{d\,y}\,S_{N}(x,y)=\sum_{j,k=0}^{2N-1}\psi_{j}^{\prime}(x)\mu_{jk}\psi_{k}^{\prime}(y).
\end{equation*}
Thus we can now write succinctly
\begin{equation}
  \label{gsekernel}
  K_{N}(x,y)=
  \left(
    \begin{array}{cc}
      S_{N}(x,y) & -\frac{d}{d\,y}\,S_{N}(x,y) \\
      \epsilon\,S_{N}(x,y) & S_{N}(y,x)
    \end{array}
  \right).
\end{equation}
To summarize, we have shown that $G_{4,N}^{2}(t,\lambda)=C_{4}^{(N)}\det(I-K_{4,N}f)$. Setting $f\equiv 0$ on both sides (where the original definition of $G_{4,N}(t,\lambda)$ as an integral is used on the left) shows that $C_{4}^{(N)}=1$. Thus
\begin{equation}
  \label{eq:5}
  G_{4,N}(t,\lambda)=\sqrt{D_{4,N}(t,\lambda)},
\end{equation}
where we define
\begin{equation}
  D_{4,N}(t,\lambda)=\det(I-K_{4,N}f),
\end{equation}
and $K_{4,N}$ is the integral operator with matrix kernel \eqref{gsekernel}.

\section{Gaussian specialization}

We would like to specialize the above results to the case of a Gaussian weight function
\begin{equation}
  \label{gseweight}
  w(x)=\exp\left(x^{2}\right)
\end{equation}
and indicator function
\begin{equation*}
  f(x)=-\lambda\,\rchi_{_{J}}\quad , \quad J=(t,\infty).
\end{equation*}
We want the matrix
\begin{equation*}
  M=\left\{\int \left(\psi_{j}(x)\,\psi_{k}^{\prime}(x) - \,\psi_{k}(x)\,\psi_{j}^{\prime}(x)\right)\,d\,x\right\}_{_{0\leq j,k\leq 2\,N-1}}
\end{equation*}
to be the direct sum of $N$ copies of
\begin{equation*}
  Z=\left(
    \begin{array}{cc}
      0 & 1 \\
      -1 & 0
    \end{array}
  \right),
\end{equation*}
so that the formulas are the simplest possible, since then $\mu_{jk}$ can only be $0$ or $\pm 1$.  In that case $M$ would be skew--symmetric so that $M^{-1}=-M$. In terms of the integrals defining the entries of $M$ this means that we would like to have
\begin{equation*}
  \int \left( \psi_{2j}(x)\,\frac{d}{d\,x}\,\psi_{2k+1}(x) - \,\psi_{2k+1}(x)\,\frac{d}{d\,x}\,\psi_{2j}(x)\right)\,d\,x = \delta_{j,k},
\end{equation*}
\begin{equation*}
  \int \left( \psi_{2j+1}(x)\,\frac{d}{d\,x}\,\psi_{2k}(x) - \,\psi_{2k}(x)\,\frac{d}{d\,x}\,\psi_{2j+1}(x)\right)\,d\,x = -\delta_{j,k}
\end{equation*}
and otherwise
\begin{equation*}
  \int \left( \psi_{j}(x)\,\frac{d}{d\,x}\,\psi_{k}(x) - \,\psi_{j}(x)\,\frac{d}{d\,x}\,\psi_{k}(x)\right)\,d\,x = 0.
\end{equation*}
It is easier to treat this last case if we replace it with three non-exclusive conditions
\begin{equation*}
  \int \left(\psi_{2j}(x)\,\frac{d}{d\,x}\,\psi_{2k}(x) -\,\psi_{2k}(x)\,\frac{d}{d\,x}\,\psi_{2j}(x)\right)\,d\,x = 0,
\end{equation*}
\begin{equation*}
  \int \left( \psi_{2j+1}(x)\,\frac{d}{d\,x}\,\psi_{2k+1}(x) - \,\psi_{2k+1}(x)\,\frac{d}{d\,x}\,\psi_{2j+1}(x)\right)\,d\,x = 0,
\end{equation*}
(so when the parity is the same for  $j,k$, which takes care of diagonal entries, among others) and
\begin{equation*}
  \int \left( \psi_{j}(x)\,\frac{d}{d\,x}\,\psi_{k}(x) - \,\psi_{j}(x)\,\frac{d}{d\,x}\,\psi_{k}(x)\right)\,d\,x = 0,
\end{equation*}
whenever $|j-k|>1$, which targets entries outside of the tridiagonal.
Define
\begin{equation}
  \label{gsevarphi}
  \varphi_{k}(x)=\frac{1}{c_{k}}H_{k}(x)\,\exp(-x^{2}/2)\quad\textrm{for}\quad c_{k}=\sqrt{2^{k}k!\sqrt{\pi}}
\end{equation}
where the $H_{k}$ are the usual Hermite polynomials defined by the orthogonality condition
\begin{equation*}
  \int_{\mathbb{R}} H_{j}(x)\,H_{k}(x)\,e^{-x^{2}}\,d\,x = c_{j}^{2}\,\delta_{j,k}.
\end{equation*}
Then it follows that
\begin{equation*}
  \int_{\mathbb{R}} \varphi_{j}(x)\,\varphi_{k}(x)\,d\,x = \delta_{j,k}.
\end{equation*}
Now let
\begin{equation*}
  \psi_{_{2j+1}}(x) = \frac{1}{\sqrt{2}}\,\varphi_{_{2j+1}}\left(x\right)\qquad \psi_{_{2j}}(x) = -\frac{1}{\sqrt{2}}\,\epsilon\,\varphi_{_{2j+1}}\left(x\right)
\end{equation*}
This definition satisfies our earlier requirement that $\psi_{j}=p_{j}\,w^{1/2}$ with $w$ defined in \eqref{gseweight}. In particular we have in this case
\begin{equation*}
  p_{2j+1}(x)=\frac{1}{c_{j}\sqrt{2}}\,H_{2j+1}\left(x\right).
\end{equation*}
Let $\epsilon$ as in \eqref{epsilonop}, and $D$ denote the operator that acts by differentiation with respect to $x$ as before, so that $D\,\epsilon=\epsilon\,D=I$. It follows that
\begin{equation*}
  \begin{aligned}
    \int_{\mathbb{R}} & \left[\psi_{_{2j}}(x)\frac{d}{d\,x}\,\psi_{_{2k+1}}(x)-\psi_{_{2k+1}}(x)\frac{d}{d\,x}\,\psi_{_{2j}}(x)\right]d\,x\\
    & = \frac{1}{2}\int_{\mathbb{R}}\left[-\epsilon\,\varphi_{_{2j+1}}\left(x\right)\frac{d}{d\,x}\,\varphi_{_{2k+1}}\left(x\right) + \varphi_{_{2k+1}}\left(x\right)\frac{d}{d\,x}\,\epsilon\,\varphi_{_{2j+1}}\left(x\right)\right]d\,x\\
    & = \frac{1}{2}\int_{\mathbb{R}}\left[ -\epsilon\,\varphi_{_{2j+1}}\left(x\right)\frac{d}{d\,x}\,\varphi_{_{2k+1}}\left(x\right) +\varphi_{_{2k+1}}\left(x\right)\varphi_{_{2j+1}}\left(x\right)\right]d\,x\\
  \end{aligned}
\end{equation*}
We integrate the first term by parts and use the fact that
\begin{equation*}
  \frac{d}{d\,x}\,\epsilon\,\varphi_{_{j}}\left(x\right) =\varphi_{_{j}}\left(x\right)
\end{equation*}
and also that $\varphi_{_{j}}$ vanishes at the boundary (i.e. $\varphi_{_{j}}\left(\pm\infty\right) =0$) to obtain
\begin{equation*}
  \begin{aligned}
    \int_{\mathbb{R}} & \left[\psi_{_{2j}}(x)\frac{d}{d\,x}\,\psi_{_{2k+1}}(x)-\psi_{_{2k+1}}(x)\frac{d}{d\,x}\,\psi_{_{2j}}(x)\right]d\,x\\
    & = \frac{1}{2} \int_{\mathbb{R}}\left[-\epsilon\,\varphi_{_{2j+1}}\left(x\right)\frac{d}{d\,x}\,\varphi_{_{2k+1}}\left(x\right) +\varphi_{_{2k+1}}\left(x\right)\varphi_{_{2j+1}}\left(x\right)\right]d\,x\\
    & = \frac{1}{2}  \int_{\mathbb{R}}\left[ \varphi_{_{2j+1}}\left(x\right)\,\varphi_{_{2k+1}}\left(x\right)+\varphi_{_{2k+1}}\left(x\right)\varphi_{_{2j+1}}\left(x\right)\right]d\,x\\
    & = \frac{1}{2}  \int_{\mathbb{R}}\left[ \varphi_{_{2j+1}}\left(x\right)\,\varphi_{_{2k+1}}\left(x\right)+\varphi_{_{2k+1}}\left(x\right)\varphi_{_{2j+1}}\left(x\right)\right]d\,x\\
    &= \frac{1}{2} \left(\delta_{j,k}+ \delta_{j,k}\right)\\
    & =\delta_{j,k},
  \end{aligned}
\end{equation*}
as desired. Similarly
\begin{equation*}
  \begin{aligned}
    \int_{\mathbb{R}} & \left[\psi_{_{2j+1}}(x)\frac{d}{d\,x}\,\psi_{_{2k}}(x)-\psi_{_{2k}}(x)\frac{d}{d\,x}\,\psi_{_{2j+1}}(x)\right]d\,x\\
    & = \frac{1}{2}\int_{\mathbb{R}}\left[-\varphi_{_{2j+1}}\left(x\right)\frac{d}{d\,x}\,\epsilon\,\varphi_{_{2k+1}}\left(x\right) +\epsilon\,\varphi_{_{2k+1}}\left(x\right)\frac{d}{d\,x}\,\varphi_{_{2j+1}}\left(x\right)\right]d\,x \\
    & = -\delta_{j,k}.
  \end{aligned}
\end{equation*}
Moreover,
\begin{equation*}
  p_{2j+1}(x)=\frac{1}{c_{j}\sqrt{2}}\,H_{2j+1}\left(x\right)
\end{equation*}
is certainly an odd function, being the multiple of and odd Hermite polynomial. On the other hand, one easily checks that $\epsilon$ maps odd functions to even functions on $L^{2}(\mathbb{R})$. Therefore
\begin{equation*}
  p_{2j}(x)=-\frac{1}{c_{j}\sqrt{2}}\,\epsilon\,H_{2j+1}\left(x\right)
\end{equation*}
is an even function, and it follows that
\begin{equation*}
  \begin{aligned}\int_{\mathbb{R}} & \left[ \psi_{_{2k}}(x)\frac{d}{d\,x}\,\psi_{_{2j}}(x)-\psi_{_{2j}}(x)\frac{d}{d\,x}\,\psi_{_{2k}}(x)\right]d\,x\\
    & = \int_{\mathbb{R}}\left[p_{2j}(x)\frac{d}{d\,x}\,p_{2k}(x)-p_{2k}(x)\frac{d}{d\,x}\,p_{2j}(x)\right]w(x)\,d\,x\\
    & = 0,
  \end{aligned}
\end{equation*}
since both terms in the integrand are odd functions, and the weight function is even. Similarly,
\begin{equation*}
  \begin{aligned}
    \int_{\mathbb{R}} & \left[\psi_{_{2k+1}}(x)\frac{d}{d\,x}\,\psi_{_{2j+1}}(x)-\psi_{_{2j+1}}(x)\frac{d}{d\,x}\,\psi_{_{2k+1}}(x)\right]d\,x\\
    & = \int_{\mathbb{R}}\left[p_{2j+1}(x)\frac{d}{d\,x}\,p_{2k+1}(x)-p_{2k+1}(x)\frac{d}{d\,x}\,p_{2j+1}(x)\right]w(x)\,d\,x\\
    & = 0.
  \end{aligned}
\end{equation*}
Finally it is easy to see that if $|j-k|>1$ then
\begin{equation*}
  \int_{\mathbb{R}} \left[\psi_{_{j}}(x)\frac{d}{d\,x}\,\psi_{_{k}}(x)-\psi_{_{j}}(x)\frac{d}{d\,x}\,\psi_{_{k}}(x)\right]d\,x=0.
\end{equation*}
Indeed both  differentiation and the action of $\epsilon$ can only ``shift'' the indices by $1$. Thus by orthogonality of the $\varphi_{j}$, this integral will always be $0$. Hence by choosing
\begin{equation*}
  \psi_{_{2j+1}}(x) = \frac{1}{\sqrt{2}}\,\varphi_{_{2j+1}}\left(x\right),\qquad \psi_{_{2j}}(x) = -\frac{1}{\sqrt{2}}\,\epsilon\,\varphi_{_{2j+1}}\left(x\right),
\end{equation*}
we force the matrix
\begin{equation*}
  M=\left\{\int_{\mathbb{R}} \left(\psi_{j}(x)\,\psi_{k}^{\prime}(x) -\,\psi_{k}(x)\,\psi_{j}^{\prime}(x)\right)\,d\,x\right\}_{_{0\leq
j,k\leq 2\,n-1}}
\end{equation*}
to be the direct sum of $N$ copies of
\begin{equation*}
  Z=\left(\begin{array}{cc} 0 & 1 \\ -1 & 0 \end{array}\right).
\end{equation*}
Hence $M^{-1}=-M$ where $M^{-1}=\{\mu_{j,k}\}_{j,k=0,2N-1}$. Moreover, with our above choice, $ \mu_{j,k}=0$ if $j,k$ have the same parity or $|j-k|>1$, and $\mu_{2k,2j+1}=\delta_{jk}=-\mu_{2j+1,2k}$ for $j,k=0.\ldots,N-1$.
Therefore
\begin{equation*}
  \begin{aligned}
    S_{N}(x,y) & =-\sum_{j,k=0}^{2N-1}\psi_{j}^{\prime}(x)\mu_{jk}\psi_{k}(y)\\
    & = -\sum_{j=0}^{N-1}\psi_{2j}^{\prime}(x)\,\psi_{2j+1}(y) + \sum_{j=0}^{N-1}\psi_{2j+1}^{\prime}(x)\,\psi_{2j}(y) \\
    & =   \frac{1}{2}\left[ \sum_{j=0}^{N-1}\varphi_{_{2j+1}}\left(x\right) \varphi_{_{2j+1}}\left(y\right)-\sum_{j=0}^{N-1}\varphi_{_{2j+1}}^{\prime}\left(x\right)\epsilon\,\varphi_{_{2j+1}}\left(y\right) \right].
  \end{aligned}
\end{equation*}
Recall that the $H_{j}$ satisfy  the differentiation formulas (see for example \cite{Andr1}, p. 280)
\begin{equation}
  \label{hermrec1}
  H_{j}^{\prime}(x) = 2\,x\,H_{j}(x)-H_{j-1}(x)\quad j=1,2,\ldots
\end{equation}
\begin{equation}
  \label{hermrec2}
  H_{j}^{\prime}(x) = 2\,j\,H_{j-1}(x)\quad j=1,2,\ldots
\end{equation}
Combining \ref{gsevarphi} and \ref{hermrec1} yields
\begin{equation}
  \label{eq:2}
  \varphi_{j}^{\prime}(x) = x\,\varphi_{j}(x) - \frac{c_{j+1}}{c_{j}}\,\varphi_{j+1}(x).
\end{equation}
Similarly, from \ref{gsevarphi} and \ref{hermrec2} we have
\begin{equation}
  \label{eq:3}
  \varphi_{j}^{\prime}(x) = -x\,\varphi_{j}(x) + 2j\,\frac{c_{j-1}}{c_{j}}\,\varphi_{j-1}(x).
\end{equation}
Combining \eqref{eq:2} and \eqref{eq:3}, we obtain
\begin{equation}
  \label{eq:4}
    \varphi_{j}^{\prime}(x) = \sqrt{\frac{j}{2}}\,\varphi_{j-1}(x) - \sqrt{\frac{j+1}{2}}\,\varphi_{j+1}(x).
\end{equation}
Let $\varphi=\left(\begin{array}{c}\varphi_{1} \\ \varphi_{2}\\ \vdots \end{array}\right)$ and $\varphi^{\prime}=\left(\begin{array}{c}\varphi_{1}^{\prime} \\ \varphi_{2}^{\prime}\\ \vdots \end{array}\right)$. Then we can rewrite \eqref{eq:4} as
\begin{equation*}
  \varphi^{\prime}=A\,\varphi
\end{equation*}
where $A=\{a_{j,k}\}$ is the infinite antisymmetric tridiagonal matrix with $a_{j,j-1}=\sqrt{\frac{j}{2}}$. Hence,
\begin{equation*}
  \varphi_{j}^{\prime}(x)=\sum_{k\geq 0}a_{jk}\varphi_{k}(x).
\end{equation*}
Moreover, using the fact that $D\,\epsilon=\epsilon\,D=I$ we also have
\begin{equation*}
  \begin{aligned}\varphi_{j}\left(x\right)
    & =\epsilon\,\varphi_{j}^{\prime}\left(x\right)=\epsilon\left(\sum_{k\geq 0}a_{jk}\varphi_{k}\left(x\right)\right)
    \\ & = \sum_{k\geq 0}a_{jk}\epsilon\,\varphi_{k}\left(x\right).
  \end{aligned}
\end{equation*}
Combining the above results, we have
\begin{equation*}
  \begin{aligned}
    \sum_{j=0}^{N-1}\varphi_{2j+1}^{\prime}\left(x\right)\epsilon\,\varphi_{2j+1}\left(y\right) & = \sum_{j=0}^{N-1}\sum_{k\geq 0} a_{2j+1,k}\varphi_{k}\left(x\right)\epsilon\,\varphi_{2j+1}\left(x\right)\\
    & = -\sum_{j=0}^{N-1}\sum_{k\geq 0} a_{k,2j+1}\varphi_{k}\left(x\right)\epsilon\,\varphi_{2j+1}\left(x\right).
  \end{aligned}
\end{equation*}
Note that $a_{k,2j+1}=0$ unless $|k-(2j+1)|=1$, that is unless $k$ is even. Thus we can rewrite the sum as
\begin{equation*}
  \begin{aligned}
    \sum_{j=0}^{N-1}\varphi_{2j+1}^{\prime}\left(x\right) \epsilon\,\varphi_{2j+1}\left(y\right)  & =-\sum_{\substack{k,j\geq 0\\ k\ \textrm{even} \\ k\leq 2N}} a_{k,j}\varphi_{k}\left(x\right)\epsilon\,\varphi_{j}\left(y\right) - a_{_{2N,2N+1}}\varphi_{_{2N}}\left(x\right)\epsilon\,\varphi_{_{2N+1}}\left(y\right)\\
    &  = -\sum_{\substack{k\geq 0\\ k\ \textrm{even} \\ k\leq 2N}} \varphi_{k}\left(x\right)\sum_{j\geq 0}a_{k,j}\epsilon\,\varphi_{j}\left(y\right) + a_{_{2N,2N+1}}\,\varphi_{_{2N}}\left(x\right)\epsilon\,\varphi_{_{2N+1}}\left(y\right)
  \end{aligned}
\end{equation*}
where the last term takes care of the fact that we are counting  an extra term in the sum that was not present before. The sum over $j$  on the right is just $\varphi_{k}(y)$, and  $a_{_{2N,2N+1}}=-\sqrt{\frac{2N+1}{2}}$. Therefore
\begin{equation*}
  \begin{aligned}
    \sum_{j=0}^{N-1}\varphi_{2j+1}^{\prime}\left(x\right)\epsilon\,\varphi_{2j+1}\left(y\right) & = \sum_{\substack{k\geq 0\\ k\ \textrm{even} \\ k\leq 2N}} \varphi_{k}\left(x\right)\varphi_{k}\left(y\right)-\sqrt{\frac{2N+1}{2}}\,\varphi_{_{2N}}\left(x\right)\epsilon\,\varphi_{_{2N+1}}\left(y\right)\\
    & = \sum_{j=0}^{N} \varphi_{2j}\left(x\right)\varphi_{2j}\left(y\right)-\sqrt{\frac{2N+1}{2}}\,\varphi_{_{2N}}\left(x\right)\epsilon\,\varphi_{_{2N+1}}\left(y\right).
  \end{aligned}
\end{equation*}
It follows that
\begin{equation*}
  \begin{aligned}
    S_{N}(x,y) =   \frac{1}{2}\left[ \sum_{j=0}^{2N}\varphi_{_{j}}\left(x\right) \varphi_{_{j}}\left(y\right)-\sqrt{\frac{2N+1}{2}}\, \varphi_{_{2N}}\left(x\right)\epsilon\,\varphi_{_{2N+1}} \left(y\right) \right].
  \end{aligned}
\end{equation*}
We redefine
\begin{equation}
  \label{newgseS}
  S_{N}(x,y) = \sum_{n=0}^{2N}\varphi_{_{n}}\left(x\right) \varphi_{_{n}}\left(y\right)=S_{N}(y,x)
\end{equation}
so that the top left entry of $K_{N}(x,y)$ is
\begin{equation*}
  S_{N}(x,y) +  \sqrt{\frac{2N+1}{2}}\, \varphi_{_{2N}}\left(x\right)\epsilon\,\varphi_{_{2N+1}}\left(y\right).
\end{equation*}
If $S_{N}$ is the operator with kernel $S_{N}(x,y)$ then integration by parts gives
\begin{equation*}
  S_{N}D f = \int_{\mathbb{R}} S(x,y) \frac{d}{d\,y}f(y)\,d\,y = \int_{\mathbb{R}}  \left(-\frac{d}{d\,y} S_{N}(x,y)\right) f(y)\,d\,y,
\end{equation*}
so that $-\frac{d}{d\,y} S_{N}(x,y)$ is in fact the kernel of $S_{N}D$. Therefore \eqref{eq:5} now holds with $K_{4,N}$ being the integral operator with matrix kernel $K_{4,N}(x,y)$ whose $(i,j)$--entry $K_{4,N}^{(i,j)}(x,y)$ is given by
\begin{equation*}
  \begin{aligned}
    K_{4,N}^{(1,1)}(x,y)& =\frac{1}{2}\left[ S_{N}(x,y) + \sqrt{\frac{2N+1}{2}}\, \varphi_{_{2N}}\left(x\right)\epsilon\,\varphi_{_{2N+1}} \left(y\right)\right], \\
    K_{4,N}^{(1,2)}(x,y) & = \frac{1}{2}\left[ SD_{N}(x,y) -\frac{d}{d\,y} \, \left(\sqrt{\frac{2N+1}{2}}\, \varphi_{_{2N}}\left(x\right)\epsilon\,\varphi_{_{2N+1}} \left(y\right) \right) \right], \\
    K_{4,N}^{(2,1)}(x,y) & = \frac{\epsilon}{2}\left[ S_{N}(x,y) +  \sqrt{\frac{2N+1}{2}}\, \varphi_{_{2N}}\left(x\right)\epsilon\,\varphi_{_{2N+1}}\left(y\right) \right],\\
    K_{4,N}^{(2,2)}(x,y) & = \frac{1}{2}\left[ S_{N}(x,y) + \sqrt{\frac{2N+1}{2}}\, \epsilon\,\varphi_{_{2N+1}}\left(x\right)\,\varphi_{_{2N}}\left(y\right)\right].
  \end{aligned}
\end{equation*}
We let $2N+1\to N$ so that $N$ is assumed to be odd from now on (this will not matter in the  end since we will take $N\to\infty$). Therefore the $K_{4,N}^{(i,j)}(x,y)$ are given by
\begin{equation*}
  \begin{aligned}
    K_{4,N}^{(1,1)}(x,y)& =\frac{1}{2}\left[ S_{N}(x,y) + \sqrt{\frac{N}{2}}\, \varphi_{_{N-1}}(x)\epsilon\,\varphi_{_{N}} (y)\right], \\
    K_{4,N}^{(1,2)}(x,y) & = \frac{1}{2}\left[ SD_{N}(x,y) -\sqrt{\frac{N}{2}}\, \varphi_{_{N-1}}(x)\,\varphi_{_{N}} (y) \right], \\
    K_{4,N}^{(2,1)}(x,y) & = \frac{\epsilon}{2}\left[ S_{N}(x,y) +  \sqrt{\frac{N}{2}}\, \varphi_{_{N-1}}(x)\epsilon\,\varphi_{_{N}}(y) \right], \\
    K_{4,N}^{(2,2)}(x,y) & = \frac{1}{2}\left[ S_{N}(x,y) + \sqrt{\frac{N}{2}}\, \epsilon\,\varphi_{_{N}} (x)\,\varphi_{_{N-1}}(y)\right],
  \end{aligned}
\end{equation*}
where
\begin{equation*}
S_{N}(x,y) = \sum_{n=0}^{N-1}\varphi_{_{n}}(x) \varphi_{_{n}}(y).
\end{equation*}
Define
\begin{equation*}
  \varphi(x)=\left(\frac{N}{2}\right)^{1/4}\varphi_{_{N}}(x)\qquad \psi(x)=\left(\frac{N}{2}\right)^{1/4}\varphi_{_{N-1}}(x),
\end{equation*}
so that
\begin{equation*}
  \begin{aligned}
    K_{4,N}^{(1,1)}(x,y)& =\frac{1}{2}\,\rchi(x)\,\left[ S_{N}(x,y) + \psi(x)\epsilon\,\varphi(y)\right] \,\rchi(y),\\
    K_{4,N}^{(1,2)}(x,y) & = \frac{1}{2}\,\rchi(x)\,\left[ SD_{N}(x,y) - \psi(x)\varphi(y)  \right]\,\rchi(y),\\
    K_{4,N}^{(2,1)}(x,y) & = \frac{1}{2}\,\rchi(x)\,\left[ \epsilon S_{N}(x,y) + \epsilon\,\psi(x)\epsilon\,\varphi(y) \right]\,\rchi(y),\\
    K_{4,N}^{(2,2)}(x,y) & = \frac{1}{2}\,\rchi(x)\,\left[ S_{N}(x,y) + \epsilon\,\varphi (x)\,\psi(y)\right]\,\rchi(y).
  \end{aligned}
\end{equation*}
Notice that
\begin{equation*}
  \begin{aligned}
    \frac{1}{2}\,\rchi\,\left( S + \psi\otimes\,\epsilon\,\varphi \right)\,\rchi & \doteq K_{4,N}^{(1,1)}(x,y),\\
    \frac{1}{2}\,\rchi\,\left(SD - \psi\otimes\,\varphi  \right)\,\rchi &\doteq K_{4,N}^{(1,2)}(x,y),\\
    \frac{1}{2}\,\rchi\,\left(\epsilon\,S + \epsilon\,\psi\otimes\,\epsilon\,\varphi \right)\,\rchi &\doteq K_{4,N}^{(2,1)}(x,y),\\
    \frac{1}{2}\,\rchi\,\left(S + \epsilon\,\varphi\otimes\,\epsilon\,\psi \right)\,\rchi &\doteq K_{4,N}^{(2,2)}(x,y).\\
  \end{aligned}
\end{equation*}
Therefore
\begin{equation}
  \label{eq:6}
  K_{4,N}=\frac{1}{2}\,\rchi\,\left(\begin{array}{cc} S + \psi\otimes\,\epsilon\,\varphi  & SD - \psi\otimes\,\varphi \\
      \epsilon\,S + \epsilon\,\psi\otimes\,\epsilon\,\varphi  &  S + \epsilon\,\varphi\otimes\,\psi  \end{array}\right)\,\rchi.
\end{equation}
Note that this is identical to the corresponding operator for $\beta=4$ obtained by Tracy and Widom in \cite{Trac2}, the only difference being that $\varphi$, $\psi$, and hence also $S$, are redefined to depend on $\lambda$. This will affect boundary conditions for the differential equations we will obtain later.

\section{Edge scaling}

\subsection{Reduction of the determinant}

We want to compute the Fredholm determinant \eqref{eq:5} with $K_{4,N}$ given by \eqref{eq:6} and $f=\rchi_{(t,\infty)}$. This is the determinant of an operator on $L^{2}(J)~\times~L^{2}(J)$. Our first task will be to rewrite the determinant as that of an operator on $L^{2}(J)$. This part follows exactly the proof in \cite{Trac2}. To begin, note that
\begin{equation}
 \label{eq:7}
  \com{S}{D}=\varphi\otimes\psi + \psi\otimes\varphi
\end{equation}
so that, using the fact that $D\,\epsilon=\epsilon\,D=I$,
\begin{eqnarray}
 \label{eq:8}
  \com{\epsilon}{S} &=&\epsilon\,S-S\,\epsilon\nonumber\\
  &=& \epsilon\,S\,D\,\epsilon-\epsilon\,D\,S\,\epsilon = \epsilon\,\com{S}{D}\,\epsilon\nonumber\\
  &=& \epsilon\,\varphi\otimes\psi\,\epsilon + \epsilon\,\psi\otimes\varphi\,\epsilon\nonumber\\
  &=&\epsilon\,\varphi\otimes\epsilon^{t}\psi + \epsilon\,\psi\otimes\epsilon^{t}\,\varphi\nonumber\\
  &=& - \epsilon\,\varphi\otimes\epsilon\,\psi -
  \epsilon\,\psi\otimes\epsilon\,\varphi,
\end{eqnarray}
where the last equality follows from the fact that $\epsilon^{t}=-\epsilon$. We thus have
\begin{eqnarray*}
  D\,\left(\epsilon\,S + \epsilon\,\psi\otimes\epsilon\,\varphi\right) & = & S + \psi\otimes\epsilon\,\varphi,\\
  D\,\left(\epsilon\,S\,D - \epsilon\,\psi\otimes\varphi\right) & =
  & S\,D - \psi\otimes\varphi.
\end{eqnarray*}
The expressions on the right side are the top matrix entries in \eqref{eq:6}. Thus the first row of $K_{4,N}$ is, as a vector,
\begin{equation*}
  D\,\left(\epsilon\,S + \epsilon\,\psi\otimes\epsilon\,\varphi, \epsilon\,S\,D - \epsilon\,\psi\otimes\varphi\right).
\end{equation*}
Now \eqref{eq:8} implies that
\begin{equation*}
  \epsilon\,S + \epsilon\,\psi\otimes\epsilon\,\varphi = S\,\epsilon -\epsilon\,\varphi\otimes\epsilon\,\psi.
\end{equation*}
Similarly \eqref{eq:7} gives
\begin{equation*}
  \epsilon\,\com{S}{D} = \epsilon\,\varphi\otimes\psi + \epsilon\,\psi\otimes\varphi,
\end{equation*}
so that
\begin{equation*}
  \epsilon\,S\,D - \epsilon\,\psi\otimes\varphi = \epsilon\,D\,S + \epsilon\,\varphi\otimes\psi = S + \epsilon\,\varphi\otimes\psi.
\end{equation*}
Using these expressions we can rewrite the first row of $K_{4,N}$ as
\begin{equation*}
  D\,\left(S\,\epsilon - \epsilon\,\varphi\otimes\epsilon\,\psi, S + \epsilon\,\varphi\otimes\psi\right).
\end{equation*}
Now use \eqref{eq:8} to show the second row of $K_{4,N}$ is
\begin{equation*}
  \left(S\,\epsilon - \epsilon\,\varphi\otimes\epsilon\,\psi, S + \epsilon\,\varphi\otimes\psi\right).
\end{equation*}
Therefore,
\begin{eqnarray*}
  K_{4,N} &=& \rchi\,\left(
    \begin{array}{cc}
      D\,\left(S\,\epsilon - \epsilon\,\varphi\otimes\epsilon\,\psi\right) & D\,\left(S + \epsilon\,\varphi\otimes\psi\right) \\
      S\,\epsilon - \epsilon\,\varphi\otimes\epsilon\,\psi &  S + \epsilon\,\varphi\otimes\psi
    \end{array}
  \right)\,\rchi \\
  & =& \left(
    \begin{array}{cc}
      \rchi\,D & 0 \\
      0 & \rchi
    \end{array}
  \right)
  \left(
    \begin{array}{cc}
      \left(S\,\epsilon - \epsilon\,\varphi\otimes\epsilon\,\psi\right)\,\rchi & \left(S + \epsilon\,\varphi\otimes\psi\right)\,\rchi \\
      \left(S\,\epsilon - \epsilon\,\varphi\otimes\epsilon\,\psi\right)\rchi & \left(S + \epsilon\,\varphi\otimes\psi\right)\rchi
    \end{array}
  \right).
\end{eqnarray*}
Since $K_{4,N}$ is of the form $A\,B$, we can use \ref{detthm}  and deduce that $D_{4,N}(s,\lambda)$ is unchanged if instead we take $ K_{4,N}$ to be
\begin{eqnarray*}
  K_{4,N} & = & \left(
    \begin{array}{cc}
      \left(S\,\epsilon - \epsilon\,\varphi\otimes\epsilon\,\psi\right)\,\rchi & \left(S + \epsilon\,\varphi\otimes\psi\right)\,\rchi \\
      \left(S\,\epsilon - \epsilon\,\varphi\otimes\epsilon\,\psi\right)\,\rchi & \left(S + \epsilon\,\varphi\otimes\psi\right)\,\rchi
    \end{array}
  \right)
  \left(
    \begin{array}{cc}
      \rchi\,D & 0 \\
      0 & \rchi
    \end{array}
  \right)\\
  & = & \left(\begin{array}{cc}
      \left(S\,\epsilon - \epsilon\,\varphi\otimes\epsilon\,\psi\right)\,\rchi\,D & \left(S + \epsilon\,\varphi\otimes\psi\right)\,\rchi \\
      \left(S\,\epsilon - \epsilon\,\varphi\otimes\epsilon\,\psi\right)\rchi\,D & \left(S + \epsilon\,\varphi\otimes\psi\right)\rchi
    \end{array}\right).
\end{eqnarray*}
Therefore
\begin{equation}
  D_{4,N}(s,\lambda)=\det\left(
    \begin{array}{cc}
      I - \frac{1}{2}\,\left(S\,\epsilon - \epsilon\,\varphi\otimes\epsilon\,\psi\right)\,\lambda\,\rchi\,D & - \frac{1}{2}\,\left(S + \epsilon\,\varphi\otimes\psi\right)\,\lambda\,\rchi \\
      - \frac{1}{2}\,\left(S\,\epsilon - \epsilon\,\varphi\otimes\epsilon\,\psi\right)\,\lambda\,\rchi\,D & I - \frac{1}{2}\,\left(S + \epsilon\,\varphi\otimes\psi\right)\,\lambda\,\rchi
    \end{array}
  \right).
\end{equation}
Now we perform row and column operations on the matrix to simplify it, which do not change the Fredholm determinant. Justification of these operations is given in \cite{Trac2}. We start by subtracting row 1 from row 2 to get
\begin{equation*}
  \left(\begin{array}{cc}
      I - \frac{1}{2}\,\left(S\,\epsilon - \epsilon\,\varphi\otimes\epsilon\,\psi\right)\,\lambda\,\rchi\,D & - \frac{1}{2}\,\left(S + \epsilon\,\varphi\otimes\psi\right)\,\lambda\,\rchi \\
      - I  & I \end{array}\right).
\end{equation*}
Next, adding column 2 to column 1 yields
\begin{equation*}
  \left(\begin{array}{cc}
      I - \frac{1}{2}\,\left(S\,\epsilon - \epsilon\,\varphi\otimes\epsilon\,\psi\right)\,\lambda\,\rchi\,D - \frac{1}{2}\,\left(S + \epsilon\,\varphi\otimes\psi\right)\,\lambda\,\rchi & - \frac{1}{2}\,\left(S + \epsilon\,\varphi\otimes\psi\right)\,\lambda\,\rchi \\
      0  & I
    \end{array}
  \right).
\end{equation*}

Thus the determinant we want equals the determinant of
\begin{equation}
  \label{eq:9}
  I - \frac{1}{2}\,\left(S\,\epsilon - \epsilon\,\varphi\otimes\epsilon\,\psi\right)\,\lambda\,\rchi\,D - \frac{1}{2}\,\left(S + \epsilon\,\varphi\otimes\psi\right)\,\lambda\,\rchi
\end{equation}
So we have reduced the problem from the computation of the Fredholm determinant of an operator on $L^{2}(J)~\times~L^{2}(J)$, to that of an operator on $L^{2}(J)$.

\subsection{Differential equations}

Next we want to write the operator in \eqref{eq:9} in the form
\begin{equation}
  \left(I-  K_{2, N}\right)\left(I - \sum_{i=1}^{L}\alpha_{i}\otimes\beta_{i}\right),
\end{equation}
where the $\alpha_{i}$  and $\beta_{i}$ are functions in $L^{2}(J)$. In other words, we want to rewrite the determinant for the GSE case as a finite dimensional perturbation of the corresponding GUE determinant. The Fredholm determinant of the product is then the product of the determinants. The limiting form for the GUE part is already known, and we can just focus on finding a limiting form for the determinant of the finite dimensional piece. It is here that the proof must be modified from that in \cite{Trac2}. A little rearrangement of \eqref{eq:9} yields (recall $\epsilon^{t}=-\epsilon$)
\begin{equation*}
  I - \frac{\lambda}{2}\,S\,\rchi  - \frac{\lambda}{2}\,S\,\epsilon\,\rchi\,D - \frac{\lambda}{2}\,\epsilon\,\varphi\,\otimes\,\rchi\,\psi - \frac{\lambda}{2}\,\epsilon\,\varphi\,\otimes\,\psi\,\epsilon\,\rchi\,D.
\end{equation*}
Writing $\epsilon\,\com{\rchi}{D}+\rchi$ for $\epsilon\,\rchi\,D$ and simplifying gives
\begin{equation*}
  I - \lambda\,S\,\rchi - \lambda\,\epsilon\,\varphi\,\otimes\,\psi\,\rchi - \frac{\lambda}{2}\,S\,\epsilon\,\com{\rchi}{D} - \frac{\lambda}{2}\,\epsilon\,\varphi\,\otimes\,\psi\,\epsilon\,\com{\rchi}{D}.
\end{equation*}
Let $\sqrt{\lambda}\,\varphi\to\varphi$, and $\sqrt{\lambda}\,\psi\to\psi$ so that $\lambda\,S\to S$ and \eqref{eq:9} goes to
\begin{equation*}
  I - S\,\rchi - \epsilon\,\varphi\,\otimes\,\psi\,\rchi - \frac{1}{2}\,S\,\epsilon\,\com{\rchi}{D} - \frac{1}{2}\,\epsilon\,\varphi\,\otimes\,\psi\,\epsilon\,\com{\rchi}{D}.
\end{equation*}
Now we define $R:=(I-S\,\rchi)^{-1}\,S\,\rchi=(I-S\,\rchi)^{-1}-I$ (the resolvent operator of $S\,\rchi$), whose kernel we denote by $R(x,y)$, and $Q_{\epsilon}:=(I-S\,\rchi)^{-1}\,\epsilon\,\varphi$. Then \eqref{eq:9} factors into
\begin{equation*}
  A =(I - S\,\rchi)\,B.
\end{equation*}
where $B$ is
\begin{equation*}
  I - Q_{\epsilon}\,\otimes\,\rchi\,\psi - \frac{1}{2}(I+R)\,S\,\epsilon\,\com{\rchi}{D} - \frac{1}{2}\,(Q_{\epsilon}\,\otimes\,\psi)\,\epsilon\,\com{\rchi}{D}
\end{equation*}
Hence
\begin{equation*}
  D_{4,N}(t,\lambda) = D_{2,N}(t,\lambda)\,\det(B).
\end{equation*}
In order to find $\det(B)$ we use the identity
\begin{equation}
  \label{eq:14}
  \epsilon\,\com{\rchi}{D}=\sum_{k=1}^{2m}(-1)^{k}\,\epsilon_{k}\otimes\delta_{k},
\end{equation}
where $\epsilon_{k}$ and $\delta_{k}$ are the functions $\epsilon(x-a_{k})$ and $\delta(x-a_{k})$ respectively, and the $a_{k}$ are the endpoints of the (disjoint) intervals considered, $J=\cup_{k=1}^{m}(a_{2\,k-1},a_{2\,k})$. In our case $m=1$ and $a_{1}=t$, $a_{2}=\infty$. We also make use of the fact that
\begin{equation}
  a\otimes b\cdot c\otimes d= \inprod{b}{c}\cdot a\otimes d
\end{equation}
where $\inprod{.}{.}$ is the usual $L^{2}$--inner product. Therefore
\begin{align*}
  (Q_{\epsilon}\otimes\psi)\,\epsilon\,\com{\rchi}{D} & = \sum_{k=1}^{2} (-1)^{k}Q_{\epsilon}\otimes\psi\cdot\epsilon_{k}\otimes\delta_{k} \\
  & =\sum_{k=1}^{2} (-1)^{k}\inprod{\psi}{\epsilon_{k}}\,Q_{\epsilon}\otimes\,\delta_{k}.
\end{align*}
It follows that
\begin{equation}
  \frac{D_{4,N}(t,\lambda)}{D_{2,N}(t,\lambda)}
\end{equation}
is the determinant of
\begin{equation}
I - Q_{\epsilon}\otimes\rchi\psi - \frac{1}{2}\,\sum_{k=1}^{2} (-1)^{k}\left[(S+R\,S)\,\epsilon_{k} +\inprod{\psi}{\epsilon_{k}}\,Q_{\epsilon}\right]\otimes\delta_{k}.
\end{equation}
We now specialize to the case of one interval $J=(t,\infty)$, so $m=1$, $a_{1}=t$ and $a_{2}=\infty$. We write $\epsilon_{t}=\epsilon_{1}$, and $\epsilon_{\infty}=\epsilon_{2}$, and similarly for $\delta_{k}$. Writing out the terms in the summation and using the fact that
\begin{equation}
  \epsilon_{\infty}=-\frac{1}{2},
\end{equation}
yields
\begin{equation}
  I - Q_{\epsilon}\otimes\rchi\psi  + \frac{1}{2}\,\left[(S+R\,S)\,\epsilon_{t}+ \inprod{\psi}{\epsilon_{t}}\,Q_{\epsilon}\right]\otimes\delta_{t} + \frac{1}{4}\,\left[(S+R\,S)\,1 + \inprod{\psi}{1}\,Q_{\epsilon}\right]\otimes\delta_{\infty}
\end{equation}
Now we can use the formula
\begin{equation}
  \label{eq:16}
  \det\left(I-\sum_{i=1}^{L}\alpha_{i}\otimes\beta_{i}\right)=\det\left(\delta_{jk}-\inprod{\alpha_{j}}{\beta_{k}}\right)_{1\leq j,k\leq L}
\end{equation}
In order to simplify the notation in preparation for the computation of the various inner products, define
\begin{align}
  Q(x,\lambda,t)&:=(I-S\,\rchi)^{-1}\,\varphi, & P(x,\lambda,t)&:=(I-S\,\rchi)^{-1}\,\psi, \nonumber\\
  Q_{\epsilon}(x,\lambda,t)&:=(I-S\,\rchi)^{-1}\,\epsilon\,\varphi, & P_{\epsilon}(x,\lambda,t)&:=(I-S\,\rchi)^{-1}\,\epsilon\,\psi,
\end{align}
\begin{align}
  \label{eq:10}
  q_{_{N}}&:=Q(t,\lambda,t), & p_{_{N}}&:=P(t,\lambda,t),\nonumber\\
  q_{\epsilon}&:=Q_{\epsilon}(t,\lambda,t), & p_{\epsilon}&:=P_{\epsilon}(t,\lambda,t),\nonumber\\
  u_{\epsilon}&:=\inprod{Q}{\rchi\,\epsilon\,\varphi}=\inprod{Q_{\epsilon}}{\rchi\,\varphi}, & v_{\epsilon}&:=\inprod{Q}{\rchi\,\epsilon\,\psi}=\inprod{P_{\epsilon}}{\rchi\,\psi}, \nonumber\\
  \tilde{v}_{\epsilon}&:=\inprod{P}{\rchi\,\epsilon\,\varphi}=\inprod{Q_{\epsilon}}{\rchi\,\varphi}, & w_{\epsilon}&:=\inprod{P}{\rchi\,\epsilon\,\psi}=\inprod{P_{\epsilon}}{\rchi\,\psi},
\end{align}
\begin{equation}
  \label{eq:11}
  \mathcal{P}_{4} := \int_{\mathbb{R}}\,\epsilon_{t}(x)\,P(x,t)\,d\,x, \quad \mathcal{Q}_{4} := \int_{\mathbb{R}}\,\epsilon_{t}(x)\,Q(x,t)\,d\,x, \quad  \mathcal{R}_{4} := \int_{\mathbb{R}}\,\epsilon_{t}(x)\,R(x,t)\,d\,x,
\end{equation}
where we remind the reader that $\epsilon_{t}$ stands for the function $\epsilon(x-t)$. Note that all quantities in \eqref{eq:10} and \eqref{eq:11} are functions of $t$ and $\lambda$ alone. Furthermore, let
\begin{equation}
  c_{\varphi} = \epsilon\,\varphi(\infty)=\frac{1}{2}\int_{-\infty}^{\infty}\varphi(x)\,d\,x, \qquad c_{\psi} = \epsilon\,\psi(\infty)=\frac{1}{2}\int_{-\infty}^{\infty}\psi(x)\,d\,x.
\end{equation}
Recall from the previous section that when $\beta=4$ we take $N$ to be odd. It follows that $\varphi$ and $\psi$ are odd and even functions respectively. Thus when $\beta=4$, $c_{\varphi}=0$ while computation using known integrals for the Hermite polynomials gives
\begin{equation}
c_{\psi}=(\pi N)^{1/4} 2^{-3/4-N/2}\,\frac{(N!)^{1/2}}{(N/2)!}\,\sqrt{\lambda}.
\end{equation}
Hence computation yields
\begin{equation}
 \lim_{N\to\infty}c_{\psi}= \sqrt{\frac{\lambda}{2}}.
\end{equation}
At $t=\infty$,
\begin{equation}
  u_{\epsilon}(\infty)=0, \quad  q_{\epsilon}(\infty)=c_{\varphi}
\end{equation}
\begin{equation}
  \mathcal{P}_{4}(\infty) = -c_{\psi}, \quad \mathcal{Q}_{4}(\infty) = -c_{\varphi},\quad \mathcal{R}_{4}(\infty) = 0.
\end{equation}
In \eqref{eq:16}, $L=3$ and if we denote $a_{4}=\inprod{\psi}{\epsilon_{t}}$, then we have explicitly
\begin{equation*}
  \alpha_{1}=Q_{\epsilon}, \quad \alpha_{2} = -\frac{1}{2}\,\left[(S+R\,S)\,\epsilon_{t} + a_{4}\,Q_{\epsilon}\right],\quad\alpha_{3} = -\frac{1}{4}\,\left[(S+R\,S)\,1 + \inprod{\psi}{1}\,Q_{\epsilon}\right],
\end{equation*}
\begin{equation*}
  \beta_{1}=\rchi\psi, \quad  \beta_{2}=\delta_{t}, \quad  \beta_{3}=\delta_{\infty}.
\end{equation*}
However notice that
\begin{equation}
  \inprod{\left(S+R\,S\right)\,\epsilon_{t}}{\delta_{\infty}}=\inprod{\epsilon_{t}}{\delta_{\infty}}=0,\quad \inprod{\left(S+R\,S\right)\,1}{\delta_{\infty}}=\inprod{1}{R_{\infty}}=0
\end{equation}
and $\inprod{Q_{\epsilon}}{\delta_{\infty}}=c_{\varphi}=0$. Therefore the terms involving $\beta_{3}=\delta_{\infty}$ are all $0$ and we can discard them reducing our computation to that of a $2\times 2$ determinant instead with
\begin{equation}
  \alpha_{1}=Q_{\epsilon}, \quad \alpha_{2} = -\frac{1}{2}\,\left[(S+R\,S)\,\epsilon_{t} + a_{4}\,Q_{\epsilon}\right], \quad \beta_{1}=\rchi\psi, \quad  \beta_{2}=\delta_{t}.
\end{equation}
Hence
\begin{align}
  \inprod{\alpha_{1}}{\beta_{1}} & = \tilde{v}_{\epsilon}, \quad \inprod{\alpha_{1}}{\beta_{2}}=q_{\epsilon}, \\
  \inprod{\alpha_{2}}{\beta_{1}} & =-\frac{1}{2} \,\left(\mathcal{P}_{4} - a_{4}\, + a_{4}\,\tilde{v}_{\epsilon}\right), \\
  \inprod{\alpha_{2}}{\beta_{2}} & = -\frac{1}{2}\,\left(\mathcal{R}_{4} + a_{4}\,q_{\epsilon}\right)\label{eq:12}.
 \end{align}
We want the limit of the determinant
\begin{equation}
  \det\left(\delta_{jk}-\inprod{\alpha_{j}}{\beta_{k}}\right)_{1\leq j,k\leq 2},
\end{equation}
as $N\to \infty$. In order to get our hands on the limits of the individual terms involved in the determinant, we will find differential equations for them first as in \cite{Trac2}. Adding $a_{4}/2$ times row 1 to row 2 shows that ${a}_{4}$ falls out of the determinant, so we will not need to find differential equations for it. Thus our determinant is now
{\large
  \begin{equation}
    \det\left(
      \begin{array}{cc}
        1 - \tilde{v}_{\epsilon} & -q_{\epsilon} \\[6pt]
        \frac{1}{2} \,\mathcal{P}_{4} & 1 + \frac{1}{2}\,\mathcal{R}_{4} \\[6pt]
      \end{array}
    \right).
  \end{equation}
}
Proceeding as in \cite{Trac2} we find the following differential equations
\begin{align}
  \frac{d}{d\,t}\,u_{\epsilon} & = - q_{_{N}}\,q_{\epsilon}, & \frac{d}{d\,t}\,q_{\epsilon} & = q_{_{N}} -q_{_{N}}\, \tilde{v}_{\epsilon} - p_{_{N}}\,u_{\epsilon},\\
 \frac{d}{d\,t}\mathcal{Q}_{4} & = -q_{_{N}}\left(\mathcal{R}_{4} + 1\right), & \frac{d}{d\,t}\mathcal{P}_{4} & = -p_{_{N}}\left(\mathcal{R}_{4} + 1\right),\label{eq:13}\\
 \frac{d}{d\,t}\mathcal{R}_{4} & = -p_{_{N}}\,\mathcal{Q}_{4}-q_{_{N}}\,\mathcal{P}_{4}.   &
\end{align}
Now we change variable from $t$ to $s$ where $t=\tau(s)= \sqrt{2\,N}+s/(\sqrt{2}\,N^{1/6})$ and take the limit $N\to \infty$, denoting the limits of $ q_{\epsilon}$, $\mathcal{P}_{4}$, $\mathcal{Q}_{4}$, $\mathcal{R}_{4},$ and the common limit of $u_{\epsilon}$ and $\tilde{v}_{\epsilon}$ respectively by $\overline{q}$, $\overline{\mathcal{P}}_{4}$, $\overline{\mathcal{Q}}_{4}$, $\overline{\mathcal{R}}_{4}$ and $\overline{u}$. Also $\overline{\mathcal{P}}_{4}$ and $\overline{\mathcal{Q}}_{4}$ differ by a constant, namely $\overline{\mathcal{Q}}_{4}=\overline{\mathcal{P}}_{4} + \sqrt{2}/{2}$. These limits hold uniformly for bounded $s$ so we can interchange $\lim_{N\to\infty}$ and $\frac{d}{d\,s}$. Also $\lim_{N\to\infty}N^{-1/6}q_{_{N}}=\lim_{N\to\infty}N^{-1/6}p_{_{N}}=q $ , where $q$ is as in \eqref{D2}. We obtain the systems
\begin{equation}
  \frac{d}{d\,s}\,\overline{u} = -\frac{1}{\sqrt{2}}\,q\,\overline{q},\qquad \frac{d}{d\,s}\,\overline{q} = \frac{1}{\sqrt{2}}\,q\,\left(1-2\,\overline{u}\right),
\end{equation}
\begin{equation}
  \frac{d}{d\,s}\overline{\mathcal{P}}_{4} = -\frac{1}{\sqrt{2}}\,q\,\left(\overline{\mathcal{R}}_{4} + 1\right), \qquad \frac{d}{d\,s}\overline{\mathcal{R}}_{4} = -\frac{1}{\sqrt{2}}\,q\,\left(2\,\overline{\mathcal{P}}_{4} + \sqrt{\frac{\lambda}{2}}\right),
\end{equation}
The change of variables $s\to\mu=\int_{s}^{\infty} q(x,\lambda)\,d\,x$ transforms these systems into constant coefficient ordinary differential equations
\begin{equation}
  \frac{d}{d\,\mu}\overline{u} = \frac{1}{\sqrt{2}}\,\overline{q}, \qquad \frac{d}{d\,\mu}\overline{q} = -\frac{1}{\sqrt{2}}\,\left(1-2\,\overline{u}\right),
\end{equation}
\begin{equation}
  \frac{d}{d\,\mu}\overline{\mathcal{P}}_{4} = \frac{1}{\sqrt{2}}\,\left(\overline{\mathcal{R}}_{4} + 1\right), \qquad \frac{d}{d\,\mu}\overline{\mathcal{R}}_{4} = \frac{1}{\sqrt{2}}\,\left(2\,\overline{\mathcal{P}}_{4}+ \sqrt{\frac{\lambda}{2}}\right).
\end{equation}
Since $\lim_{s\to \infty}\mu=0$, corresponding to the boundary values at $t=\infty$ which we found earlier for $\mathcal{P}_{4}, \mathcal{R}_{4}$, we now have initial values at $\mu=0$. Therefore
\begin{equation}
  \overline{u}(\mu=0)=\overline{q}(\mu=0)=0,
\end{equation}
\begin{equation}
  \overline{\mathcal{P}}_{4}(\mu=0)=-\sqrt{\frac{\lambda}{2}}, \qquad \overline{\mathcal{R}}_{4}(\mu=0)=0.
\end{equation}
We use this to solve the systems and get
\begin{align}
  \overline{q} & = \frac{1}{2\sqrt{2}}\,\left(e^{-\mu}-e^{\mu}\right)\\
 \overline{u} & = \frac{1}{2}\,\left(1 - \frac{1}{2}\,e^{\mu} - \frac{1}{2}\,e^{-\mu}\right)\\
 \overline{\mathcal{P}}_{4} & = \frac{1}{2\,\sqrt{2}}\,\left(\frac{2-\sqrt{\lambda}}{2}\,e^{\mu} - \frac{2+\sqrt{\lambda}}{2}\,e^{-\mu} - \sqrt{\lambda}\right) \label{eq:24}\\
 \overline{\mathcal{R}}_{4} & = \frac{2-\sqrt{\lambda}}{4}\,e^{\mu} + \frac{2+\sqrt{\lambda}}{4}\,e^{-\mu} - 1\label{eq:25}
\end{align}
Substituting these expressions into the determinant gives
\eqref{gsedet}, namely
\begin{equation}
  \label{eq:26}
  D_{4}(s,\lambda)= D_{2}(s,\lambda)\,\cosh^{2}\left(\frac{\mu(s,\lambda)}{2}\right),
\end{equation}
where $D_{\beta}=\lim_{N\to\infty}D_{\beta,N}$.
Note that even though there are $\lambda$--terms in \eqref{eq:24} and \eqref{eq:25}, these do not appear in the final result \eqref{eq:26}, making it similar to the GUE case where the main conceptual difference between the $m=1$ (largest eigenvalue) case and the general $m$ is the dependence of the function $q$ on $\lambda$. The right hand side of the above formula clearly reduces to the $\beta=4$ Tracy-Widom distribution when we set $\lambda=1$.
Note that where we have $D_{4}(s,\lambda)$ above, Tracy and Widom (and hence many RMT references) write $D_{4}(s/\sqrt{2},\lambda)$ instead. Tracy and Widom applied the change of variable $s\to s/\sqrt{2}$ in their derivation in \cite{Trac2} so as to agree with Mehta's form of the $\beta=4$ joint eigenvalue density (see \eqref{eq:37} and footnote in the Introduction), which has $-2x^{2}$ in the exponential in the weight function, instead of $-x^{2}$ in our case. To switch back to the other convention, one just needs to substitute in the argument $s/\sqrt{2}$ for $s$ everywhere in our results. At this point this is just a cosmetic discrepancy, and it does not change anything in our derivations since all the differentiations are done with respect to $\lambda$ anyway. It \underline{does} change conventions for rescaling data while doing numerical work though (see Section~\ref{sec:stand-devi-matt}).

\chapter{The distribution of the $m^{th}$ largest eigenvalue in the GOE}

\chaptermark{The $m^{th}$ largest eigenvalue in the GOE}

\section{The distribution function as a Fredholm determinant}
The GOE corresponds case corresponds to the specialization $\beta=1$ in \eqref{jointdensity} so that
\begin{equation}
\label{eq:17}
  G_{1,N}(t,\lambda)=C_{1}^{(N)}\underset{x_{i}\in\mathbb{R}}{\int\cdots\int} \prod_{j<k}\left|x_{j}-x_{k}\right|\,\prod_{j}^{N}w(x_{j})\,\prod_{j}^{N}\left(1+f(x_{j})\right)\,d\,x_{1}\cdots d\,x_{N}
\end{equation}
where $w(x)=\exp\left(-x^{2}\right)$, $f(x)=-\lambda\,\rchi_{J}(x)$, and $C_{1}^{(N)}$ depends only on $N$. As in the GSE case, we will lump into $C_{1}^{(N)}$ any constants depending only on $N$ that appear in the derivation. A simple argument  at the end will show that the final constant is $1$. These calculations more or less faithfully follow and expand on \cite{Trac1}. We want to use \eqref{debruijn2}, which requires an ordered space. Note that the above integrand is symmetric under permutations, so the integral is $n!$ times the same integral over ordered pairs $x_{1}\leq\ldots\leq x_{N}$. So we can rewrite \ref{eq:17} as
\begin{equation*}
  (N!)\,\underset{x_{1}\leq\ldots\leq x_{N}\in \mathbb R\,\,}{\int\cdots\int} \prod_{j<k}(x_{k}-x_{j})\,\prod_{i=1}^{N}w(x_{k})\prod_{i=1}^{N}(1+f(x_{k}))\,d\,x_{1}\cdots d\,x_{N},
\end{equation*}
where we can remove the absolute values since the ordering insures that $(x_{j}-x_{i})\geq 0$ for $i<j$. Recall that the Vandermonde determinant is
\begin{equation*}
  \Delta_{N}(x)=\det(x^{j-1}_{k})_{_{1\leq j,k\leq N}}=(-1)^{\frac{N\,(N-1)}{2}}\prod_{j<k}(x_{j}-x_{k}).
\end{equation*}
Therefore what we have inside the integrand above is, up to sign
\begin{equation*}
  \det(x^{j-1}_{k}\,w(x_{k})\,(1+f(x_{k})))_{_{1\leq j,k\leq N}}.
\end{equation*}
Note that the sign depends only on $N$. Now we can use \eqref{debruijn2} with $\varphi_{j}(x)=x^{j-1}\,w(x)\,(1+f(x))$. In using \eqref{debruijn2} we square both sides so that the right hand side is now a determinant instead of a Pfaffian. Therefore $G_{1,N}^{2}(t,\lambda)$ equals
\begin{equation*}
  C^{(N)}_{1}\,\det\left(\int\int\sgn(x-y)x^{j-1}\,y^{k-1}\,(1+f(x))\,w(x)\,w(y)\,d\,x\,d\,y \right)_{_{1\leq j,k\leq N}}.
\end{equation*}
Shifting indices, we can write it as
\begin{equation}
  C^{(N)}_{1}\,\det\left(\int\int\sgn(x-y)x^{j}\,y^{k}\,(1+f(x))\,w(x)\,w(y)\,d\,x\,d\,y \right)_{_{1\leq j,k\leq N-1}}
\end{equation}
where $C^{(N)}_{1}$ is a constant depending only on $N$, and is such that the right side is $1$ if $f\equiv 0$. Indeed this would correspond to the probability that $\lambda_{p}^{GOE(N)}<\infty$, or equivalently to the case where the excluded set $J$ is empty. We can replace $x^{j}$ and $y^{k}$ by any arbitrary polynomials $p_{j}(x)$ and $p_{k}(x)$, of degree $j$ and $k$ respectively, which are obtained by row operations on the matrix. Indeed such operations would not change the determinant. We also replace $\sgn(x-y)$ by $\epsilon(x-y)=\frac{1}{2}\,\sgn(x-y)$ which just produces a factor of $2$ that we absorb in $C^{(N)}_{1}$. Thus $G_{1,N}^{2}(t,\lambda)$ now equals
\begin{equation}
 C^{(N)}_{1}\,\det\left(\int\int\epsilon(x-y)\,p_{j}(x)\,p_{k}(y)\,(1 + f(x))\,(1 + f(y))\,w(x)\,w(y)\,d\,x\,d\,y \right)_{_{0\leq j,k\leq N-1}}.
\end{equation}
Let $\psi_{j}(x)=p_{j}(x)\,w(x)$  so the above integral becomes
\begin{equation}
 C^{(N)}_{1}\,\det\left(\int\int\epsilon(x-y)\,\psi_{j}(x)\,\psi_{k}(y)\,(1+f(x)+f(y)+f(x)\,f(y))\,d\,x\,d\,y \right)_{_{0\leq j,k\leq N-1}}.
\end{equation}
Partially multiplying out the term we obtain
\begin{align*}
  C^{(N)}_{1}\,&\det\left(\int\int\epsilon(x-y)\,\psi_{j}(x)\,\psi_{k}(y)\,d\,x\,d\,y \right. \\
  & + \left. \int\int\epsilon(x-y)\,\psi_{j}(x)\,\psi_{k}(y)\,(f(x)+f(y)+f(x)\,f(y))\,d\,x\,d\,y \right)_{_{0\leq j,k\leq N-1}}.
\end{align*}
Define
\begin{equation}
  \label{eq:18}
  M=\left(\int\int\epsilon(x-y)\,\psi_{j}(x)\,\psi_{k}(y)\,d\,x\,d\,y \right)_{_{0\leq j,k\leq N-1}},
\end{equation}
so that $G_{1,N}^{2}(t,\lambda)$ is now
\begin{equation*}
  C^{(N)}_{1}\,\det\left(M + \int\int\epsilon(x-y)\,\psi_{j}(x)\,\psi_{k}(y)\,(f(x)+f(y)+f(x)\,f(y))\,d\,x\,d\,y \right)_{_{0\leq j,k\leq N-1}}.
\end{equation*}
Let $\epsilon$ be the operator defined in \eqref{epsilonop}. We can use operator notation to simplify the expression for $G_{1,N}^{2}(t,\lambda)$ a great deal by rewriting the double integrals as single integrals. Indeed
\begin{eqnarray*}
  \int\int\epsilon(x-y)\,\psi_{j}(x)\,\psi_{k}(y)\,f(x)\,d\,x\,d\,y & = & \int f(x)\,\psi_{j}(x)\,\int\epsilon(x-y)\,\psi_{k}(y)\,d\,y\,d\,x\\
  & = & \int f\,\psi_{j}\,\epsilon\,\psi_{k}\,d\,x.
\end{eqnarray*}
Similarly,
\begin{eqnarray*}
  \int\int\epsilon(x-y)\,\psi_{j}(x)\,\psi_{k}(y)\,f(y)\,d\,x\,d\,y & = & - \int\int\epsilon(y-x)\,\psi_{j}(x)\,\psi_{k}(y)\,f(y)\,d\,x\,d\,y \\
  & = & - \int f(y)\,\psi_{k}(y)\,\int\epsilon(y-x)\,\psi_{j}(x)\,d\,x\,d\,y \\
  & = &  - \int f(x)\,\psi_{k}(x)\,\int\epsilon(x-y)\,\psi_{j}(y)\,d\,y\,d\,x \\
  & = &  - \int f\,\psi_{k}\,\epsilon\,\psi_{j}\,d\,x.
\end{eqnarray*}
Finally,
\begin{equation}
  \begin{aligned}
    \int\int \epsilon(x-y)\,\psi_{j}(x)\,&\psi_{k}(y)\,f(x)\,f(y)\,d\,x\,d\,y  \\
    &  = - \int\int\epsilon(y-x)\,\psi_{j}(x)\,\psi_{k}(y)\,f(x)\,f(y)\,d\,x\,d\,y \\
    &  = - \int f(y)\,\psi_{k}(y)\,\int\epsilon(y-x)\,f(x)\,\psi_{j}(x)\,d\,x\,d\,y \\
    &  = - \int f(x)\,\psi_{k}(x)\,\int\epsilon(x-y)\,f(y)\,\psi_{j}(y)\,d\,y\,d\,x \\
    &  = - \int f\,\psi_{k}\,\epsilon\,(f\,\psi_{j})\,d\,x.
  \end{aligned}
\end{equation}
It follows that
\begin{equation}
  G_{1,N}^{2}(t,\lambda)=C^{(N)}_{1}\,\det\left(M + \int\left[ f\,\psi_{j}\,\epsilon\,\psi_{k} - f\,\psi_{k}\,\epsilon\,\psi_{j} -  f\,\psi_{k}\,\epsilon\,(f\,\psi_{j}) \right]\,d\,x \right)_{_{0\leq j,k\leq N-1}}.
\end{equation}
If we let $M^{-1}=\left(\mu_{j\,k} \right)_{_{0\leq j,k\leq N-1}}$, and factor $\det(M)$ out, then $G_{1,N}^{2}(t,\lambda) $ equals
\begin{align}
  C^{(N)}_{1}\,& \det(M)\,\det\left(\vphantom{M^{-1}\cdot\left(\int\left[ f\,\psi_{j}\,\epsilon\,\psi_{k} - f\,\psi_{k}\,\epsilon\,\psi_{j} -  f\,\psi_{k}\,\epsilon\,(f\,\psi_{j}) \right]\,d\,x\right)_{_{0\leq j,k\leq N-1}}} I + \right. \nonumber\\ & \left. M^{-1}\cdot\left(\int\left[ f\,\psi_{j}\,\epsilon\,\psi_{k} - f\,\psi_{k}\,\epsilon\,\psi_{j} -  f\,\psi_{k}\,\epsilon\,(f\,\psi_{j}) \right]\,d\,x\right)_{_{0\leq j,k\leq N-1}} \right)_{_{0\leq j,k\leq N-1}}
\end{align}
where the dot denotes matrix multiplication of $M^{-1}$ and the matrix with the integral as its $(j, k)$--entry. define $\eta_{j}=\sum_{j}\mu_{j\,k}\,\psi_{k}$ and use it to simplify the result of carrying out the matrix multiplication. From \eqref{eq:18} it follows that $\det(M)$ depends only on $N$ we lump it into $C^{(N)}_{1}$. Thus $G_{1,N}^{2}(t,\lambda)$ equals
\begin{equation}
  C^{(N)}_{1}\,\det\left(I + \left(\int\left[ f\,\eta_{j}\,\epsilon\,\psi_{k} - f\,\psi_{k}\,\epsilon\,\eta_{j} -  f\,\psi_{k}\,\epsilon\,(f\,\eta_{j}) \right]\,d\,x\right)_{_{0\leq j,k\leq N-1}} \right)_{_{0\leq j,k\leq N-1}}.
\end{equation}
Recall our remark at the very beginning of the section that if $f\equiv 0$ then the integral we started with evaluates to $1$ so that
\begin{equation}
  C^{(N)}_{1}\det(I) = C^{(N)}_{1},
\end{equation}
which implies that $C^{(N)}_{1}=1$. Now $G_{1,N}^{2}(t,\lambda)$ is of the form $\det(I+AB)$ where $A:L^{2}(J)\times L^{2}(J) \to \mathbb C^{N}$ is a $N\times 2$ matrix
\begin{equation*}
  A=\left(
    \begin{array}{c}
      A_{1} \\
      A_{2} \\
      \vdots \\
      A_{N}
    \end{array}
  \right),
\end{equation*}
whose $j^{th}$ row is given by
\begin{equation*}
  A_{j}= A_{j}(x)=\left(-f\,\epsilon\,\eta_{j}-f\,\epsilon\,\left(f\,\eta_{j}\right)\qquad f\,\eta_{j}\right).
\end{equation*}
Therefore, if
\begin{equation*}
  g=\left(
    \begin{array}{c}
      g_{1} \\
      g_{2}
    \end{array}
  \right)\in L^{2}(J)\times L^{2}(J),
\end{equation*}
then $A\,g$ is a column vector whose $j^{th}$ row is $\inprod{A_{j}}{g}_{_{L^{2}\times L^{2}}}$
\begin{equation*}
  \left(A\,g\right)_{j}=  \int\left[-f\,\epsilon\,\eta_{j}-f\,\epsilon\,\left(f\,\eta_{j}\right)\right]\,g_{1}\,d\,x + \int f\,\eta_{j}\,g_{2}\,d\,x.
\end{equation*}
Similarly, $B:\mathbb C^{N} \to L^{2}(J)\times L^{2}(J)$ is a $2\times N$ matrix
\begin{equation*}
  B=\left(
    \begin{array}{cccc}
      B_{1} & B_{2} & \ldots & B_{N}
    \end{array}\right),
\end{equation*}
whose $j^{th}$ column is given by
\begin{equation*}
  B_{j}=B_{j}(x)=\left(\begin{array}{c} \psi_{j} \\ \epsilon\,\psi_{j}
    \end{array}
  \right).
\end{equation*}
Thus if
\begin{equation*}
  h=\left(
    \begin{array}{c}
      h_{1}  \\
      \vdots \\
      h_{N}
    \end{array}
  \right)\in \mathbb C^{N},
\end{equation*}
then $B\,h$ is the column vector of $L^{2}(J)\times L^{2}(J)$ given by
\begin{equation*}
  B\,h=\left(
    \begin{array}{c}
      \sum_{j}h_{i}\,\psi_{j}  \\
      \sum_{j}h_{i}\,\epsilon\,\psi_{j}
    \end{array}
  \right).
\end{equation*}
Clearly $A\,B:\mathbb C^{N} \to \mathbb C^{N}$ and $B\,A: L^{2}(J)\times L^{2}(J) \to L^{2}(J)\times L^{2}(J) $ with kernel
\begin{equation*}
  \left(
    \begin{array}{cc}
      -\sum_{j}\psi_{j}\otimes f\,\epsilon\,\eta_{j} - \sum_{j}\psi_{j}\otimes f\,\epsilon\,(f\,\eta_{j}) & \sum_{j}\psi_{j}\otimes f\,\eta_{j} \\
      -\sum_{j}\epsilon\,\psi_{j}\otimes f\,\epsilon\,\eta_{j} - \sum_{j}\epsilon\,\psi_{j}\otimes f\,\epsilon\,(f\,\eta_{j}) &
      \sum_{j}\epsilon\,\psi_{j}\otimes f\,\eta_{j}
    \end{array}
  \right).
\end{equation*}
Hence $I+BA$ has kernel
\begin{equation*}
  \left(
    \begin{array}{cc}
      I -\sum_{j}\psi_{j}\otimes f\,\epsilon\,\eta_{j} - \sum_{j}\psi_{j}\otimes f\,\epsilon\,(f\,\eta_{j}) &  \sum_{j}\psi_{j}\otimes f\,\eta_{j} \\
      -\sum_{j}\epsilon\,\psi_{j}\otimes f\,\epsilon\,\eta_{j} - \sum_{j}\epsilon\,\psi_{j}\otimes f\,\epsilon\,(f\,\eta_{j}) & I + \sum_{j}\epsilon\,\psi_{j}\otimes f\,\eta_{j}
    \end{array}
  \right),
\end{equation*}
which can be written as
\begin{equation*}
  \left(
    \begin{array}{cc}
      I -\sum_{j}\psi_{j}\otimes f\,\epsilon\,\eta_{j} &   \sum_{j}\psi_{j}\otimes f\,\eta_{j} \\
      -\sum_{j}\epsilon\,\psi_{j}\otimes f\,\epsilon\,\eta_{j} - \epsilon\,f & I + \sum_{j}\epsilon\,\psi_{j}\otimes f\,\eta_{j}
    \end{array}
  \right)\cdot \left(
    \begin{array}{cc}
      I & 0 \\
      \epsilon\,f & I
    \end{array}
  \right).
\end{equation*}
Since we are taking the determinant of this operator expression, and the determinant of the second term is just 1, we can drop it. Therefore
\begin{align*}
  G_{1,N}^{2}(t,\lambda) & =\det\left(
    \begin{array}{cc}
      I -\sum_{j}\psi_{j}\otimes f\,\epsilon\,\eta_{j} &  \sum_{j}\psi_{j}\otimes f\,\eta_{j} \\
      -\sum_{j}\epsilon\,\psi_{j}\otimes f\,\epsilon\,\eta_{j} - \epsilon\,f & I + \sum_{j}\epsilon\,\psi_{j}\otimes f\,\eta_{j}
    \end{array}
  \right) \\
  &= \det(I+\,K_{1,N}\,f),
\end{align*}
where
\begin{align*}
  K_{1,N}& = \left(
    \begin{array}{cc}
      -\sum_{j}\psi_{j}\otimes \,\epsilon\,\eta_{j}  &   \sum_{j}\psi_{j}\otimes \,\eta_{j} \\
      -\sum_{j}\epsilon\,\psi_{j}\otimes \,\epsilon\,\eta_{j} - \epsilon & \sum_{j}\epsilon\,\psi_{j}\otimes \,\eta_{j}
    \end{array}
  \right)\\
  &= \left(
    \begin{array}{cc}
      -\sum_{j,k}\psi_{j}\otimes \,\mu_{jk}\,\epsilon\,\psi_{k}  &   \sum_{j,k}\psi_{j}\otimes \,\mu_{jk}\,\psi_{k} \\
      -\sum_{j,k}\epsilon\,\psi_{j}\otimes \,\mu_{jk}\,\epsilon\,\psi_{k} - \epsilon & \sum_{j,k}\epsilon\,\psi_{j}\otimes \,\mu_{jk}\,\psi_{k}
    \end{array}
  \right)
\end{align*}
and $K_{1,N}$ has matrix kernel
\begin{align*}
  K_{1,N}(x,y)&= \left(
    \begin{array}{cc}
      -\sum_{j,k}\psi_{j}(x)\,\mu_{jk}\,\epsilon\,\psi_{k}(y) &   \sum_{j,k}\psi_{j}(x)\,\mu_{jk}\,\psi_{k}(y) \\
      -\sum_{j,k}\epsilon\,\psi_{j}(x) \,\mu_{jk}\,\epsilon\,\psi_{k}(y) - \epsilon(x-y) & \sum_{j,k}\epsilon\,\psi_{j}(x)\,\mu_{jk}\,\psi_{k}(y)
    \end{array}
  \right).
\end{align*}
We define
\begin{equation*}
  S_{N}(x,y)=- \sum_{j,k=0}^{N-1}\psi_{j}(x)\mu_{jk}\epsilon\,\psi_{k}(y).
\end{equation*}
Since $M$ is antisymmetric,
\begin{equation*}
  S_{N}(y,x)=- \sum_{j,k=0}^{N-1}\psi_{j}(y)\mu_{jk}\epsilon\,\psi_{k}(x)=\sum_{j,k=0}^{N-1}\psi_{j}(y)\mu_{kj}\epsilon\,\psi_{k}(x)= \sum_{j,k=0}^{N-1}\epsilon\,\psi_{j}(x)\mu_{jk}\psi_{k}(y).
\end{equation*}
Note that
\begin{equation*}
  \epsilon\,S_{N}(x,y)= \sum_{j,k=0}^{N-1}\epsilon\,\psi_{j}(x)\,\mu_{jk}\,\epsilon\,\psi_{k}(y),
\end{equation*}
whereas
\begin{equation*}
  -\frac{d}{d\,y}\,S_{N}(x,y)=\sum_{j,k=0}^{N-1}\psi_{j}(x)\mu_{jk}\psi_{k}(y).
\end{equation*}
So we can now write succinctly
\begin{equation}
  \label{eq:19}
  K_{1,N}(x,y)=\left(
    \begin{array}{cc}
      S_{N}(x,y) & -\frac{d}{d\,y}\,S_{N}(x,y) \\
      \epsilon\,S_{N}(x,y) - \epsilon & S_{N}(y,x)
    \end{array}
  \right)
\end{equation}
So we have shown that
\begin{equation}
\label{eq:21}
   G_{1,N}(t,\lambda)=\sqrt{D_{1,N}\left(t,\lambda\right)}
\end{equation}
where
\begin{equation*}
   D_{1,N}(t,\lambda)=\det\left(I+K_{1,N}\,f\right)
\end{equation*}
where $K_{1,N}$ is the integral operator with matrix kernel $K_{1,N}(x,y)$ given in \eqref{eq:19}.

\section{Gaussian specialization}
\label{sec:gauss-spec}

We specialize the results above to the case of a Gaussian weight function
\begin{equation}
  \label{eq:36}
  w(x)=\exp\left(-x^{2}/2\right)
\end{equation}
and indicator function
\begin{equation*}
  f(x)=-\lambda\,\rchi_{_{J}}\qquad J=(t,\infty)
\end{equation*}
Note that this does not agree with the weight function in \eqref{jointdensity}. However it is a necessary choice if we want the technical convenience of working with exactly the same orthogonal polynomials (the Hermite functions) as in the $\beta=2,4$ cases. In turn the Painlev\'e function in the limiting distribution will be unchanged. The discrepancy is resolved by the choice of standard deviation. Namely here the standard deviation on the diagonal matrix elements is taken to be $1$, corresponding to the weight function \eqref{eq:36}. In the $\beta=2,4$ cases the standard deviation on the diagonal matrix elements is $1/\sqrt{2}$, giving the weight function \eqref{gseweight}. We expand on this in Section~\ref{sec:stand-devi-matt}.
Now we again want the matrix
\begin{equation*}
  M=\left(\int\int\epsilon(x-y)\,\psi_{j}(x)\,\psi_{k}(y)\,d\,x\,d\,y \right)_{_{0\leq j,k\leq N-1}} = \left(\int\,\psi_{j}(x)\,\epsilon\,\psi_{k}(x)\,d\,x \right)_{_{0\leq j,k\leq N-1}}
\end{equation*}
to be the direct sum of $\frac{N}{2}$ copies of
\begin{equation*}
  Z=\left(\begin{array}{cc} 0 & 1 \\ -1 & 0 \end{array}\right)
\end{equation*}
so that the formulas are the simplest possible, since then $\mu_{jk}$ can only be $0$ or $\pm 1$.  In that case $M$ would be skew--symmetric so that $M^{-1}=-M$. In terms of the integrals defining the entries of $M$ this means that we would like to have
\begin{equation*}
  \int\psi_{2m}(x)\,\epsilon\,\psi_{2n+1}(x)\,d\,x = \delta_{m,n},
\end{equation*}
\begin{equation*}
  \int \psi_{2m+1}(x)\,\epsilon\,\psi_{2n}(x)\,d\,x = -\delta_{m,n},
\end{equation*}
and otherwise
\begin{equation*}
  \int \psi_{j}(x)\,\frac{d}{d\,x}\,\psi_{k}(x) \,d\,x = 0.
\end{equation*}
It is easier to treat this last case if we replace it with three non-exclusive conditions
\begin{equation*}
  \int \psi_{2m}(x)\,\epsilon\,\psi_{2n}(x)d\,x = 0,
\end{equation*}
\begin{equation*}
  \int \psi_{2m+1}(x)\,\epsilon\,\psi_{2n+1}(x)\,d\,x = 0
\end{equation*}
(so when the parity is the same for  $j,k$, which takes care of diagonal entries, among others), and
\begin{equation*}
  \int \psi_{j}(x)\,\epsilon\,\psi_{k}(x)\,d\,x = 0.
\end{equation*}
whenever $|j-k|>1$, which targets entries outside of the tridiagonal. Define
\begin{equation*}
  \varphi_{n}(x)=\frac{1}{c_{n}}H_{n}(x)\,\exp(-x^{2}/2)\quad\textrm{for}\quad c_{n}=\sqrt{2^{n}n!\sqrt{\pi}}
\end{equation*}
where the $H_{n}$ are the usual Hermite polynomials defined by the orthogonality condition
\begin{equation*}
  \int H_{j}(x)\,H_{k}(x)\,e^{-x^{2}}\,d\,x = c_{j}^{2}\,\delta_{j,k}.
\end{equation*}
It follows that
\begin{equation*}
  \int \varphi_{j}(x)\,\varphi_{k}(x)\,d\,x = \delta_{j,k}.
\end{equation*}
Now let
\begin{equation}
  \label{eq:20}
  \psi_{_{2n+1}}(x) = \frac{d}{d\,x}\,\varphi_{_{2n}}\left(\,x\right)\qquad
  \psi_{_{2n}}(x) = \varphi_{_{2n}}\left(\,x\right).
\end{equation}
This definition satisfies our earlier requirement that $\psi_{j}=p_{j}\,w$ for
\begin{equation*}
  w(x)=\exp\left(-x^{2}/2\right).
\end{equation*}
In this case for example
\begin{equation*}
  p_{2n}(x)=\frac{1}{c_{n}}\,H_{2n}\left(x\right).
\end{equation*}
With $\epsilon$ defined as in \eqref{epsilonop}, and recalling that, if $D$ denote the operator that acts by differentiation with respect to $x$, then $D\,\epsilon=\epsilon\,D=I$, it follows that
\begin{equation*}
  \begin{aligned}
    \int \psi_{_{2m}}(x)\epsilon\,\psi_{_{2n+1}}(x)d\,x & =\int \varphi_{_{2m}}\left(x\right)\epsilon\,\frac{d}{d\,x}\, \varphi_{_{2n+1}}\left(x\right)d\,x \\
    &=\int \varphi_{_{2m}}\left(x\right) \varphi_{_{2n+1}}\left(x\right)d\,x \\
    &= \int \varphi_{_{2m}}\left(x\right) \varphi_{_{2n+1}}\left(x\right)d\left(x\right) \\
    &=\delta_{m,n},
  \end{aligned}
\end{equation*}
as desired. Similarly, integration by parts gives
\begin{equation*}
  \begin{aligned}
    \int \psi_{_{2m+1}}(x)\epsilon\,\psi_{_{2n}}(x)d\,x & =\int \frac{d}{d\,x}\,\varphi_{_{2m}}\left(x\right)\,\epsilon\,\varphi_{_{2n}} \left(x\right)d\,x \\
    &= -\int \,\varphi_{_{2m}}\left(x\right)\,\varphi_{_{2n}} \left(x\right)d\,x \\
    &=-\int \varphi_{_{2m}}\left(x\right) \varphi_{_{2n+1}}\left(x\right)d\left(x\right) \\& = -\delta_{m,n}.
  \end{aligned}
\end{equation*}
Also $\psi_{_{2n}}$ is even since $H_{2n}$ and $\varphi_{_{2n}}$ are. Similarly, $\psi_{_{2n+1}}$ is odd.  It follows  that $\epsilon\,\psi_{_{2n}}$, and $\epsilon\,\psi_{_{2n+1}}$, are respectively odd and even functions. From these observations, we obtain
\begin{equation*}
  \begin{aligned}\int \psi_{_{2n}}(x)\,\epsilon\,\psi_{_{2m}}(x)d\,x  = 0,
  \end{aligned}
\end{equation*}
since the integrand is a product of an odd and an even function. Similarly
\begin{equation*}
  \begin{aligned}
    \int \psi_{_{2n+1}}(x)\,\epsilon\,\psi_{_{2m+1}}(x)d\,x= 0.
  \end{aligned}
\end{equation*}
Finally it is easy to see that if $|j-k|>1$, then
\begin{equation*}
  \int\psi_{_{j}}(x)\,\epsilon\,\psi_{_{k}}(x)d\,x=0.
\end{equation*}
Indeed both  differentiation and the action of $\epsilon$ can only ``shift'' the indices by $1$. Thus by orthogonality of the $\varphi_{j}$, this integral will always be $0$. Thus by our choice in \eqref{eq:20}, we force the matrix
\begin{equation*}
M= \left(\int\,\psi_{j}(x)\,\epsilon\,\psi_{k}(x)\,d\,x \right)_{_{0\leq j,k\leq N-1}}
\end{equation*}
to be the direct sum of $\frac{N}{2}$ copies of
\begin{equation*}
  Z=\left(
    \begin{array}{cc}
      0 & 1 \\
      -1 & 0
    \end{array}
  \right)
\end{equation*}
This means $M^{-1}=-M$ where $M^{-1}=\{\mu_{j,k}\}$. Moreover, $\mu_{j,k}=0$ if $j,k$ have the same parity or $|j-k|>1$, and $\mu_{2j,2k+1}=\delta_{jk}=-\mu_{2k+1,2j}$ for $j,k=0.\ldots,\frac{N}{2}-1$. Therefore
\begin{equation*}
  \begin{aligned}S_{N}(x,y) & =- \sum_{j,k=0}^{N-1}\psi_{j}(x)\mu_{jk}\epsilon\,\psi_{k}(y)\\
    & = -\sum_{j=0}^{N/2-1}\,\psi_{2j}(x)\,\epsilon\,\psi_{2j+1}(y) + \sum_{j=0}^{N/2-1}\,\psi_{2j+1}(x)\,\epsilon\,\psi_{2j}(y) \\
    & =  \left[ \sum_{j=0}^{N/2-1}\varphi_{_{2j}}\left(\frac{x}{\,\sigma}\right) \,\varphi_{_{2j}}\left(\frac{y}{\,\sigma}\right)-\sum_{j=0}^{N/2-1}\frac{d}{d\,x}\, \varphi_{_{2j}}\left(\frac{x}{\,\sigma}\right)\,\epsilon\,\varphi_{_{2j}} \left(\frac{y}{\,\sigma}\right) \right].
  \end{aligned}
\end{equation*}
Manipulations similar to those in the $\beta=4$ case (see \eqref{hermrec1} through \eqref{newgseS}) yield
\begin{equation*}
  \begin{aligned}
    S_{N}(x,y) =  \left[ \sum_{j=0}^{N-1}\varphi_{_{j}}\left(x\right) \varphi_{_{j}}\left(y\right)-\sqrt{\frac{N}{2}}\, \varphi_{_{N-1}}\left(x\right)\left(\epsilon\,\varphi_{_{N}}\right) \left(y\right) \right].
  \end{aligned}
\end{equation*}
We redefine
\begin{equation*}
  S_{N}(x,y) = \sum_{j=0}^{N-1}\varphi_{_{j}}\left(x\right) \varphi_{_{j}}\left(y\right)=S_{N}(y,x),
\end{equation*}
so that the top left entry of $K_{1,N}(x,y)$ is
\begin{equation*}
  S_{N}(x,y) +  \sqrt{\frac{N}{2}}\, \varphi_{_{N-1}}\left(x\right)\left(\epsilon\,\varphi_{_{N}}\right) \left(y\right).
\end{equation*}
If $S_{N}$ is the operator with kernel $S_{N}(x,y)$ then integration by parts gives
\begin{equation*}
  S_{N}D f = \int S(x,y) \frac{d}{d\,y}f(y)\,d\,y = \int  \left(-\frac{d}{d\,y} S_{N}(x,y)\right) f(y)\,d\,y,
\end{equation*}
so that $-\frac{d}{d\,y} S_{N}(x,y)$ is in fact the kernel of $S_{N}D$. Therefore \eqref{eq:21} now holds with $K_{1,N}$ being the integral operator with matrix kernel $K_{1,N}(x,y)$ whose $(i,j)$--entry $K_{1,N}^{(i,j)}(x,y)$ is given by
\begin{equation*}
  \begin{aligned}
    K _{1,N}^{(1,1)}(x,y) & =\left[ S_{N}(x,y) + \sqrt{\frac{N}{2}}\, \varphi_{_{N-1}}\left(x\right)\left(\epsilon\,\varphi_{_{N}}\right) \left(y\right)\right], \\
    K_{1,N}^{(1,2)}(x,y) & =\left[ SD_{N}(x,y) -\frac{d}{d\,y} \, \left(\sqrt{\frac{N}{2}}\, \varphi_{_{N-1}}\left(x\right)\left(\epsilon\,\varphi_{_{N}}\right) \left(y\right) \right) \right], \\
    K_{1,N}^{(2,1)}(x,y) & = \epsilon\,\left[ S_{N}(x,y) +  \sqrt{\frac{N}{2}}\, \varphi_{_{N-1}}\left(x\right)\left(\epsilon\,\varphi_{_{N}}\right)\left(y\right)-1 \right], \\
    K_{1,N}^{(2,2)}(x,y) & =\left[ S_{N}(x,y) + \sqrt{\frac{N}{2}}\,\left(\epsilon\,\varphi_{_{N}}\right)\left(x\right)\, \varphi_{_{N-1}} \left(y\right)\right].
  \end{aligned}
\end{equation*}
Define
\begin{equation*}
  \varphi(x)=\left(\frac{N}{2}\right)^{1/4}\varphi_{_{N}}(x),\qquad \psi(x)=\left(\frac{N}{2}\right)^{1/4}\varphi_{_{N-1}}(x),
\end{equation*}
so that
\begin{equation*}
  \begin{aligned}
    K_{1,N}^{(1,1)}(x,y)& =\rchi(x)\,\left[ S_{N}(x,y) + \psi(x)\,\epsilon\,\varphi(y)\right] \,\rchi(y),\\
    K_{1,N}^{(1,2)}(x,y) & = \rchi(x)\,\left[ SD_{N}(x,y) - \psi(x)\,\varphi(y)  \right]\,\rchi(y), \\
    K_{1,N}^{(2,1)}(x,y) & = \rchi(x)\,\left[ \epsilon S_{N}(x,y) + \epsilon\,\psi(x)\,\epsilon\,\varphi(y) - \epsilon(x-y) \right]\,\rchi(y), \\
    K_{1,N}^{(2,2)}(x,y) & = \rchi(x)\,\left[ S_{N}(x,y) + \epsilon\,\varphi(x)\,\psi(y)\right]\,\rchi(y).
  \end{aligned}
\end{equation*}
Note that
\begin{equation*}
  \begin{aligned}
    \rchi\,\left( S + \psi\otimes\,\epsilon\,\varphi \right)\,\rchi &\doteq K_{1,N}^{(1,1)}(x,y),   \\
    \rchi\,\left( SD - \psi\otimes\,\varphi  \right)\,\rchi &\doteq K_{1,N}^{(1,2)}(x,y),  \\
    \rchi\,\left( \epsilon\,S + \epsilon\,\psi\otimes\,\epsilon\,\varphi - \epsilon \right)\,\rchi &\doteq K_{1,N}^{(2,1)}(x,y) , \\
    \rchi\,\left( S + \epsilon\,\varphi\otimes\,\epsilon\,\psi \right)\,\rchi &\doteq   K_{1,N}^{(2,2)}(x,y). \\
  \end{aligned}
\end{equation*}
Hence
\begin{equation*}
  K_{1,N}=\rchi\,\left(
    \begin{array}{cc}
      S + \psi\otimes\,\epsilon\,\varphi  & SD - \psi\otimes\,\varphi \\
      \epsilon\,S + \epsilon\,\psi\otimes\,\epsilon\,\varphi  - \epsilon &  S + \epsilon\,\varphi\otimes\,\epsilon\,\psi
    \end{array}
  \right)\,\rchi.
\end{equation*}

Note that this is identical to the corresponding operator for $\beta=1$ obtained by Tracy and Widom in \cite{Trac2}, the only difference being that $\varphi$, $\psi$, and hence also $S$, are redefined to depend on $\lambda$.

\section{Edge scaling}

\subsection{Reduction of the determinant}

The above determinant is that of an operator on $L^{2}(J)~\times~L^{2}(J)$. Our first task will be to rewrite these determinants as those of operators on $L^{2}(J)$. This part follows exactly the proof in \cite{Trac2}. To begin, note that \begin{equation}\com{S}{D}=\varphi\otimes\psi + \psi\otimes\varphi\label{sdcom}\end{equation} so that (using the fact that $D\,\epsilon=\epsilon\,D=I$ )
\begin{eqnarray}
  \com{\epsilon}{S} &=&\epsilon\,S-S\,\epsilon\nonumber\\
  &=& \epsilon\,S\,D\,\epsilon-\epsilon\,D\,S\,\epsilon = \epsilon\,\com{S}{D}\,\epsilon\nonumber\\
  &=& \epsilon\,\varphi\otimes\psi\,\epsilon + \epsilon\,\psi\otimes\varphi\,\epsilon\nonumber\\
  &=&\epsilon\,\varphi\otimes\epsilon^{t}\psi + \epsilon\,\psi\otimes\epsilon^{t}\,\varphi\nonumber\\
  &=& - \epsilon\,\varphi\otimes\epsilon\,\psi -
  \epsilon\,\psi\otimes\epsilon\,\varphi,
  \label{escom}\end{eqnarray}
where the last equality follows from the fact that $\epsilon^{t}=-\epsilon$. We thus have

\begin{eqnarray*}
  D\,\left(\epsilon\,S + \epsilon\,\psi\otimes\epsilon\,\varphi\right) & = & S + \psi\otimes\epsilon\,\varphi,\\
  D\,\left(\epsilon\,S\,D - \epsilon\,\psi\otimes\varphi\right) & =
  & S\,D - \psi\otimes\varphi.
\end{eqnarray*}
The expressions on the right side are the top entries of $K_{1,N}$. Thus the first row of $K_{1,N}$ is, as a vector,
\begin{equation*}
  D\,\left(\epsilon\,S + \epsilon\,\psi\otimes\epsilon\,\varphi, \epsilon\,S\,D - \epsilon\,\psi\otimes\varphi\right).
\end{equation*}
Now \eqref{escom} implies that
\begin{equation*}
\epsilon\,S + \epsilon\,\psi\otimes\epsilon\,\varphi = S\,\epsilon -\epsilon\,\varphi\otimes\epsilon\,\psi.
\end{equation*}
Similarly \eqref{sdcom} gives
\begin{equation*}
  \epsilon\,\com{S}{D} = \epsilon\,\varphi\otimes\psi + \epsilon\,\psi\otimes\varphi,
\end{equation*}
so that
\begin{equation*}
  \epsilon\,S\,D - \epsilon\,\psi\otimes\varphi = \epsilon\,D\,S + \epsilon\,\varphi\otimes\psi = S + \epsilon\,\varphi\otimes\psi.
\end{equation*}
Using these expressions we can rewrite the first row of $K_{1,N}$ as
\begin{equation*}
D\,\left(S\,\epsilon - \epsilon\,\varphi\otimes\epsilon\,\psi, S + \epsilon\,\varphi\otimes\psi\right).
\end{equation*}
Applying $\epsilon$ to this expression shows the second row of $K_{1,N}$ is given by
\begin{equation*}
\left(\epsilon\,S - \epsilon + \epsilon\,\psi\otimes\epsilon\,\varphi, S + \epsilon\,\varphi\otimes\psi\right)
\end{equation*}
Now use
\eqref{escom} to show the second row of $K_{1,N}$ is
\begin{equation*}
\left(S\,\epsilon - \epsilon + \epsilon\,\varphi\otimes\epsilon\,\psi, S + \epsilon\,\varphi\otimes\psi\right).
\end{equation*}
Therefore,
\begin{eqnarray*}
  K_{1,N} &=& \rchi\,
  \left(
    \begin{array}{cc} D\,\left(S\,\epsilon - \epsilon\,\varphi\otimes\epsilon\,\psi\right) & D\,\left(S + \epsilon\,\varphi\otimes\psi\right) \\
      S\,\epsilon - \epsilon + \epsilon\,\varphi\otimes\epsilon\,\psi &  S + \epsilon\,\varphi\otimes\psi
    \end{array}
  \right)
  \,\rchi \\
  & =&
  \left(
    \begin{array}{cc}
      \rchi\,D & 0 \\
      0 & \rchi
    \end{array}
  \right)
  \left(
    \begin{array}{cc}
      \left(S\,\epsilon - \epsilon\,\varphi\otimes\epsilon\,\psi\right)\,\rchi & \left(S + \epsilon\,\varphi\otimes\psi\right)\,\rchi \\
      \left(S\,\epsilon - \epsilon + \epsilon\,\varphi\otimes\epsilon\,\psi\right)\rchi & \left(S + \epsilon\,\varphi\otimes\psi\right)\rchi
    \end{array}
  \right).
\end{eqnarray*}
Since $K_{1,N}$ is of the form $A\,B$, we can use the fact that $\det(I-A\,B)=~\det(I-~B\,A)$ and deduce that $D_{1,N}(s,\lambda)$ is unchanged if instead we take $K_{1,N}$ to be
\begin{eqnarray*}
  K_{1,N} & = &
  \left(
    \begin{array}{cc}
      \left(S\,\epsilon - \epsilon\,\varphi\otimes\epsilon\,\psi\right)\,\rchi & \left(S + \epsilon\,\varphi\otimes\psi\right)\,\rchi \\
      \left(S\,\epsilon - \epsilon + \epsilon\,\varphi\otimes\epsilon\,\psi\right)\,\rchi & \left(S + \epsilon\,\varphi\otimes\psi\right)\,\rchi
    \end{array}
  \right)
  \left(
    \begin{array}{cc}
      \rchi\,D & 0 \\
      0 & \rchi
    \end{array}
  \right)\\
  & = &
  \left(
    \begin{array}{cc}
      \left(S\,\epsilon - \epsilon\,\varphi\otimes\epsilon\,\psi\right)\,\rchi\,D & \left(S + \epsilon\,\varphi\otimes\psi\right)\,\rchi \\
      \left(S\,\epsilon - \epsilon + \epsilon\,\varphi\otimes\epsilon\,\psi\right)\rchi\,D & \left(S + \epsilon\,\varphi\otimes\psi\right)\rchi
    \end{array}\right).
\end{eqnarray*}
Therefore
\begin{equation}
  D_{1,N}(s,\lambda)=\det
  \left(
    \begin{array}{cc}
      I - \left(S\,\epsilon - \epsilon\,\varphi\otimes\epsilon\,\psi\right)\,\lambda\,\rchi\,D & - \left(S + \epsilon\,\varphi\otimes\psi\right)\,\lambda\,\rchi \\
      - \left(S\,\epsilon - \epsilon + \epsilon\,\varphi\otimes\epsilon\,\psi\right)\,\lambda\,\rchi\,D & I - \left(S + \epsilon\,\varphi\otimes\psi\right)\,\lambda\,\rchi
    \end{array}
  \right).
\end{equation}
Now we perform row and column operations on the matrix to simplify it, which do not change the Fredholm determinant. Justification of these operations is given in \cite{Trac2}. We start by subtracting row 1 from row 2 to get

\begin{equation*}  \left(\begin{array}{cc}
I - \left(S\,\epsilon - \epsilon\,\varphi\otimes\epsilon\,\psi\right)\,\lambda\,\rchi\,D & - \left(S + \epsilon\,\varphi\otimes\psi\right)\,\lambda\,\rchi \\
  - I + \epsilon\,\lambda\,\rchi\,D  & I \end{array}\right). \end{equation*}
Next, adding column 2 to column 1 yields
\begin{equation*}
  \left(
    \begin{array}{cc}
      I - \left(S\,\epsilon - \epsilon\,\varphi\otimes\epsilon\,\psi\right)\,\lambda\,\rchi\,D - \left(S + \epsilon\,\varphi\otimes\psi\right)\,\lambda\,\rchi & - \left(S + \epsilon\,\varphi\otimes\psi\right)\,\lambda\,\rchi \\
      \epsilon\,\lambda\,\rchi\,D  & I
    \end{array}
  \right).
\end{equation*}
Then right-multiply column 2 by $-\epsilon\,\lambda\,\rchi\,D$ and add it to column 1, and multiply row 2 by $S + \epsilon\,\varphi\otimes\psi$ and add it to row 1 to arrive at
\begin{equation*}
  \det
  \left(
    \begin{array}{cc}
      I - \left(S\,\epsilon - \epsilon\,\varphi\otimes\epsilon\,\psi\right)\,\lambda\,\rchi\,D +  \left(S + \epsilon\,\varphi\otimes\psi\right)\,\lambda\,\rchi\,\left(\epsilon\,\lambda\,\rchi\,D - I\right) & 0 \\
      0  & I
    \end{array}
  \right).
\end{equation*}
Thus the determinant we want equals the determinant of
\begin{equation}
  I - \left(S\,\epsilon - \epsilon\,\varphi\otimes\epsilon\,\psi\right)\,\lambda\,\rchi\,D + \left(S +   \epsilon\,\varphi\otimes\psi\right)\,\lambda\,\rchi\,\left(\epsilon\,\lambda\,\rchi\,D - I\right)
  \label{operator}.
\end{equation}
So we have reduced the problem from the computation of the Fredholm determinant of an operator on $L^{2}(J)~\times~L^{2}(J)$, to that of an operator on $L^{2}(J)$.

\subsection{Differential equations}

Next we want to write the operator in \eqref{operator} in the form
\begin{equation} \left(I- K_{2, N}\right)\left(I -
\sum_{i=1}^{L}\alpha_{i}\otimes\beta_{i}\right), \end{equation}
where the $\alpha_{i}$  and $\beta_{i}$ are functions in
$L^{2}(J)$. In other words, we want to rewrite the determinant for
the GOE case as a finite dimensional perturbation of the
corresponding GUE determinant. The Fredholm determinant of the
product is then the product of the determinants. The limiting form
for the GUE part is already known, and we can just focus on
finding a limiting form for the determinant of the finite
dimensional piece. It is here that the proof must be modified from
that in \cite{Trac2}. A little simplification of \eqref{operator}
yields
\begin{equation*}
I -
\lambda\,S\,\rchi-\lambda\,S\,\left(1-\lambda\,\rchi\right)\,\epsilon\,\rchi\,D
- \lambda\,\left(\epsilon\,\varphi\,\otimes\,\rchi\,\psi\right) -
\lambda\,\left(\epsilon\,\varphi\,\otimes\,\psi\right)\left(1-\lambda\,\rchi\right)\,\epsilon\,\rchi\,D.
\end{equation*}
Writing $\epsilon\,\com{\rchi}{D}+\rchi$ for $\epsilon\,\rchi\,D$
and simplifying $\left(1-\lambda\,\rchi\right)\,\rchi$ to
$\left(1-\lambda\right)\,\rchi$ gives

 \begin{align*}& I -
\lambda\,S\,\rchi - \lambda\,\left(1-\lambda\right)\,S\,\rchi -
\lambda\,\left(\epsilon\,\varphi\,\otimes\,\rchi\,\psi\right)
-\lambda\,\left(1-\lambda\right)\,\left(\epsilon\,\varphi\,\otimes\,\rchi\,\psi\right)
\\ & \qquad
-\lambda\,S\,\left(1-\lambda\,\rchi\right)\,\epsilon\,\com{\rchi}{D}
-
\lambda\,\left(\epsilon\,\varphi\,\otimes\,\psi\right)\,\left(1-\lambda\,\rchi\right)\,\epsilon\,\com{\rchi}{D}
 \\ & = I -  (2\lambda-\lambda^{2})\,S\,\rchi - (2\lambda-\lambda^{2})\,(\epsilon\,\varphi\,\otimes\,\rchi\,\psi)
  - \lambda\,S\,(1-\lambda\,\rchi)\,\epsilon\,\com{\rchi}{D} \\ &
  \qquad
  -  \lambda\,(\epsilon\,\varphi\,\otimes\,\psi)\,(1-\lambda\,\rchi)\,\epsilon\,\com{\rchi}{D}.\end{align*}
Define $\tilde{\lambda}=2\,\lambda-\lambda^{2}$ and let
$\sqrt{\tilde{\lambda}}\,\varphi\to\varphi$, and
$\sqrt{\tilde{\lambda}}\,\psi\to\psi$ so that
$\tilde{\lambda}\,S\to S$ and \eqref{operator} goes to
\begin{align*} I -  & S\,\rchi - (\epsilon\,\varphi\,\otimes\,\rchi\,\psi)
  - \frac{\lambda}{\tilde{\lambda}}\,S\,(1-\lambda\,\rchi)\,\epsilon\,\com{\rchi}{D} \\ & -
 \frac{\lambda}{\tilde{\lambda}}\,(\epsilon\,\varphi\,\otimes\,\psi)\,(1-\lambda\,\rchi)\,\epsilon\,\com{\rchi}{D}.\end{align*}
Now we define $R:=(I-S\,\rchi)^{-1}\,S\,\rchi=(I-S\,\rchi)^{-1}-I$
(the resolvent operator of $S\,\rchi$), whose kernel we denote by
$R(x,y)$, and
$Q_{\epsilon}:=(I-S\,\rchi)^{-1}\,\epsilon\,\varphi$. Then
\eqref{operator} factors into

\begin{equation*}A =(I - S\,\rchi)\,B.\end{equation*}
where $B$ is

\begin{align*}I - & (Q_{\epsilon}\,\otimes\,\rchi\,\psi)
  - \frac{\lambda}{\tilde{\lambda}}\,(I+R)\,S\,(1-\lambda\,\rchi)\,\epsilon\,\com{\rchi}{D}\\ & -
  \frac{\lambda}{\tilde{\lambda}}\,(Q_{\epsilon}\,\otimes\,\psi)\,(1-\lambda\,\rchi)\,\epsilon\,\com{\rchi}{D},\qquad \lambda\neq 1.\end{align*} Hence
\begin{equation*}D_{1,N}(s,\lambda)
=D_{2,N}(s,\tilde{\lambda})\,\det(B).\end{equation*}
Note that because of the change of variable $\tilde{\lambda}\,S\to S$, we are in effect factoring $I-(2\lambda-\lambda^{2})\,S$, rather that $I-\lambda\,S$ as we did in the $\beta=4$ case. The fact that we factored $I -  (2\lambda-\lambda^{2})\,S\,\rchi$ as opposed to $I - \lambda\,S\,\rchi$ is crucial here for it is what makes $B$ finite rank. If we had factored $I - \lambda\,S\,\rchi$ instead, $B$ would have been
\[\begin{split} B =  I & -
\lambda\,\sum_{k=1}^{2}(-1)^{k}\left(S + R\,S\right)\,\left(I-\lambda\,\rchi\right)\epsilon_{k}\otimes\delta_{k} - \lambda\,\left(I + R\right)\,\epsilon\,\varphi\,\otimes\,\rchi\,\psi \\
& - \lambda\,\sum_{k=1}^{2}(-1)^{k}\inprod{\psi}{\left(I -\lambda\,\rchi\right)\,\epsilon_{k}}\,\left(\left(I + R\right)\,\epsilon\,\varphi\right)\otimes\delta_{k} \\ & - \lambda\,(1-\lambda)\,\left(S + R\,S\right)\,\rchi - \lambda\,(1-\lambda)\,\left(\left(I + R\right)\,\epsilon\,\varphi\right)\otimes\rchi\,\psi
\end{split}\]
The first term on the last line is not finite rank, and the methods we have used previously in the $\beta=4$ case would not work here. It is also interesting to note that these complications disappear when we are dealing with the case of the largest eigenvalue; then is no differentiation with respect to $\lambda$, and we just set $\lambda=1$ in all these formulae. All the new troublesome terms vanish!
\par\noindent In order to find $\det(B)$ we use the identity
\begin{equation}\epsilon\,\com{\rchi}{D}=\sum_{k=1}^{2m}
(-1)^{k}\,\epsilon_{k}\otimes\delta_{k},\end{equation} where
$\epsilon_{k}$ and $\delta_{k}$ are the functions
$\epsilon(x-a_{k})$ and $\delta(x-a_{k})$ respectively, and the
$a_{k}$ are the endpoints of the (disjoint) intervals considered,
$J=\cup_{k=1}^{m}(a_{2\,k-1},a_{2\,k})$. We also make use of the
fact that
\begin{equation}a\otimes b\cdot c\otimes d= \inprod{b}{c}\cdot a\otimes d\end{equation}
where $\inprod{.}{.}$ is the usual $L^{2}$--inner product.
Therefore

\begin{align*}(Q_{\epsilon}\otimes\psi)\,(1-\lambda\,\rchi)\,\epsilon\,\com{\rchi}{D} &=
\sum_{k=1}^{2m} (-1)^{k}Q_{\epsilon}\otimes\psi\cdot
(1-\lambda\,\rchi)\,\epsilon_{k}\otimes\delta_{k}
\\ &=\sum_{k=1}^{2m}
(-1)^{k}\inprod{\psi}{(1-\lambda\,\rchi)\,\epsilon_{k}}\,Q_{\epsilon}\otimes\,\delta_{k}.
\end{align*}
It follows that

\begin{equation*}\frac{D_{1,N}(s,\lambda)}{D_{2,N}(s,\tilde{\lambda})}\end{equation*}
equals the determinant of
\begin{equation}
I  - Q_{\epsilon}\otimes\rchi\psi - \frac{\lambda}{\tilde{\lambda}}\,\sum_{k=1}^{2m} (-1)^{k}\left[(S+R\,S)\,(1-\lambda\,\rchi)\,\epsilon_{k} +\inprod{\psi}{(1-\lambda\,\rchi)\,\epsilon_{k}}\,Q_{\epsilon}\right]\otimes\delta_{k}.
\end{equation}
We now specialize to the case of one interval $J=(t,\infty)$, so
$m=1$, $a_{1}=t$ and $a_{2}=\infty$. We write
$\epsilon_{t}=\epsilon_{1}$, and $\epsilon_{\infty}=\epsilon_{2}$,
and similarly for $\delta_{k}$. Writing the terms in the summation
and using the facts that
\begin{equation}\epsilon_{\infty}=-\frac{1}{2},\end{equation}
 and
 \begin{equation}(1-\lambda\,\rchi)\,\epsilon_{t}=-\frac{1}{2}\,(1-\lambda\,\rchi)+(1-\lambda\,\rchi)\,\rchi,\end{equation}
 then yields
\begin{align*}
I - Q_{\epsilon}\otimes\rchi\psi  -
\frac{\lambda}{2\tilde{\lambda}}\,
\left[(S+R\,S)\,(1-\lambda\,\rchi)+
\inprod{\psi}{(1-\lambda\,\rchi)}\,Q_{\epsilon}\right]\otimes(\delta_{t}-\delta_{\infty})\\
\qquad   +\frac{\lambda}{\tilde{\lambda}}
\left[(S+R\,S)\,(1-\lambda\,\rchi)\,\rchi+
\inprod{\psi}{(1-\lambda\,\rchi)\,\rchi}\,Q_{\epsilon}\right]\otimes\delta_{t}
\end{align*}
which, to simplify notation, we write as
\begin{align*}
 I - Q_{\epsilon}\otimes\rchi\psi  - \frac{\lambda}{2\tilde{\lambda}}\,
\left[(S+R\,S)\,(1-\lambda\,\rchi)+
a_{1,\lambda}\,Q_{\epsilon}\right]\otimes(\delta_{t}-\delta_{\infty}) \\
\qquad   +\frac{\lambda}{\tilde{\lambda}}
\left[(S+R\,S)\,(1-\lambda\,\rchi)\,\rchi+
\tilde{a}_{1,\lambda}\,Q_{\epsilon}\right]\otimes\delta_{t},
\end{align*}
where
\begin{equation}a_{1,\lambda} = \inprod{\psi}{(1-\lambda\,\rchi)},\qquad \tilde{a}_{1,\lambda}= \inprod{\psi}{(1-\lambda\,\rchi)\,\rchi}.\end{equation}
 Now we can use the formula:
\begin{equation}\det\left(I-\sum_{i=1}^{L}\alpha_{i}\otimes\beta_{i}\right)=\det\left(\delta_{jk}-\inprod{\alpha_{j}}{\beta_{k}}\right)_{1\leq j,k\leq L}\end{equation}
In this case, $L=3$, and
\begin{align}\alpha_{1}=Q_{\epsilon}, &\qquad
 \alpha_{2} =\frac{\lambda}{\tilde{\lambda}}\,
\left[(S+R\,S)\,(1-\lambda\,\rchi)+ a_{1,\lambda}\,Q_{\epsilon}\right]\nonumber, \\
\alpha_{3}= & -\frac{\lambda}{\tilde{\lambda}}
\left[(S+R\,S)\,(1-\lambda\,\rchi)\,\rchi+
\tilde{a}_{1,\lambda}\,Q_{\epsilon}\right],\nonumber \\ &
\beta_{1}=\rchi\psi, \qquad  \beta_{2}=\delta_{t}-\delta_{\infty},
\qquad \beta_{3}=\delta_{t}.\end{align} In order to simplify the
notation, define

\begin{align}
Q(x,\lambda,t)&:=(I-S\,\rchi)^{-1}\,\varphi, &
P(x,\lambda,t)&:=(I-S\,\rchi)^{-1}\,\psi, \nonumber\\
Q_{\epsilon}(x,\lambda,t)&:=(I-S\,\rchi)^{-1}\,\epsilon\,\varphi,
& P_{\epsilon}(x,\lambda,t)&:=(I-S\,\rchi)^{-1}\,\epsilon\,\psi,
\end{align}

\begin{align}
q_{_{N}}&:=Q(t,\lambda,t), & p_{_{N}}&:=P(t,\lambda,t),\nonumber\\
q_{\epsilon}&:=Q_{\epsilon}(t,\lambda,t), & p_{\epsilon}&:=P_{\epsilon}(t,\lambda,t),\nonumber\\
u_{\epsilon}&:=\inprod{Q}{\rchi\,\epsilon\,\varphi}=\inprod{Q_{\epsilon}}{\rchi\,\varphi},
&
v_{\epsilon}&:=\inprod{Q}{\rchi\,\epsilon\,\psi}=\inprod{P_{\epsilon}}{\rchi\,\psi}, \nonumber\\
\tilde{v}_{\epsilon}&:=\inprod{P}{\rchi\,\epsilon\,\varphi}=\inprod{Q_{\epsilon}}{\rchi\,\varphi},
&
w_{\epsilon}&:=\inprod{P}{\rchi\,\epsilon\,\psi}=\inprod{P_{\epsilon}}{\rchi\,\psi},
\label{inproddef1}\end{align}

\begin{align}
\mathcal{P}_{1,\lambda}&:= \int(1-\lambda\,\rchi)\,P\,d\,x, &
\tilde{\mathcal{P}}_{1,\lambda} &:=
\int(1-\lambda\,\rchi)\,\rchi\,P\,d\,x, \nonumber\\
\mathcal{Q}_{1,\lambda} &:= \int(1-\lambda\,\rchi)\,Q\,d\,x, &
\tilde{\mathcal{Q}}_{1,\lambda} &:=
\int(1-\lambda\,\rchi)\,\rchi\,Q\,d\,x,\nonumber\\
\mathcal{R}_{1,\lambda}&:= \int(1-\lambda\,\rchi)\,R(x,t)\,d\,x, &
\tilde{\mathcal{R}}_{1,\lambda} &:=
\int(1-\lambda\,\rchi)\,\rchi\,R(x,t)\,d\,x.
\label{inproddef2}\end{align} Note that all quantities in
\eqref{inproddef1} and \eqref{inproddef2} are functions of $t$
alone. Furthermore, let
\begin{equation} c_{\varphi} =
\epsilon\,\varphi(\infty)=\frac{1}{2}\int_{-\infty}^{\infty}\varphi(x)\,d\,x,
\qquad c_{\psi} =
\epsilon\,\psi(\infty)=\frac{1}{2}\int_{-\infty}^{\infty}\psi(x)\,d\,x.\end{equation}
Recall from the previous section that when $\beta=1$ we take $N$ to be even. It follows that $\varphi$ and $\psi$ are even and odd functions respectively. Thus $c_{\psi}=0$ for $\beta=1$, and computation gives
\begin{equation}
c_{\varphi}=(\pi N)^{1/4} 2^{-3/4-N/2}\,\frac{(N!)^{1/2}}{(N/2)!}\,\sqrt{\lambda}.
\end{equation}
Hence computation yields
\begin{equation}
 \lim_{N\to\infty}c_{\varphi}= \sqrt{\frac{\lambda}{2}},
\end{equation}
and at $t=\infty$ we have
\begin{equation*}
  u_{\epsilon}(\infty)=0, \quad  q_{\epsilon}(\infty)=c_{\varphi}
\end{equation*}
\begin{equation*}
\mathcal{P}_{1,\lambda}(\infty) = 2\,c_{\psi}, \quad
\mathcal{Q}_{1,\lambda}(\infty) =2\,c_{\varphi},\quad
\mathcal{R}_{1,\lambda}(\infty) = 0 ,
\end{equation*}
\begin{equation*}
\tilde{\mathcal{P}}_{1,\lambda}(\infty) =
\tilde{\mathcal{Q}}_{1,\lambda}(\infty) =\tilde{
\mathcal{R}}_{1,\lambda}(\infty) = 0.
\end{equation*}
Hence
\begin{align}
\inprod{\alpha_{1}}{\beta_{1}} & = \tilde{v}_{\epsilon}, \quad
\inprod{\alpha_{1}}{\beta_{2}}=q_{\epsilon}-c_{\varphi}, \quad
\inprod{\alpha_{1}}{\beta_{3}}=q_{\epsilon},
\\
\inprod{\alpha_{2}}{\beta_{1}} &
=\frac{\lambda}{2\,\tilde{\lambda}}
\,\left[\mathcal{P}_{1,\lambda}-a_{1,\lambda}\,(1-\tilde{v}_{\epsilon})\right], \\
\inprod{\alpha_{2}}{\beta_{2}} & =
\frac{\lambda}{2\,\tilde{\lambda}}\,\left[\mathcal{R}_{1,\lambda}
+ a_{1,\lambda}\,(q_{\epsilon}-c_{\varphi})\right]\label{example}, \\
\inprod{\alpha_{2}}{\beta_{3}} & =
\frac{\lambda}{2\,\tilde{\lambda}}\,\left[\mathcal{R}_{1,\lambda}
+ a_{1,\lambda}\,q_{\epsilon}\right], \\
\inprod{\alpha_{3}}{\beta_{1}} & =
-\frac{\lambda}{\tilde{\lambda}}\,\left[\tilde{\mathcal{P}}_{1,\lambda}
- \tilde{a}_{1,\lambda}\,(1-\tilde{v}_{\epsilon})\right], \\
\inprod{\alpha_{3}}{\beta_{2}} & =
-\frac{\lambda}{\tilde{\lambda}}\,\left[\tilde{\mathcal{R}}_{1,\lambda}+\tilde{a}_{1,\lambda}\,(q_{\epsilon}
- c_{\varphi})\right], \\
 \inprod{\alpha_{3}}{\beta_{3}} & =
-\frac{\lambda}{\tilde{\lambda}}\,
\left[\tilde{\mathcal{R}}_{1,\lambda}+\tilde{a}_{1,\lambda}\,q_{\epsilon}\right].
\end{align}
As an illustration, let us do the computation that led to
\eqref{example} in detail. As in \cite{Trac2}, we use the facts
that $S^{t}=S$, and $(S+S\,R^{t})\,\rchi=R$ which can be easily
seen by writing $R=\sum_{k=1}^{\infty}(S\,\rchi)^{k}$. Furthermore
we write $R(x,a_{k})$ to mean
\begin{equation*} \lim_{\substack{y\to a_{k}\\y\in J}}R(x,y).\end{equation*}
In general, since all evaluations are done by taking the limits
from within $J$, we can use the identity
$\rchi\,\delta_{k}=\delta_{k}$ inside the inner products. Thus
\begin{align*}
\inprod{\alpha_{2}}{\beta_{2}} & =
\frac{\lambda}{\tilde{\lambda}}\,\left[
\inprod{(S+R\,S)\,(1-\lambda\,\rchi)}{\delta_{t}-\delta_{\infty}}
+
a_{1,\lambda}\,\inprod{Q_{\epsilon}}{\delta_{t}-\delta_{\infty}}\right]\\
&=\frac{\lambda}{\tilde{\lambda}}\,\left[\inprod{(1-\lambda\,\rchi)}{(S+R^{t}\,S)\,\left(\delta_{t}-\delta_{\infty}\right)}
+ a_{1,\lambda}\left(Q_{\epsilon}(t)-Q_{\epsilon}(\infty)\right)\right]\\
&
=\frac{\lambda}{\tilde{\lambda}}\,\left[\inprod{(1-\lambda\,\rchi)}{(S+R^{t}\,S)\,\rchi\,\left(\delta_{t}-\delta_{\infty}\right)}
+ a_{1,\lambda}\left(q_{\epsilon}-c_{\varphi}\right)\right]\\
&
=\frac{\lambda}{\tilde{\lambda}}\,\left[\inprod{(1-\lambda\,\rchi)}{R(x,t)-R(x,\infty)}
+ a_{1,\lambda}\left(q_{\epsilon}-c_{\varphi}\right)\right]\\
&
=\frac{\lambda}{\tilde{\lambda}}\,\left[\mathcal{R}_{1,\lambda}(t)-\mathcal{R}_{1,\lambda}(\infty)
+ a_{1,\lambda}\left(q_{\epsilon}-c_{\varphi}\right)\right]\\
&
=\frac{\lambda}{\tilde{\lambda}}\,\left[\mathcal{R}_{1,\lambda}(t)
+ a_{1,\lambda}\left(q_{\epsilon}-c_{\varphi}\right)\right].
\end{align*}
We want the limit of the determinant
\begin{equation}\det\left(\delta_{jk}-\inprod{\alpha_{j}}{\beta_{k}}\right)_{1\leq
j,k\leq 3},\end{equation} as $N\to \infty$. In order to get our
hands on the limits of the individual terms involved in the
determinant, we will find differential equations for them first as
in \cite{Trac2}. Row operation on the matrix show that $a_{1,\lambda}$ and
$\tilde{a}_{1,\lambda}$ fall out of the determinant; to see this
add $\lambda\,a_{1,\lambda}/(2\,\tilde{\lambda})$ times row 1 to
row 2 and $\lambda\,\tilde{a}_{1,\lambda}/\tilde{\lambda}$ times
row 1 to row 3. So we will not need to find differential equations
for them. Our determinant is
{\large
  \begin{equation}
    \det\left(
      \begin{array}{ccc}
        1 - \tilde{v}_{\epsilon} & -(q_{\epsilon}-c_{\varphi}) & -q_{\epsilon} \\[6pt]
        -\frac{\lambda\,\mathcal{P}_{1,\lambda}}{2\,\tilde{\lambda}} & 1 - \frac{\lambda\,\mathcal{R}_{1,\lambda}}{2\,\tilde{\lambda}}  & -\frac{\lambda\,\mathcal{R}_{1,\lambda}}{2\,\tilde{\lambda}} \\[6pt]
        \frac{\lambda\,\tilde{\mathcal{P}}_{1,\lambda}}{\tilde{\lambda}} & \frac{\lambda\,\tilde{\mathcal{R}}_{1,\lambda}}{\tilde{\lambda}} & 1 + \frac{\lambda\,\tilde{\mathcal{R}}_{1,\lambda}}{\tilde{\lambda}}
      \end{array}
    \right).
  \end{equation}
}
Proceeding as in \cite{Trac2} we find the following differential
equations
\begin{align} \frac{d}{d\,t}\,u_{\epsilon} & = q_{_{N}}\,q_{\epsilon}, & \frac{d}{d\,t}\,q_{\epsilon} & = q_{_{N}} -q_{_{N}}\, \tilde{v}_{\epsilon} - p_{_{N}}\,u_{\epsilon},\\
 \frac{d}{d\,t}\mathcal{Q}_{1,\lambda} & =
q_{_{N}}\left(\lambda-\mathcal{R}_{1,\lambda}\right), &
\frac{d}{d\,t}\mathcal{P}_{1,\lambda} & =
p_{_{N}}\left(\lambda-\mathcal{R}_{1,\lambda}\right),\label{ex2}\\
 \frac{d}{d\,t}\mathcal{R}_{1,\lambda} & =
-p_{_{N}}\,\mathcal{Q}_{1,\lambda}-q_{_{N}}\,\mathcal{P}_{1,\lambda},
  & \frac{d}{d\,t}\tilde{\mathcal{R}}_{1,\lambda} & =
-p_{_{N}}\,\tilde{\mathcal{Q}}_{1,\lambda}-q_{_{N}}\,\tilde{\mathcal{P}}_{1,\lambda},
\\  \frac{d}{d\,t}\tilde{\mathcal{Q}}_{1,\lambda} & =
q_{_{N}}\left(\lambda-1-\tilde{\mathcal{R}}_{1,\lambda}\right),
 & \frac{d}{d\,t}\tilde{\mathcal{P}}_{1,\lambda}  & =
p_{_{N}}\left(\lambda-1-\tilde{\mathcal{R}}_{1,\lambda}\right).\end{align}
Let us derive the first equation in \eqref{ex2} for example. From
\cite{Trac3} (equation $2.17$), we have
\begin{equation*}
\frac{\partial Q}{\partial t}=-R(x,t)\,q_{_{N}}.
\end{equation*}
Therefore
\begin{align*}
\frac{\partial \mathcal{Q}_{1,\lambda}}{\partial t} & =
\frac{d}{d\,t}\left[\int_{-\infty}^{t}Q(x,t)\,d\,x-(1-\lambda)\,\int_{\infty}^{t}Q(x,t)\,d\,x\right]\\
& = q_{_{N}} + \int_{-\infty}^{t}\frac{\partial Q}{\partial
t}\,d\,x -
(1-\lambda)\left[q_{_{N}}+\int_{\infty}^{t}\frac{\partial
Q}{\partial t}\,d\,x\right]
\\
& = q_{_{N}} - q_{_{N}}\int_{-\infty}^{t}R(x,t)\,d\,x -
(1-\lambda)\,q_{_{N}}+(1-\lambda)\,q_{_{N}}\,\int_{\infty}^{t}R(x,t)\,d\,x \\
& = \lambda\,q_{_{N}}-q_{_{N}}\,\int_{-\infty}^{\infty}(1-\lambda)\,R(x,t)\,d\,x\\
& = \lambda\,q_{_{N}}-q_{_{N}}\,\mathcal{R}_{1,\lambda} =
q_{_{N}}\left(\lambda-\mathcal{R}_{1,\lambda}\right).
\end{align*}
 Now we change variable from $t$ to $s$ where $t=\tau(s)=
2\,\sigma\,\sqrt{N}+\frac{\sigma\,s}{N^{1/6}}$. Then we take the
limit $N\to \infty$, denoting the limits of $ q_{\epsilon},
\mathcal{P}_{1,\lambda},
\mathcal{Q}_{1,\lambda},\mathcal{R}_{1,\lambda},
\tilde{\mathcal{P}}_{1,\lambda} ,\tilde{\mathcal{Q}}_{1,\lambda} ,
\tilde{\mathcal{R}}_{1,\lambda}$ and the common limit of
$u_{\epsilon}$ and $\tilde{v}_{\epsilon}$ respectively by
$\overline{q}, \overline{\mathcal{P}}_{1,\lambda},
\overline{\mathcal{Q}}_{1,\lambda},
\overline{\mathcal{R}}_{1,\lambda},
\overline{\overline{\mathcal{P}}}_{1,\lambda} ,
\overline{\overline{\mathcal{Q}}}_{1,\lambda} ,
\overline{\overline{\mathcal{R}}}_{1,\lambda}$ and $\overline{u}$.
We eliminate $ \overline{\mathcal{Q}}_{1,\lambda}$ and
$\overline{\overline{\mathcal{Q}}}_{1,\lambda}$ by using the facts
that
$\overline{\mathcal{Q}}_{1,\lambda}=\overline{\mathcal{P}}_{1,\lambda}+\lambda\,\sqrt{2}$
and $\overline{\overline{\mathcal{Q}}}_{1,\lambda}=
\overline{\overline{\mathcal{P}}}_{1,\lambda}$. These limits hold
uniformly for bounded $s$ so we can interchange $\lim$ and
$\frac{d}{d\,s}$. Also
$\lim_{N\to\infty}N^{-1/6}q_{_{N}}=\lim_{N\to\infty}N^{-1/6}p_{_{N}}=q
$ , where $q$ is as in \eqref{D2}. We obtain the systems
\begin{equation}\frac{d}{d\,s}\,\overline{u} = -\frac{1}{\sqrt{2}}\,q\,\overline{q}
,
 \qquad
\frac{d}{d\,s}\,\overline{q} =
\frac{1}{\sqrt{2}}\,q\,\left(1-2\,\overline{u}\right),
\end{equation}
\begin{equation}\frac{d}{d\,s}\overline{\mathcal{P}}_{1,\lambda} =
-\frac{1}{\sqrt{2}}\,q\,\left(\overline{\mathcal{R}}_{1,\lambda}-\lambda\right),
\qquad \frac{d}{d\,s}\overline{\mathcal{R}}_{1,\lambda} =
-\frac{1}{\sqrt{2}}\,q\,\left(2\,\overline{\mathcal{P}}_{1,\lambda}+\sqrt{2\tilde{\lambda}}\right),
\end{equation}
\begin{equation}\frac{d}{d\,s}\overline{\overline{\mathcal{P}}}_{1,\lambda} =
\frac{1}{\sqrt{2}}\,q\,\left(1-\lambda-\overline{\overline{\mathcal{R}}}_{1,\lambda}\right),
\qquad
 \frac{d}{d\,s}\overline{\overline{\mathcal{R}}}_{1,\lambda} =
 -q\,\sqrt{2}\,\overline{\overline{\mathcal{P}}}_{1,\lambda}.
\end{equation}
 The change of variables $s\to\mu=\int_{s}^{\infty} q(x, \lambda)\,d\,x$
transforms these systems into constant coefficient ordinary differential equations

\begin{equation}\frac{d}{d\,\mu}\overline{u} =
\frac{1}{\sqrt{2}}\,\overline{q}, \qquad
\frac{d}{d\,\mu}\overline{q} =
-\frac{1}{\sqrt{2}}\,\left(1-2\,\overline{u}\right),
\end{equation}
\begin{equation}\frac{d}{d\,\mu}\overline{\mathcal{P}}_{1,\lambda} =
\frac{1}{\sqrt{2}}\,\left(\overline{\mathcal{R}}_{1,\lambda}-\lambda\right),
\qquad \frac{d}{d\,\mu}\overline{\mathcal{R}}_{1,\lambda} =
\frac{1}{\sqrt{2}}\,\left(2\,\overline{\mathcal{P}}_{1,\lambda}+\sqrt{2\tilde{\lambda}}\right),
\end{equation}
\begin{equation}\frac{d}{d\,\mu}\overline{\overline{\mathcal{P}}}_{1,\lambda} =
-\frac{1}{\sqrt{2}}\,\left(1-\lambda-\overline{\overline{\mathcal{R}}}_{1,\lambda}\right),
\qquad
 \frac{d}{d\,\mu}\overline{\overline{\mathcal{R}}}_{1,\lambda} =
 \sqrt{2}\,\overline{\overline{\mathcal{P}}}_{1,\lambda}.
\end{equation} Since $\lim_{s\to \infty}\mu=0$, corresponding to the boundary
values at $t=\infty$ which we found earlier for
$\mathcal{P}_{1,\lambda}, \mathcal{R}_{1,\lambda},
\tilde{\mathcal{P}}_{1,\lambda} ,
\tilde{\mathcal{R}}_{1,\lambda}$, we now have initial values at
$\mu=0$. Therefore
\begin{equation} \overline{\mathcal{P}}_{1,\lambda}(0) = \overline{\mathcal{R}}_{1,\lambda}(0)=
\overline{\overline{\mathcal{P}}}_{1,\lambda}(0)=
\overline{\overline{\mathcal{R}}}_{1,\lambda}(0)=0.\end{equation}
We use this to solve the systems and get
\begin{align}\overline{q} & = \frac{\sqrt{\tilde{\lambda}}-1}{2\,\sqrt{2}}\,e^{\mu} + \frac{\sqrt{\tilde{\lambda}} + 1}{2\,\sqrt{2}}\,e^{-\mu},\\
 \overline{u} & =
\frac{\sqrt{\tilde{\lambda}}-1}{4}\,e^{\mu} -
\frac{\sqrt{\tilde{\lambda}}+1}{4}\,e^{-\mu}+\frac{1}{2},\\
\overline{\mathcal{P}}_{1,\lambda} & =
\frac{\sqrt{\tilde{\lambda}}-\lambda}{2\,\sqrt{2}}\,e^{\mu} +
\frac{\sqrt{\tilde{\lambda}} + \lambda}{2\,\sqrt{2}}\,e^{-\mu} -
\sqrt{\frac{\tilde{\lambda}}{2}},\\
 \overline{\mathcal{R}}_{1,\lambda} & =
\frac{\sqrt{\tilde{\lambda}}-\lambda}{2}\,e^{\mu} -
\frac{\sqrt{\tilde{\lambda}} + \lambda}{2}\,e^{-\mu} +
\lambda, \\
\overline{\overline{\mathcal{P}}}_{1,\lambda} & =
\frac{1-\lambda}{2\,\sqrt{2}}\,(e^{\mu}-e^{-\mu}), \qquad
\overline{\overline{\mathcal{R}}}_{1,\lambda}=\frac{1-\lambda}{2}\,(e^{\mu}+e^{-\mu}-2).
\end{align}
Substituting these expressions into the determinant gives
\eqref{goedet}, namely
\begin{equation} D_{1}(s,\lambda)= D_{2}(s,\tilde{\lambda})\,\frac{\lambda - 1
- \cosh{\mu(s,\tilde{\lambda})} +
\sqrt{\tilde{\lambda}}\,\sinh{\mu(s,\tilde{\lambda})}}{\lambda -
2},\end{equation} where $D_{\beta}=\lim_{N\to\infty}D_{\beta,N}$.
As mentioned in the Introduction, the functional form of the $\beta=1$ limiting determinant is very different from what one would expect, unlike in the $\beta=4$ case. Also noteworthy is the dependence on $\tilde{\lambda}=2\lambda-\lambda^{2}$ instead of just $\lambda$. However one should also note that when $\lambda$ is set equal to $1$, then $\tilde{\lambda}=\lambda=1$. Hence in the largest eigenvalue case, where there is no prior differentiation with respect to $\lambda$, and $\lambda$ is just set to $1$, a great deal of simplification occurs. The above formula then nicely reduces to the $\beta=1$ Tracy-Widom distribution.

\chapter{Applications}

\section{An Interlacing property}

The following series of lemmas establish
Corollary~\eqref{interlacingcor}:
\begin{lemma}
Define
\begin{equation}a_{j}=\frac{d^{j}}{d\,\lambda^{j}}\,\sqrt{\frac{\lambda}{2-\lambda}}\,\,\bigg{\vert}_{\lambda=1}.\label{ajdef}\end{equation}
Then $a_{j}$ satisfies the following recursion
\begin{equation} a_{j} = \begin{cases}
          \quad 1 &\text{if} \quad j=0, \\
           \quad (j-1)\,a_{j-1} &\text{for $j\geq 1$, $j$ even,}  \\
           \quad j\,a_{j-1} &\text{for $j\geq 1$, $j$ odd.} \\
         \end{cases}\end{equation}
\label{ajlemma}\end{lemma}
\begin{proof}
Consider the expansion of the generating function $f(\lambda)=\sqrt{\frac{\lambda}{2-\lambda}}$ around $\lambda=1$
\begin{equation*}
f(\lambda)= \sum_{j\geq 0}\frac{a_{j}}{j!}\,(\lambda-1)^{j}=\sum_{j\geq 0} b_{j}\,(\lambda-1)^{j}.
\end{equation*}
Since $a_{j}=j!\,b_{j}$, the statement of the lemma reduces to proving the following recurrence for the $b_{j}$
\begin{equation}
  b_{j} =
  \begin{cases}
    \quad 1 &\text{if} \quad j=0, \\
    \quad \frac{j-1}{j}\,b_{j-1} &\text{for $j\geq 1$, $j$ even,}  \\
    \quad b_{j-1} &\text{for $j\geq 1$, $j$ odd.} \\
  \end{cases}\label{bjrecurrence}
\end{equation}
Let
\begin{equation*}
  f^{even}(\lambda)= \frac{1}{2}\left(\sqrt{\frac{\lambda}{2-\lambda}}+\sqrt{\frac{2-\lambda}{\lambda}}\right), \qquad f^{odd}(\lambda)= \frac{1}{2}\left(\sqrt{\frac{\lambda}{2-\lambda}}-\sqrt{\frac{2-\lambda}{\lambda}}\right).
\end{equation*}
These are the even and odd parts of $f$ relative to the reflection $\lambda-1\to-(\lambda-1)$ or $\lambda\to 2-\lambda$. Recurrence \eqref{bjrecurrence} is equivalent to
\begin{equation*}
  \frac{d}{d\,\lambda}\,f^{even}(\lambda)=(\lambda-1)\,\frac{d}{d\,\lambda}\,f^{odd}(\lambda),
\end{equation*}
which is easily shown to be true.
\end{proof}
\begin{lemma}
  \label{flemma}
  Define
  \begin{equation}
    f(s,\lambda)=1-\sqrt{\frac{\lambda}{2-\lambda}}\,\,\tanh{\frac{\mu(s,\tilde{\lambda})}{2}},
  \end{equation}
  for $\tilde{\lambda}=2\,\lambda-\lambda^{2}$. Then
  \begin{equation}\frac{\partial^{2\,n}}{\partial\,\lambda^{2\,n}}\,f(s,\lambda)\,\,\bigg{\vert}_{\lambda=1}-
    \frac{1}{2\,n+1}\,\frac{\partial^{2\,n+1}}{\partial\,\lambda^{2\,n+1}}\,f(s,\lambda)\,\,\bigg{\vert}_{\lambda=1}=
    \begin{cases}
      \quad 1 &\text{if $n=0$,}\\
      \quad0 &\text{if $n\geq 1$.}  \\
    \end{cases} \label{flemmaeq}\end{equation}
\end{lemma}
\begin{proof}
The case $n=0$ is readily checked. The main ingredient for the
general case is Fa\'a di Bruno's formula
\begin{equation}
\frac{d^{n}}{d t^{n}}g(h(t))=\sum\frac{n!}{k_{1}!\cdots k_{n}!}
\left(\frac{d^{k}g}{dh^{k}}(h(t))\right)\left(\frac{1}{1!}\frac{dh}{d
t}\right)^{k_{1}}\cdots\left(\frac{1}{n!}\frac{d^{n} h}{d
t^{n}}\right)^{k_{n}}, \label{faa}\end{equation} where
$k=\sum_{i=1}^{n}k_{i}$ and the above sum is over all partitions
of $n$, that is all values of $k_{1},\ldots, k_{n}$ such that
$\sum_{i=1}^{n} i\,k_{i}=n$. We apply Fa\'a di Bruno's formula to
derivatives of the function
$\tanh{\frac{\mu(s,\tilde{\lambda})}{2}}$, which we treat as some
function $g(\tilde{\lambda}(\lambda))$. Notice that for $j\geq 1$,
$\frac{d^{j}\tilde{\lambda}}{d\,\lambda^{j}}\,\,\big{\vert}_{\lambda=1}$
is nonzero only when $j=2$, in which case it equals $-2$. Hence,
in \eqref{faa}, the only term that survives is the one
corresponding to the partition all of whose parts equal $2$. Thus
we have
\begin{align*}\frac{\partial^{2n-k}}{\partial\,\lambda^{2n-k}}\,&\tanh{\frac{\mu(s,\tilde{\lambda})}{2}}\,\,\bigg{\vert}_{\lambda=1}&&\\
 &= \begin{cases}
           0 &\text{if $k=2j+1$, $j\geq 0$}\\
          \frac{(-1)^{n-j}\,(2\,n-k)!}{(n-j)!} \frac{\partial^{n-j}}{\partial\,\tilde{\lambda}^{n-j}}\,\tanh{\frac{\mu(s,\tilde{\lambda})}{2}}\,\,\bigg{\vert}_{\tilde{\lambda}=1} &\text{for $k=2j$, $j\geq 0$}  \\
         \end{cases} \end{align*}

 \begin{align*}\frac{\partial^{2n-k+1}}{\partial\,\lambda^{2n+1-k}}\,&\tanh{\frac{\mu(s,\tilde{\lambda})}{2}}\,\,\bigg{\vert}_{\lambda=1}
 \\
 &= \begin{cases}
           0 &\text{if $k=2j$, $j\geq 0$}\\
         \frac{(-1)^{n-j}\,(2\,n+1-k)!}{(n-j)!} \frac{\partial^{n-j}}{\partial\,\tilde{\lambda}^{n-j}}\,\tanh{\frac{\mu(s,\tilde{\lambda})}{2}}\,\,\bigg{\vert}_{\tilde{\lambda}=1} &\text{for $k=2j+1$, $j\geq 0$}  \\
         \end{cases} \end{align*}
Therefore, recalling the definition of $a_{j}$ in \eqref{ajdef}
and setting $k=2\,j$, we obtain
\begin{eqnarray*}
\frac{\partial^{2\,n}}{\partial\,\lambda^{2\,n}}\,f(s,\lambda)\,\,\bigg{\vert}_{\lambda=1}&=&
\sum_{k=0}^{2\,n}\binom{2\,n}{k}\,\frac{\partial^{k}}{\partial\,\lambda^{k}}\,\sqrt{\frac{\lambda}{2-\lambda}}\,\frac{\partial^{2\,n-k}}{\partial\,\lambda^{2\,n-k}}\,\tanh{\frac{\mu(s,\tilde{\lambda})}{2}}\,\,\bigg{\vert}_{\lambda=1}\\
&=&
\sum_{j=0}^{n}\frac{(2\,n)!\,(-1)^{n-j}}{(2\,j)!\,(n-j)!}\,a_{2\,j}\,\frac{\partial^{n-j}}{\partial\,\tilde{\lambda}^{n-j}}\,\tanh{\frac{\mu(s,\tilde{\lambda})}{2}}\,\,\bigg{\vert}_{\tilde{\lambda}=1}.
\end{eqnarray*}
Similarly, using $k=2\,j+1$ instead yields
\begin{eqnarray*}
  \frac{\partial^{2\,n+1}}{\partial\,\lambda^{2\,n+1}} \,f(s,\lambda)\,\bigg{\vert}_{\lambda=1} &=& \sum_{k=0}^{2\,n+1}\binom{2\,n+1}{k}\,\frac{\partial^{k}}{\partial\,\lambda^{k}}\,\sqrt{\frac{\lambda}{2-\lambda}}\,\frac{\partial^{2\,n+1-k}}{\partial\,\lambda^{2\,n+1-k}}\,\tanh{\frac{\mu(s,\tilde{\lambda})}{2}}\,\,\bigg{\vert}_{\lambda=1}\\
  &=& (2\,n+1)\,\sum_{j=0}^{n}\frac{(2\,n)!\,(-1)^{n-j}}{(2\,j)!\,(n-j)!}\,\frac{a_{2\,j+1}}{2\,j+1}\,\frac{\partial^{n-j}}{\partial\,\tilde{\lambda}^{n-j}}\,\tanh{\frac{\mu(s,\tilde{\lambda})}{2}}\,\,\bigg{\vert}_{\tilde{\lambda}=1}\\
  &=& (2\,n+1)\,\frac{\partial^{2\,n}}{\partial\,\lambda^{2\,n}}\,f(s,\lambda)\,\,\bigg{\vert}_{\lambda=1},
\end{eqnarray*}
since $a_{_{2\,j+1}}/(2\,j+1)=a_{_{2j}}$. Rearranging this last
equality leads to \eqref{flemmaeq}.
\end{proof}

\begin{lemma}
Let $D_{1}(s,\lambda)$ and $D_{4}(s,\tilde{\lambda})$ be as in
\eqref{goedet} and \eqref{gsedet}. Then
\begin{equation}D_{1}(s,\lambda)=D_{4}(s,\tilde{\lambda})\,
\left(1-\sqrt{\frac{\lambda}{2-\lambda}}\,\tanh{\frac{\mu(s,\tilde{\lambda})}{2}}\right)^{2}.\end{equation}
\end{lemma}
\begin{proof}
Using the facts that $-1-\cosh{x}=-2\,\cosh^{2}\frac{x}{2}$,
$1=\cosh^{2}{x}-\sinh^{2}{x}$ and
$\sinh{x}=2\,\sinh\frac{x}{2}\,\cosh\frac{x}{2}$ we get

\begin{eqnarray*}
  D_{1}(s,\lambda) &=& \frac{-1-\cosh{\mu(s,\tilde{\lambda})}}{\lambda-2}\,D_{2}(s,\tilde{\lambda}) + D_{2}(s,\tilde{\lambda})\,\frac{\lambda + \sqrt{\tilde{\lambda}}\,\sinh{\mu(s,\tilde{\lambda})}}{\lambda-2}\\
  &=&\frac{-2}{\lambda-2}\,D_{4}(s,\tilde{\lambda}) \\
  & & \qquad + D_{2}(s,\tilde{\lambda})\,\frac{\lambda\left(\cosh^{2}{\frac{\mu(s,\tilde{\lambda})}{2}} - \sinh^{2}{\frac{\mu(s,\tilde{\lambda})}{2}}\right) + \sqrt{\tilde{\lambda}}\,\sinh{\mu(s,\tilde{\lambda})}}{\lambda-2}\\
  &=& D_{4}(s,\tilde{\lambda}) + \frac{D_{4}(s,\tilde{\lambda})}{\cosh^{2}\left(\frac{\mu(s,\lambda)}{2}\right)}\,\frac{\lambda\,\sinh^{2}{\frac{\mu(s,\tilde{\lambda})}{2}} - \sqrt{\tilde{\lambda}}\,\sinh{\mu(s,\tilde{\lambda})}}{2 - \lambda}\\
  &=& D_{4}(s,\tilde{\lambda})\,\left(1 + \frac{\lambda\,\sinh^{2}{\frac{\mu(s,\tilde{\lambda})}{2}} - 2\,\sqrt{\tilde{\lambda}}\,\sinh\left(\frac{\mu(s,\lambda)}{2}\right)\,\cosh\left(\frac{\mu(s,\lambda)}{2}\right)}{(2-\lambda)\,\cosh^{2}\left(\frac{\mu(s,\lambda)}{2}\right)}\right)\\
  &=& D_{4}(s,\tilde{\lambda})\,\left(1 - 2\,\sqrt{\frac{\lambda}{2-\lambda}}\,\tanh{\frac{\mu(s,\tilde{\lambda})}{2}} + \frac{\lambda}{2-\lambda}\,\tanh^{2}{\frac{\mu(s,\tilde{\lambda})}{2}} \right)\\
  &=& D_{4}(s,\tilde{\lambda})\,\left(1-\sqrt{\frac{\lambda}{2-\lambda}}\,\,\tanh{\frac{\mu(s,\tilde{\lambda})}{2}}\right)^{2}.
\end{eqnarray*}
\end{proof}
For notational convenience, define $d_{1}(s,\lambda)=D_{1}^{1/2}(s,\lambda)$, and $d_{4}(s,\lambda)=D_{4}^{1/2}(s,\lambda)$. Then

\begin{lemma}For $n\geq 0$,
\begin{equation*}\left[-\frac{1}{(2\,n+1)!}\,\frac{\partial^{2\,n+1}}{\partial\,\lambda^{2\,n+1}}
+\frac{1}{(2\,n)!}\,\frac{\partial^{2\,n}}{\partial\,\lambda^{2\,n}}\right]\,d_{1}(s,\lambda)\,\,\bigg{\vert}_{\lambda=1}=\frac{(-1)^{n}}{n!}\,\frac{\partial^{n}}{\partial\,\lambda^{n}}\,d_{4}(s,\lambda)\,\,\bigg{\vert}_{\lambda=1}.\end{equation*}
\label{dlemma}\end{lemma}
\begin{proof}
Let
\begin{equation*}f(s,\lambda)=1-\sqrt{\frac{\lambda}{2-\lambda}}\,\,\tanh{\frac{\mu(s,\tilde{\lambda})}{2}}\end{equation*}
by the previous lemma, we need to show that
\begin{equation}
  \label{eq:38}
  \left[-\frac{1}{(2\,n+1)!}\,\frac{\partial^{2\,n+1}}{\partial\,\lambda^{2\,n+1}} + \frac{1}{(2\,n)!}\,\frac{\partial^{2\,n}}{\partial\,\lambda^{2\,n}}\right]\,d_{4}(s,\tilde{\lambda})\,f(s,\lambda)\,\,\bigg{\vert}_{\lambda=1}
\end{equation}
equals
\begin{equation*}
\frac{(-1)^{n}}{n!}\,\frac{\partial^{n}}{\partial\,\tilde{\lambda}^{n}}\,d_{4}(s,\tilde{\lambda})\,\,\bigg{\vert}_{\lambda=1}.
\end{equation*} Now formula \eqref{faa} applied to $d_{4}(s,\tilde{\lambda})$
gives
\begin{equation*}\frac{\partial^{k}}{\partial\,\lambda^{k}}\,d_{4}(s,\tilde{\lambda})\,\,\bigg{\vert}_{\lambda=1}=
\begin{cases}
          \quad 0 &\text{if $k=2j+1$, $j\geq 0$,}\\
           \quad\frac{(-1)^{j}\,k!}{j!} \frac{\partial^{j}}{\partial\,\tilde{\lambda}^{j}}\,d_{4}(s,\tilde{\lambda}) &\text{if $k=2j$, $j\geq 0$.}  \\
         \end{cases} \end{equation*}
Therefore
\begin{align*}
-\frac{1}{(2\,n+1)!}\,\frac{\partial^{2\,n+1}}{\partial\,\lambda^{2\,n+1}}&\,d_{4}(s,\tilde{\lambda})\,f(s,\lambda)\,\,\bigg{\vert}_{\lambda=1}\\
&=-\frac{1}{(2\,n+1)!}\,\sum_{k=0}^{2\,n+1}\binom{2\,n+1}{k}\,\frac{\partial^{k}}{\partial\,\lambda^{k}}\,d_{4\,}
\frac{\partial^{2\,n+1-k}}{\partial\,\lambda^{2\,n+1-k}}\,f
\,\,\bigg{\vert}_{\lambda=1}\\
&=
-\sum_{j=0}^{n}\frac{(-1)^{j}}{(2\,n-2\,j+1)!\,j!}\,\frac{\partial^{j}}{\partial\,\tilde{\lambda}^{j}}\,d_{4\,}
\frac{\partial^{2\,n-2\,j+1}}{\partial\,\lambda^{2\,n-2\,j+1}}\,f
\,\,\bigg{\vert}_{\lambda=1}
\end{align*}
Similarly,
\begin{align*}
  \frac{1}{(2\,n)!}\,\frac{\partial^{2\,n}}{\partial\,\lambda^{2\,n}}\,d_{4}(s,\tilde{\lambda})\,f(s,\lambda)\,\,\bigg{\vert}_{\lambda=1}
  &=\frac{1}{(2\,n)!}\,\sum_{k=0}^{2\,n}\binom{2\,n}{k}\,\frac{\partial^{k}}{\partial\,\lambda^{k}}\,d_{4\,}\frac{\partial^{2\,n-k}}{\partial\,\lambda^{2\,n-k}}\,f\,\,\bigg{\vert}_{\lambda=1}\\
  &=\sum_{j=0}^{n}\frac{(-1)^{j}}{(2\,n-2\,j)!\,j!}\,\frac{\partial^{j}}{\partial\,\tilde{\lambda}^{j}}\,d_{4\,}\frac{\partial^{2\,n-2\,j}}{\partial\,\lambda^{2\,n-2\,j}}\,f\,\,\bigg{\vert}_{\lambda=1}.
\end{align*}
Therefore the expression in \eqref{eq:38} equals
\begin{align*}
  \sum_{j=0}^{n}\frac{(-1)^{j}}{(2\,n-2\,j)!\,j!}\,\frac{\partial^{j}}{\partial\,\tilde{\lambda}^{j}}\,d_{4\,}(s,\tilde{\lambda}) \left[\frac{\partial^{2\,n-2\,j}}{\partial\,\lambda^{2\,n-2\,j}}\,f- \frac{1}{2\,n-2\,j+1}\,\frac{\partial^{2\,n-2\,j+1}}{\partial\,\lambda^{2\,n-2\,j+1}}\,f\right]\,\,\bigg{\vert}_{\lambda=1}.
\end{align*}
\begin{center}
{\large
\begin{table}
\label{tab:1}
\newcommand\T{\rule{0pt}{2.6ex}}
\newcommand\B{\rule[-1.2ex]{0pt}{0pt}}
\begin{tabular}[htb]{c|c|c|c|c|}
\cline{2-5}
\T \B & \multicolumn{4}{c|}{\textbf{Statistic}} \\
\hline
\multicolumn{1}{|c|}{\textbf{Eigenvalue}}  \T \B & $\mu$ & $\sigma$ & $\gamma_{1}$ & $\gamma_{2}$ \\
\hline
\multicolumn{1}{|c|}{$\lambda_{1}$} \T & $-1.206548$ & $1.267941$ & $0.293115$ & $0.163186$ \\[6pt]
\multicolumn{1}{|c|}{$\lambda_{2}$} \T & $-3.262424$ & $1.017574$ & $0.165531$ & $0.049262$ \\[6pt]
\multicolumn{1}{|c|}{$\lambda_{3}$} \T & $-4.821636$ & $0.906849$ & $0.117557$ & $0.019506$ \\[6pt]
\multicolumn{1}{|c|}{$\lambda_{4}$} \T \B & $-6.162036$ & $0.838537$ & $0.092305$ & $0.007802$ \\
\hline
\end{tabular}
\caption{Mean, standard deviation, skewness and kurtosis data for first four edge--scaled eigenvalues in the $\beta=1$ Gaussian ensemble. Corollary~\ref{interlacingcor} implies that rows $2$ and $4$ give corresponding data for the two largest eigenvalues in the $\beta=4$ Gaussian ensemble. Compare to Table~$1$ in \cite{Trac4}, keeping in mind that the discrepancy in the $\beta=4$ data is caused by the different normalization in the definition of $F_{4}(s,1)$; see footnote in the Introduction, as well as the comments following Equation~\eqref{eq:26}, and Section~\ref{sec:stand-devi-matt}.}
\end{table}
}
\end{center}
Now Lemma~\ref{flemma} shows that the square bracket inside the
summation is zero unless $j=n$, in which case it is $1$. The
result follows.
\end{proof}
In an inductive proof of Corollary~\ref{interlacingcor}, the base case $F_{4}(s,2)=F_{1}(s,1)$ is easily checked by direct calculation. Lemma~\ref{dlemma} establishes the inductive step in the proof since, with the assumption $F_{4}(s,n)=F_{1}(s,2\,n)$, it is equivalent to the statement
\begin{equation*}
  F_{4}(s,n+1)=F_{1}(s,2\,n+2).
\end{equation*}
\begin{center}
{\large
\begin{table}
\label{tab:2}
\newcommand\T{\rule{0pt}{2.6ex}}
\newcommand\B{\rule[-1.2ex]{0pt}{0pt}}
\begin{tabular}[htb]{c|c|c|c|c|}
\cline{2-5}
\T \B & \multicolumn{4}{c|}{\textbf{Statistic}} \\
\hline
\multicolumn{1}{|c|}{\textbf{Eigenvalue}}  \T \B & $\mu$ & $\sigma$ & $\gamma_{1}$ & $\gamma_{2}$ \\
\hline
\multicolumn{1}{|c|}{$\lambda_{1}$} \T & $-1.771087$ & $0.901773$ & $0.224084$ & $0.093448$ \\[6pt]
\multicolumn{1}{|c|}{$\lambda_{2}$} \T \B & $-3.675440$ & $0.735214$ & $0.125000$ & $0.021650$ \\
\hline
\end{tabular}
\caption{Mean, standard deviation, skewness and kurtosis data for first two edge--scaled eigenvalues in the $\beta=2$ Gaussian ensemble. Compare to Table~$1$ in \cite{Trac4}.}
\end{table}
}
\end{center}


\section{The Wishart ensembles}

Consider an $n\times p$ data matrix $X$ whose rows $x_{i}$ are independent Gaussian $N_{p}(0,\Sigma)$. The product $X^{t}\,X$, of $X$ and its transpose is (up to factor $1/n$) a sample covariance matrix. In this setting, the matrices $A=X^{t}\,X$ are said to have Wishart distribution $W_{p}(n,\Sigma)$. The so-called ``Null Case'' corresponds to $\Sigma=\textrm{Id}_{n}$. Let $\lambda_{1}>\cdots >\lambda_{p}$ be eigenvalues of $X^{t}\,X$ and define
\begin{equation*}
  \mu_{_{np}}=\left(\sqrt{n-1}+\sqrt{p}\right)^{2}\quad,\quad \sigma_{_{np}}=\left(\sqrt{n-1}+\sqrt{p}\right)\left(\frac{1}{\sqrt{n-1}}+\frac{1}{\sqrt{p}}\right)^{1/3}
\end{equation*}
Johnstone proved the following theorem (see \cite{John1}).
\begin{thm}
\label{sec:wishart-ensembles-3}
  If $n/p\to\gamma\geq 0$ then
  \begin{equation*}
    \frac{\lambda_{1}-\mu_{_{np}}}{\sigma_{_{np}}}\xrightarrow{\mathscr{D}}F_{1}(s,1).
  \end{equation*}
\end{thm}
Soshnikov generalized this to the $m^{th}$ largest eigenvalue (see\cite{Sosh2}).
\begin{thm}
\label{sec:wishart-ensembles-4}
  If $n/p\to\gamma\geq 0$ then
  \begin{equation*}
    \frac{\lambda_{m}-\mu_{_{np}}}{\sigma_{_{np}}}\xrightarrow{\mathscr{D}}F_{1}(s,m).
  \end{equation*}
\end{thm}
Note that Soshnikov redefines $\mu_{_{np}},\sigma_{_{np}}$ by letting $n\to n+1$ in Johnstone definition, but this does not affect the limiting distribution. El Karoui proves \ref{sec:wishart-ensembles-3} for $0\leq \gamma\leq \infty$ (see \cite{Elka1}). It is an open problem to show Soshnikov's theorem~\ref{sec:wishart-ensembles-4} holds in this generality, but if it does then our results will also have wider applicability. Soshnikov went further and disposed of the Gaussian assumption, showing the universality of this result. Let us redefine the matrices $X=\left\{x_{i,j}\right\}$ to satisfy
\begin{itemize}
\item $\textrm{\textbf{E}}x_{ij}=0$, $\textrm{\textbf{E}}(x_{ij})^{2}=1$.
\item the r.v's $x_{ij}$ have symmetric laws of distribution
\item all moments of these r.v's are finite
\item the distributions of the $x_{ij}$ decay at least as fast as  a Gaussian at infinity: $\textrm{\textbf{E}}(x_{ij})^{2m}\leq (\textrm{const}\,m)^{m}$
\item $n-p=O(P^{1/3})$
\end{itemize}
With these new assumptions, Soshnikov proved the following (see\cite{Sosh2}).
\begin{thm}
  If $n/p\to\gamma\geq 0$ then
  \begin{equation*}
    \frac{\lambda_{m}-\mu_{_{np}}}{\sigma_{_{np}}}\xrightarrow{\mathscr{D}}F_{1}(s,m).
  \end{equation*}
\end{thm}
Our results are immediately applicable here and thus explicitly give the distribution of the $m^{th}$ largest eigenvalue of the appropriate Wishart distribution.
Table~\ref{table} shows a comparison of percentiles of the $F_{1}$ distribution to corresponding percentiles of empirical Wishart distributions. Here $\lambda_{i}$ denotes the $i^{th}$ largest eigenvalue in the Wishart Ensemble. The percentiles in the $\lambda_{i}$ columns were obtained by finding the ordinates corresponding to the $F_{1}$--percentiles listed in the first column, and computing the proportion of eigenvalues lying to the left of that ordinate in the empirical distributions for the $\lambda_{i}$. The bold entries correspond to the levels of confidence most commonly used in statistical applications. The reader should compare this table to a similar one in \cite{John1}.

\begin{center}
{\large
\begin{table}
\newcommand\T{\rule{0pt}{2.6ex}}
\newcommand\B{\rule[-1.2ex]{0pt}{0pt}}
\begin{tabular}[htb]{c|ccc|ccc|}
\cline{2-7}
\T \B & \multicolumn{3}{c|}{$\mathbf{100\times 100}$} & \multicolumn{3}{c|}{$\mathbf{100\times 400}$} \\
\hline
\multicolumn{1}{|c|}{\textbf{$F_{1}$-Percentile}} \T \B & $\lambda_{1}$ & $\lambda_{2}$ & $\lambda_{3}$ & $\lambda_{1}$ & $\lambda_{2}$ & $\lambda_{3}$ \\
\hline
\multicolumn{1}{|c|}{$0.01$} \T & $0.008$ & $0.005$  & $0.004$  & $0.008$ & $0.006$ & $0.004$ \\[6pt]
\multicolumn{1}{|c|}{$0.05$} & $0.042$ & $0.033$  & $0.025$  & $0.042$ & $0.037$ & $0.032$ \\[6pt]
\multicolumn{1}{|c|}{$0.10$} & $0.090$ & $0.073$  & $0.059$  & $0.088$ & $0.081$ & $0.066$ \\[6pt]
\multicolumn{1}{|c|}{$0.30$} & $0.294$ & $0.268$  & $0.235$  & $0.283$ & $0.267$ & $0.254$ \\[6pt]
\multicolumn{1}{|c|}{$0.50$} & $0.497$ & $0.477$  & $0.440$  & $0.485$ & $0.471$ & $0.455$ \\[6pt]
\multicolumn{1}{|c|}{$0.70$} & $0.699$ & $0.690$  & $0.659$  & $0.685$ & $0.679$ & $0.669$ \\[6pt]
\multicolumn{1}{|c|}{$\mathbf{0.90}$} & $\mathbf{0.902}$  & $\mathbf{0.891}$  & $\mathbf{0.901}$ & $\mathbf{0.898}$ & $\mathbf{0.894}$ & $\mathbf{0.884}$ \\[6pt]
\multicolumn{1}{|c|}{$\mathbf{0.95}$} & $\mathbf{0.951}$  & $\mathbf{0.948}$  & $\mathbf{0.950}$ & $\mathbf{0.947}$ & $\mathbf{0.950}$ & $\mathbf{0.941}$ \\[6pt]
\multicolumn{1}{|c|}{$\mathbf{0.99}$} \B & $\mathbf{0.992}$  & $\mathbf{0.991}$  & $\mathbf{0.991}$ & $\mathbf{0.989}$ & $\mathbf{0.991}$ & $\mathbf{0.989}$ \\
\hline
\end{tabular}
\caption{Percentile comparison of $F_{1}$ vs. empirical
distributions for $100\times 100$ and $100\times 400$ Wishart
matrices with identity covariance.}\label{table}
\end{table}
}
\end{center}

\chapter{Numerics}

\section{Partial derivatives of $q(x,\lambda)$}

Let \begin{equation}q_{n}(x)=\frac{\partial^{n}}{\partial\lambda^{n}}\,q(x,\lambda)\bigg\vert_{\lambda=1},\end{equation}
so that $q_{0}$ equals $q$ from \eqref{pII}.  In order to compute
$F_{\beta}(s,m)$ it is crucial to know $q_{n}$ with $0\leq n\leq m$ accurately.
Asymptotic expansions for $q_{n}$ at $-\infty$ are given in
\cite{Trac3}. We outline how to compute $q_{0}$ and $ q_{1}$ as an
illustration. From \cite{Trac3}, we know that, as $t\to+\infty$, $q_{0}(-t/2)$ is given by
\begin{equation*}
  \frac{1}{2}\sqrt{t}\left(1-\frac{1}{t^{3}}-\frac{73}{2t^{6}}-\frac{10657}{2t^{9}}-\frac{13912277}{8t^{12}}+\textrm{O}\left(\frac{1}{t^{15}}\right)\right),
\end{equation*}
whereas $q_{1}(-t/2)$ can be expanded as
\begin{equation*}
  \frac{\exp{(\frac{1}{3}t^{3/2})}}{2\sqrt{2\pi}\,t^{1/4}}\left(1+\frac{17}{24t^{3/2}}+\frac{1513}{2^{7}3^{2}t^{3}}+\frac{850193}{2^{10}3^{4}t^{9/2}}-\frac{407117521}{2^{15}3^{5}t^{6}}+\textrm{O}\left(\frac{1}{t^{15/2}}\right)\right).
\end{equation*}
Quantities needed to compute $F_{\beta}(s,m), m=1,2,$ are not only
$q_{0}$ and $q_{1}$ but also integrals involving $q_{0}$, such as
\begin{equation}
I_{0}=\int_{s}^{\infty}(x-s)\,q_{0}^{2}(x)\,d\,x,\quad
J_{0}=\int_{s}^{\infty}q_{0}(x)\,d\,x.
\end{equation}
Instead of computing these integrals afterward, it is better to
include them as variables in a system together with $q_{0}$, as
suggested in \cite{Pers1}. Therefore all quantities needed are
computed in one step, greatly reducing errors, and taking full
advantage of the powerful numerical tools in MATLAB\texttrademark\,. Since
\begin{equation}
I_{0}'=-\int_{s}^{\infty}q_{0}^{2}(x)\,d\,x,\quad
I_{0}''=q_{0}^{2},\qquad J_{0}'=-q_{0},
\end{equation}
the system closes, and can be concisely written
\begin{equation}
\frac{d}{ds}\left(\begin{array}{c}q_{0} \\ q_{0}' \\ I_{0} \\ I_{0}' \\
J_{0}
\end{array}\right) =\left(\begin{array}{c}q_{0}' \\ s\,q_{0}+2q_{0}^3 \\ I_{0}' \\ q_{0}^{2} \\ -q_{0}
\end{array}\right).
\label{q0sys}\end{equation} We first use the MATLAB\texttrademark\, built--in
Runge--Kutta--based ODE solver \texttt{ode45} to obtain a first
approximation to  the solution of \eqref{q0sys} between $x=6$, and
$x=-8$, with an initial values obtained using the Airy function on
the right hand side. Note that it is not possible to extend the
range to the left due to the high instability of the solution a
little after $-8$; (This is where the transition region between
the three different regimes in the so--called ``connection
problem'' lies. We circumvent this limitation by patching up our
solution with the asymptotic expansion to the left of $x=-8$.).
The approximation obtained is then used as a trial solution in the
MATLAB\texttrademark\, boundary value problem solver \texttt{bvp4c}, resulting in
an accurate solution vector between $x=6$ and $x=-10$.

Similarly, if we define
\begin{equation}
I_{1}=\int_{s}^{\infty}(x-s)\,q_{0}(x)\,q_{1}(x)\,d\,x,\quad
J_{1}=\int_{s}^{\infty}q_{0}(x)\,q_{1}(x)\,d\,x,
\end{equation}
then we have the first--order system
\begin{equation}
\frac{d}{ds}\left(\begin{array}{c}q_{1} \\ q_{1}' \\ I_{1} \\ I_{1}' \\
J_{1}
\end{array}\right) =\left(\begin{array}{c}q_{1}' \\ s\,q_{1}+6q_{0}^2\,q_{1} \\ I_{1}' \\ q_{0}\,q_{1} \\ -q_{0}\,q_{1}
\end{array}\right),
\end{equation}
which can be implemented using \texttt{bvp4c} together with a
``seed'' solution obtained in the same way as for $q_{0}$.

\section{MATLAB\texttrademark\, code}
\label{sec:matl-code}
In this section we gather the MATLAB\texttrademark\, code used to evaluate and plot the distributions. It is broken down into twenty-five functions/files for easier use. They are implemented in version 7 of MATLAB\texttrademark\,, and may not be compatible with earlier versions. The function ``twmakenew'' should be ran first. It initializes everything and creates some data files which will be needed. There will be a few warning messages that can be safely ignored. Once that is done, the distributions are ready to be used. The only functions directly needed to evaluate them are ``twdens'' and ``twdist'' which, as their names indicate, give the Tracy-Widom density and distribution functions. They both take 3 arguments; the first is ``beta'' which is the beta of RMT so it can be 1, 2 or 4; then ``n'' which is the eigenvalue need; finally ``s'' which is the value where you want to evaluate the function. All these functions are vectorized, meaning that they can take a vector as well as scalar argument. This is convenient for plotting.
One does not need to run ``twmakenew'' every time you want to work with these functions. In fact one should only need to run it once ever, unless some files are deleted or moved. Once it is ran the first time, it creates data files that  contain all you need to evaluate the functions. On subsequent uses one can just start by running ``tw'' instead, which will initialize by reading all the data that ``twmakenew'' produced the first time it ran. This is much faster. The only reason for running ``twmakenew'' again is to clear everything and refresh all the data. Otherwise one can just start with ``tw'' and then start using ``twdens'' and ``twdist''. As a summary here is a list of commands that will plot the density of the $3^{rd}$ GOE eigenvalue between $0$ and $-9$:
\begin{verbatim}
twmakenew
s=-9:0.1:0
plot(s, twdens(1,3,s))
\end{verbatim}

The MATLAB\texttrademark\, functions needed are below. Note that the digits at the beginning of each line are not part of the code and are added only for (ease of) reference purposes. They should be omitted when coding. These functions are available electronically by writing to the author.

\vfill\eject

\begin{verbatim}
%
%       (c) 2004-2005 Momar Dieng, All rights reserved.
%
% Comments? e-mail momar@math.ucdavis.edu
%  This software was developed with support of the
%      National Science Foundation/ Grant DMS--0304414

%  In EVERY case, the software is COPYRIGHTED by the original author.
%      For permissions, see below.
%
%  Comments/Permissions/Information: e-mail
%      momar@math.ucdavis.edu or tracy@math.ucdavis.edu
%
% COPYLEFT NOTICE.
% Permission is granted to make and distribute verbatim copies of this
% entire software package, provided that all files are copied
% **together** as a **unit**.
%
% Here **copying as a unit** means: all the files listed  are copied.
% Permission is granted to make and distribute modified copies of this
% software, under the conditions for verbatim copying, provided that
% the entire resulting derived work is distributed under the terms
% of a permission notice identical to this one. Names of new authors,
% their affiliations and information about their improvements may be
% added to the files.
%
% The purpose of this permission notice is to allow you to make
% copies of the software and distribute them to others, for free
% or for a fee, subject to the constraint that you maintain the
% collection of tools here as a unit.  This enables people you
% give the software to to be aware of its origins, to ask questions
% of us by e-mail, to request improvements, obtain later releases,
% etc ...
%
% If you seek permission to copy and distribute translations of
% this software into another language, please e-mail a specific
% request to momar@math.ucdavis.edu or tracy@math.ucdavis.edu.
%
% If you seek permission to excerpt a **part** of the software for
% example to appear in a scientific publication, please e-mail a
% specific request to momar@math.ucdavis.edu or tracy@math.ucdavis.edu.
%

\end{verbatim}
\listinginput{1}{twnext/q0bcjac.m}
\bigskip
\listinginput{1}{twnext/q0bc.m}
\bigskip
\listinginput{1}{twnext/q0bvpinit.m}
\bigskip
\listinginput{1}{twnext/q0jac.m}
\bigskip
\listinginput{1}{twnext/q0leftexp.m}
\bigskip
\listinginput{1}{twnext/q0leftexpprime.m}
\bigskip
\listinginput{1}{twnext/q0.m}
\bigskip
\listinginput{1}{twnext/q0ode.m}
\bigskip
\listinginput{1}{twnext/q0prime.m}
\bigskip
\listinginput{1}{twnext/q0seed.m}
\bigskip
\listinginput{1}{twnext/q1bcjac.m}
\bigskip
\listinginput{1}{twnext/q1bc.m}
\bigskip
\listinginput{1}{twnext/q1bvpinit.m}
\bigskip
\listinginput{1}{twnext/q1jac.m}
\bigskip
\listinginput{1}{twnext/q1leftexp.m}
\bigskip
\listinginput{1}{twnext/q1.m}
\bigskip
\listinginput{1}{twnext/q1ode.m}
\bigskip
\listinginput{1}{twnext/q1prime.m}
\bigskip
\listinginput{1}{twnext/q1seed.m}
\bigskip
\listinginput{1}{twnext/twdens.m}
\bigskip
\listinginput{1}{twnext/twdist.m}
\bigskip
\listinginput{1}{twnext/tw.m}
\bigskip
\listinginput{1}{twnext/twmakenew.m}
\bigskip
\listinginput{1}{twnext/matsim.m}
\bigskip
\listinginput{1}{twnext/matsimplot.m}

\appendix
\chapter{The matrix ensembles}
\label{cha:matrix-ensembles}

As mentioned in the Introduction, the three matrix ensembles discussed in this work were originally (and maybe more naturally) introduced as probability spaces of (self dual) unitary, orthogonal and symplectic matrices matrices (see \cite{Meht1} for more history and an expanded treatment along these lines). We quickly review those set-ups here only insofar as they help motivate the MATLAB\texttrademark\, code we used to generate the eigenvalues, and set the appropriate background for the recent powerful generalizations alluded to earlier (namely in \cite{Dumi1} and \cite{Edel1}).

\section{The Gaussian orthogonal and unitary ensembles}
\label{sec:gauss-orth-unit}

The Gaussian orthogonal (resp. unitary) ensemble GOE (resp. GUE) is defined as a probability space on the set $\mathcal H_{n}$ of $N\times N$ real symmetric (resp. complex hermitian) matrices. If one requires that:
\begin{enumerate}
\item the probability $P(H)\,d\,H$ that a matrix $H$ in $\mathcal H_{n}$ belongs to the volume element $d\,H$ be invariant under similarity transformations $H\to A\,H\,A^{\dag}$ by orthogonal (resp. unitary) matrices,
\item the various $H_{jk}$ elements of the matrix $H$ be statistically independent so that $P(H)$ factors in a product of functions each depending on a single $H_{jk}$,
\end{enumerate}
then $P(H)$ is restricted to the form
\begin{equation}
\label{eq:39}
  P(H)=C\,\exp\left(-a\tr{H^{2}+b\,\tr{H}+c}\right),\quad a>0
\end{equation}
Converting to spectral variables and integrating out the angular part (the Vandermonde is the Jacobian) yields
\begin{equation}
   P_{\beta}(\theta_{1},\ldots,\theta_{N})= P_{\beta}^{(N)}(\vec{\theta}\,) = C\,\exp\left[-\sum_{j=1}^{N}a\,\theta_{j}^{2}-b\,\theta_{j}-c\right]\prod_{j<k}|\theta_{j}-\theta_{k}|^{\beta},
\end{equation}
where $\beta=1$ in the real symmetric case, and $\beta=2$ in the complex hermitian case. If we let $\theta_{j}\to x_{j}\sqrt{\beta/2a}+b/2a$ then
\begin{equation}
  \label{eq:37}
  P_{\beta}(x_{1},\ldots,x_{N})= P_{\beta}^{(N)}(\vec{x}\,) = C_{\beta}^{(N)}\,\exp\left[-\frac{\beta}{2}\,\sum_{j=1}^{N}x_{j}^{2}\right]\prod_{j<k}|x_{j}-x_{k}|^{\beta}.
\end{equation}
with \eqref{eq:33} yielding
\begin{equation}
  \label{eq:35}
  C_{\beta}^{(N)} = (2\,\pi)^{-N/2}\,\beta^{-N/2-\beta\,N(N-1)/4}\cdot\prod_{j=1}^{N}\frac{\Gamma(1 + \gamma)\,\Gamma(1 + \frac{\beta}{2})}{\Gamma(1 + \frac{\beta}{2}\,j)}
\end{equation}
The following MATLAB\texttrademark\, code produces a GOE--distributed matrix $A$ of size $n\times n$ with diagonal entries distributed $N(0,\sigma^{2})$,  and off diagonal entries distributed as $N(0,\sigma^{2}/2)$. It is adapted from \cite{Dumi1}; see \cite{Edel1}, and the code listing for function \textsc{matsim} in Section~\ref{sec:matl-code} for \underline{much} more efficient schemes.
\begin{verbatim}
A = sigma*randn(n);
A = (A + A')/2;
\end{verbatim}
whereas the lines below produce a GUE--distributed matrix $A$ of size $n\times n$ also with diagonal entries distributed $N(0,\sigma^{2})$,  and off diagonal entries whose real and imaginary parts are distributed as $N(0,\sigma^{2}/2)$:
\begin{verbatim}
X = sigma*randn(n) + i*sigma*randn(n);
Y = sigma*randn(n) + i*sigma*randn(n);
A = (A + A')/2;
\end{verbatim}

\section{The Gaussian symplectic ensemble}
\label{sec:gauss-sympl-ensemble}

The GSE is defined as a probability space on the set of all $2n\times 2n$ hermitian matrices which are self--dual (in a way to be defined soon) when considered as $n\times n$ quaternion matrices. Arguments similar to those in the previous section lead to the representations \eqref{eq:39} for $P(H)$ and \eqref{eq:37} for the joint probability distribution of eigenvalues with $\beta$ set to $4$ this time (see \cite{Meht1} for details). The GSE is best defined in terms of quaternion matrices, since this makes its similarity to the GUE and GSE plain, and yields a MATLAB\texttrademark\, algorithm to produce the matrices. We expand on the construction of the ensemble here in order to motivate the MATLAB\texttrademark\, code. A good introduction to quaternions, and the classical groups can be found in \cite{Chev1}, which we use extensively in the following. Consider the algebra $\mathcal Q$ of dimension $4$ over $\mathbb R$ generated by the basis elements $e_{0}, e_{1}, e_{2}, e_{3}$ and relations
\begin{equation}
  e_{0}\,e_{i}=e_{i}\,e_{0}=e_{i}; \qquad e_{i}^{2}=-e_{0}; \qquad e_{i}\,e_{j}=-e_{j}\,e_{i}=e_{k}
\end{equation}
where $1\leq i,j,k \leq 3$ and $\left(\begin{array}{ccc} i & j & k \\ 1 & 2 & 3 \end{array}\right)$ is an even permutation. Elements of $\mathcal Q$ can be thought of as ``complex numbers with three imaginary parts''. In fact $\mathcal Q$ contains a sub-algebra isomorphic to the field $\mathbb C$ (the sub-algebra of quaternions of the form $q_{0}\,e_{0} + q_{1}\,e_{1}$, each of which is mapped to the complex number $q_{0} + i\,q_{1}$). Therefore $\mathcal Q$ can be considered as a vector space of dimension $2$ over $\mathbb C$. Indeed every quaternion
\begin{equation}
  q = q_{0}\,e_{0} + q_{1}\,e_{1} + q_{2}\,e_{2}+ q_{3}\,e_{3}
\end{equation}
can be rewritten, using the above rules, as
\begin{equation}
  q = e_{0}(q_{0}\,e_{0} + q_{1}\,e_{1}) + e_{2}(q_{2}\,e_{0} - q_{3}\,e_{1}),
\end{equation}
which under the map to $\mathbb C$ mentioned above becomes
\begin{equation}
  q=e_{0}(q_{0} + i\,q_{1}) + e_{2}(q_{2} - i\,q_{3}).
\end{equation}
Hence every quaternion can be written as a linear combination of the two basis elements $e_{0}$ and $e_{2}$ over $\mathbb C$. It is easy to check that the action of $\mathcal Q$ on itself (considered as a $2$--dimensional $\mathbb C$--vector space) gives a representation $\rho:~\mathcal Q~\to~GL(2,\mathbb C)$ defined by
\begin{equation}
  \begin{array}{ccc}
    \rho(e_{0})= \left(
      \begin{array}{cc}
        1 & 0 \\
        0 & 1
      \end{array}
    \right)
    &  & \rho(e_{1})=
    \left(
      \begin{array}{cc} i & 0 \\
        0 & -i
      \end{array}
    \right) \\  \\ \rho(e_{2})=
    \left(
      \begin{array}{cc} 0 & -1 \\
        1 & 0
      \end{array}
    \right)
    &  & \rho(e_{3})=
    \left(
      \begin{array}{cc} 0 & -i \\
        -i & 0
      \end{array}
    \right)
  \end{array}.
\end{equation}
This representation is often taken as the starting definition for the basis of $\mathcal Q$. In this language, $\mathcal Q$, is then thought of as an algebra (over $\mathbb C$) of matrices with entries in $\mathbb C$. We will adopt that point of view from here on (i.e. we will stop distinguishing between $e_{i}$ and $\rho(e_{i})$ as we have been so far). It follows from the above discussion that a matrix in $GL(2\,n,\mathbb C)$ can be identified with one in $GL(n, \mathcal Q)$ in the following way; think of the complex $2\,n\times 2\,n$ matrix as a block matrix, with $n^{2}$ blocks each made up of a $2\times 2$ matrices of the form
\begin{equation}
  \left(
    \begin{array}{cc}
      a & b \\
      c & d
    \end{array}
  \right),
\end{equation}
where $a, b, c, d\in\mathbb C$. Using the matrix representation of $\mathcal Q$ we can decompose the above $2\times 2$ matrix as
\begin{equation}
  \left(
    \begin{array}{cc} a & b \\
      c & d
    \end{array}
  \right)
  = \frac{1}{2}\,(a + d)\,e_{0} - \frac{i}{2}\,(a - d)\,e_{1} + \frac{1}{2}\,(b - c)\,e_{2} - \frac{i}{2}\,(b + c)\,e_{3}.
\end{equation}
If each block of this form is identified with the quaternion on the right--hand side, the complex $2\,n\times 2\,n$ matrix then becomes a quaternionic $n \times n$ matrix. We can define three distinct analogues of usual complex conjugation in $\mathcal Q$. Let
\begin{equation}
  q = q_{0}\,e_{0} + q_{1}^{}\,e_{1} + q_{2}\,e_{2}+ q_{3}\,e_{3}.
\end{equation}
Then following \cite{Meht1} we define
\begin{eqnarray}
  \overline{q} & = & q_{0}\,e_{0} - q_{1}^{}\,e_{1} - q_{2}\,e_{2} - q_{3}\,e_{3}\\
  q^{*} & = & q_{0}^{*}\,e_{0} + q_{1}^{*}\,e_{1} + q_{2}^{*}\,e_{2} + q_{3}^{*}\,e_{3}\\
  q^{\dag} & = & \overline{q}^{\,*}
\end{eqnarray}
where $^{*}$ denotes usual complex conjugation (recall $q_{i}\in\mathbb C$). We extend this notation to matrices in $\mathcal Q^{n}$ in the usual fashion: for example if $Q=\{q_{jk}\}\in \mathcal Q^{n}$ then $\overline{Q}=\{\overline{q}_{jk}\}$ and $Q^{\dag}=\{q^{\dag}_{kj}\}$. Let us determine what the usual matrix transformations (such as transposition and usual complex conjugation) of a matrix $Q\in GL(2\,n,\mathbb C)$ correspond to when $Q$ is instead considered as a matrix in $\mathcal Q^{n}$. For example if we transpose $Q$ as a complex matrix, what happens to its quaternion entry $q_{ij}$? We write $\left(Q^{t}\right)_{ij}$ for the $(i,j)$ quaternionic element of the matrix $Q\in\mathcal Q^{n}$ whose transpose is taken when it is is considered as an element of $GL(2\,n,\mathbb C)$. Under transposition, the block
\begin{equation}
  \left(\begin{array}{cc} a & b \\ c & d \end{array}\right)
\end{equation}
in position $(j,k)$ of $Q$ goes to
\begin{equation}
  \left(\begin{array}{cc} a & c \\ b & d \end{array}\right)
\end{equation}
in position $(k,j)$. The first block corresponds to the quaternion
\begin{equation}
  q_{jk} = \frac{1}{2}\,(a + d)\,e_{0} - \frac{i}{2}\,(a - d)\,e_{1} + \frac{1}{2}\,(b - c)\,e_{2} - \frac{i}{2}\,(b + c)\,e_{3}.
\end{equation}
The second block is the quaternion
\begin{eqnarray}
  \left(Q^{t}\right)_{kj} & = & \frac{1}{2}\,(a + d)\,e_{0} - \frac{i}{2}\,(a - d)\,e_{1} + \frac{1}{2}\,(c - b)\,e_{2} - \frac{i}{2}\,(b + c)\,e_{3} \nonumber\\
  & = & -e_{2}\,\left(\frac{1}{2}\,(a + d)\,e_{0} + \frac{i}{2}\,(a - d)\,e_{1} - \frac{1}{2}\,(b - c)\,e_{2} + \frac{i}{2}\,(b + c)\,e_{3}\right)\,e_{2} \nonumber\\
  & = & -e_{2}\,\overline{q}_{jk}\,e_{2}.\label{eq:27}
\end{eqnarray}
So transposing $Q$ as a complex matrix is the equivalent of left  and right multiplying by $e_{2}$ the quaternionic transpose of $-\overline{Q}$.  Similarly, doing this for regular hermitian conjugation, we see that if
\begin{equation}
  q_{jk} = \left(\begin{array}{cc} a & b \\ c & d \end{array}\right) = \frac{1}{2}\,(a + d)\,e_{0} - \frac{i}{2}\,(a - d)\,e_{1} + \frac{1}{2}\,(b - c)\,e_{2} - \frac{i}{2}\,(b + c)\,e_{3},
\end{equation}
then
\begin{eqnarray}
  \left(Q^{\dag}\right)_{kj} & = & \left(\begin{array}{cc} \overline{a} & \overline{c} \nonumber\\
      \overline{b} & \overline{d} \end{array}\right) \nonumber\\
  & = & \frac{1}{2}\,(\overline{a} + \overline{d})\,e_{0} - \frac{i}{2}\,(\overline{a} - \overline{d})\,e_{1} + \frac{1}{2}\,(\overline{c} - \overline{b})\,e_{2} - \frac{i}{2}\,(\overline{b} + \overline{c})\,e_{3}\nonumber\\
  & = & \overline{\frac{1}{2}\,(a + d)}\,e_{0} + \overline{\frac{i}{2}\,(a - d)}\,e_{1} - \overline{\frac{1}{2}\,(b - c)}\,e_{2} + \overline{ \frac{i}{2}\,(b + c)}\,e_{3}\nonumber\\
  & = & q^{\dag}_{jk}.\label{eq:28}
\end{eqnarray}
Finally we define the dual $Q^{R}$ of $Q\in\mathcal Q^{n}$ to be
\begin{equation}
  Q^{R}=e_{2}\,\left(Q^{t}\right)\,e_{2}^{-1}.
\end{equation}
Using the relations found above, we see that
\begin{equation}
  \left(Q^{R}\right)_{kj}=e_{2}\,\left(Q^{t}\right)_{kj}\,e_{2}^{-1}=e_{2}\,\left(-e_{2}\,\overline{q}_{jk}\,e_{2}\right)\,e_{2}^{-1}=\overline{q}_{jk}.\label{eq:29}
\end{equation}
So this is the true analog in $\mathcal Q^{n}$ of transpose conjugation of complex matrices. Now we can define the GSE properly
\begin{defn}
The GSE consist of all $2n\times 2n$ hermitian matrices which are self--dual when considered as $n\times n$ quaternion matrices. Equivalently, these are $2n\times 2n$ hermitian matrices $B$ that are invariant under $B\to Z\,B^{t}\,Z^{-1}$ where $Z$ is $2n\times 2n$ block diagonal matrix with $n$ blocks $e_{2}$ down the diagonal.
\end{defn}
Let the elements of $H$ be denoted by $\{h_{rs}\}$ if $H$ is viewed in $GL(2\,n,\mathbb C)$, and by $\{q_{jk}\}$ if it is viewed in $\mathcal Q^{n}$ so that $1\leq r,s\leq 2n$ and $1\leq j,k\leq n$. Furthermore, in order to be able to specify the individual quaternion components of the matrix elements, let us switch to the notation
\begin{equation}
  q_{jk} = q_{jk}^{(0)}\,e_{0} + q_{jk}^{(1)}\,e_{1} + q_{jk}^{(2)}\,e_{2}+ q_{jk}^{(3)}\,e_{3}
\end{equation}
Then $H$ is hermitian self--dual if
\begin{equation}
  H^{\dag}=H^{R}=H,
\end{equation}where we abuse notation since technically the first and second operations on $H$ are in different spaces. It follows from that in terms of the quaternion entries $q_{jk}$ of $H$, being self--dual implies
\begin{equation}
  q^{\dag}_{kj}=\overline{q}_{kj}=q_{jk}\label{eq:30}
\end{equation}
Combining \eqref{eq:30} with \eqref{eq:28}, and \eqref{eq:29} yields
\begin{equation}
  q_{jk}^{(0)}=q_{kj}^{(0)}=\left(q_{kj}^{(0)}\right)^{*}
\end{equation}
The last equality implies reality, and the first symmetry. Therefore if $H$ is hermitian self--dual, $\left\{q_{jk}^{(0)}\right\}, 1\leq j,k\leq n$ forms a real symmetric sub-matrix. Similarly for $l=1,2,3$ we get
\begin{equation}
  q_{jk}^{(l)}=-q_{kj}^{(l)}=-\left(q_{kj}^{(l)}\right)^{*}
\end{equation}
Again the last equality implying reality, and the first antisymmetry. Therefore if $H$ is hermitian self--dual, $\left\{q_{jk}^{(l)}\right\},1\leq j,k\leq n$ forms a real antisymmetric sub-matrix for each $l$.
Finally we need to know how to go back from the $q_{jk}^{(l)}$ to the $h_{rs}$. It is easy to see that if
\begin{equation}
  \left(\begin{array}{cc} a & b \\ c & d \end{array}\right)
\end{equation}
corresponds to
\begin{equation}
  q_{jk} = \frac{1}{2}\,(a + d)\,e_{0} - \frac{i}{2}\,(a - d)\,e_{1} + \frac{1}{2}\,(b - c)\,e_{2} - \frac{i}{2}\,(b + c)\,e_{3},
\end{equation}
then
\begin{equation}
  q_{jk} = q_{jk}^{(0)}\,e_{0} + q_{jk}^{(1)}\,e_{1} + q_{jk}^{(2)}\,e_{2}+ q_{jk}^{(3)}\,e_{3}
\end{equation}
corresponds to
\begin{equation}
  \left(\begin{array}{cc} a & b \\ c & d \end{array}\right)
\end{equation}
with
\begin{eqnarray}
  a & =q_{jk}^{(0)} + i\,q_{jk}^{(1)} \\
  b & =q_{jk}^{(2)} + i\,q_{jk}^{(3)}\\
  c & =-\left(q_{jk}^{(2)} - i\,q_{jk}^{(3)}\right)\\
  d & =q_{jk}^{(0)} - i\,q_{jk}^{(1)}\\
\end{eqnarray}
which yields the following recipe to construct a random $2n\times 2n$ GSE matrix $M$:
\begin{enumerate}
\item generate a random $n\times n$ real symmetric matrix $W$ and 3 random $n\times n$ real antisymmetric matrices $X, Y, Z$.
\item the sub-matrix $A$ of $M$ which in our notation corresponds to the entry $a$ in each complex $2\times 2$ block will be given by $A = W + i\, X$. Similarly $B = Y + i\,Z$, $C = -(Y-i\,Z)$, $D = W -  i\,X$
\item then $M = A \otimes I_{1} + B \otimes I_{2} + C \otimes I_{3} + D \otimes I_{4}$
where
\begin{equation}
  I_{k}= \left(\begin{array}{cc} \delta_{1,k} & \delta_{2,k} \\ \delta_{3,k} & \delta_{4,k} \end{array}\right)
\end{equation}

\end{enumerate}
This recipe is simplified by noting that the matrix $M$ above has the same eigenvalues as the matrix
\begin{equation}
  H=\left(\begin{array}{cc} A & B \\ C  & D \end{array}\right).
\end{equation}
Note that if $\overline{A}$ denotes the complex conjugate of $A$, and $A^{t}$ its transpose, then by construction $A=A^{\dag}$, $A=\overline{D}$, $B^{t}=-B$, $-\overline{B}=C$. Hence we can equivalently let
\begin{equation}
  H=\left(\begin{array}{cc} A & B \\ -\overline{B}  & \overline{A} \end{array}\right)
\end{equation}
where $A$ and $B$ are complex matrices satisfy $A^{\dag}=A$ and $B^{t}=-B$. Hence the following MATLAB\texttrademark\, code produces a GSE--distributed matrix $A$ of size $2\,n\times 2\,n$ with diagonal entries distributed $N(0,\sigma^{2})$,  and off diagonal entries distributed as $N(0,\sigma^{2}/2)$. It is adapted from \cite{Dumi1}; see \cite{Edel1}, and the code listing for function \textsc{matsim} in Section~\ref{sec:matl-code} for \underline{much} more efficient schemes.
\begin{verbatim}
X = sigma*randn(n) + i*sigma*randn(n);
Y = sigma*randn(n) + i*sigma*randn(n);
A = [X Y; -conj(Y) conj(X)];
A = (A + A')/2;
\end{verbatim}

\section{Selberg's integral}
\label{sec:selbergs-integral}

For any positive integer $n$ let $d\,x=d\,x_{1}\cdots d\,x_{n}$ and define
\begin{equation*}
  \Delta(x) = \Delta(x_{1}, \ldots,x_{n}) = \prod_{1\leq j<k\leq n}(x_{j} - x_{k}),
\end{equation*}
if $n>1$, and $\Delta(x)=1$ if $n=1$. Also let
\begin{equation*}
  \Phi(x) = \Phi(x_{1}, \ldots,x_{n}) = \left|\Delta(x)\right|^{2\,\gamma}\,\prod_{j=1}^{n}x_{j}^{\alpha-1}\,(1 - x_{j})^{\beta-1}.
\end{equation*}
Then,
\begin{thm}[Selberg, Aomoto]
  \begin{equation}
\label{eq:31}
    I(\alpha,\beta,\gamma,n) = \int_{0}^{1}\cdots\int_{0}^{1}\Phi(x)\,d\,x = \prod_{j=0}^{n-1}\frac{\Gamma(1 + \gamma + j\,\gamma)\,\Gamma(\alpha + j\,\gamma)\,\Gamma(\beta + j\,\gamma)}{\Gamma(1 + \gamma)\,\Gamma(\alpha + \beta + (n+j-1)\,\gamma)},
  \end{equation}
and for $1\leq m \leq n$,
\begin{equation}
  \label{eq:32}
  \int_{0}^{1}\cdots\int_{0}^{1}x_{1}\,x_{2}\cdots x_{m}\,\Phi(x)\,d\,x =\prod_{j=0}^{n-1}\frac{\alpha + (n-1)\,\gamma}{\alpha + \beta + (2\,n - j - 1)\,\gamma}\, \int_{0}^{1}\cdots\int_{0}^{1}\Phi(x)\,d\,x,
\end{equation}
valid for integer $n$, and complex $\,\alpha$, $\beta$, $\gamma$ with $\re{\alpha}>0$, $\re{\beta}>0$, and
\begin{equation*}
  \re{\gamma}>-\min\left(\frac{1}{n}, \frac{\re{\alpha}}{n-1},  \frac{\re{\beta}}{n-1}\right)
\end{equation*}
\end{thm}
The identity \eqref{eq:31} is the one commonly referred to as Selberg's integral, and \eqref{eq:32} is Aomoto's extension of it. See \cite{Meht1} for proofs and an expanded discussion of its consequences and applications. The following specialization is of particular interest here; put $x_{j}=y_{j}/L$, $\alpha=\beta=a\,L^{2}+1$, and take the limit $L\to\infty$ in \eqref{eq:31} to obtain
\begin{align}
   \int_{-\infty}^{\infty}\cdots & \int_{-\infty}^{\infty}\left|\Delta(x)\right|^{2\,\gamma}\,\prod_{j=1}^{n}\exp\left(-a\,x_{j}^{2}\right)\,d\,x_{j} \nonumber\\ & =(2\,\pi)^{n/2}\,(2\,a)^{-n(\gamma\,(n-1)+1)/2}\cdot\prod_{j=1}^{n}\frac{\Gamma(1 + j\,\gamma)}{\Gamma(1 + \gamma)\,\Gamma(1 + \gamma)}.\label{eq:33}
\end{align}
Set $\gamma=\beta/2$ and $a=1$ to justify \eqref{eq:34}, whereas $a=\beta/2$ gives \eqref{eq:35}.

\section{Standard deviation matters}
\label{sec:stand-devi-matt}
As mentioned in Section~\ref{sec:gauss-spec}, the choice of standard deviation is not the same in all three ensembles. This affects the way we scale data obtained in MATLAB\texttrademark\, simulations with the code give in the previous sections in this appendix. We treated the $\beta=2,4$ cases with weight function $\exp(-x^{2})$ corresponding to a standard deviation $\sigma_{d}=1/\sqrt{2}$ on the diagonal matrix elements, and the scaling
\begin{equation*}
  t=\sqrt{2\,N}+\frac{s}{\sqrt{2}\,N^{1/6}}.
\end{equation*}
Therefore if we want to rescale data with an arbitrary standard deviation $\sigma_{d}$ on diagonal elements, we need to let $t\to t/\sqrt{2}\sigma_{d}$ resulting in the more general change of variables
\begin{equation*}
  t=2\,\sigma_{d}\,\sqrt{N}+\frac{\sigma_{d}\,s}{N^{1/6}},\quad\beta=2,4.
\end{equation*}
In the $\beta=1$ case however our choice of weight function $\exp(-x^{2}/2)$ corresponds to $\sigma_{d}=1$. A similar argument would thus yield
\begin{equation*}
  t=\sigma_{d}\,\sqrt{2\,N}+\frac{\sigma_{d}\,s}{\sqrt{2}\,N^{1/6}},\quad\beta=1.
\end{equation*}
In order to preserve some uniformity in the formulas, we can use the standard deviation of the off diagonal elements instead in the $\beta=1$ case. Indeed, if we denote by $\sigma_{o}$ the standard deviation on the off--diagonal elements then $\sigma_{d}=\sqrt{2}\,\sigma_{o}$ by construction in all three ensembles. Hence in our treatment, $\sigma_{o}=1/\sqrt{2}$ in the $\beta=1$ case, and for general $\sigma_{o}$ we have, similar to the $\beta=2,4$ change of variables
\begin{equation*}
  t=2\,\sigma_{o}\,\sqrt{N}+\frac{\sigma_{o}\,s}{N^{1/6}},\quad\beta=1.
\end{equation*}

\begin{thebibliography}{10}

\bibitem{Andr1}
G.~E. Andrews, Askey R., and R.~Ranjan.
\newblock {\em {Special Functions}}, volume~71 of {\em {Encyclopedia of
  Mathematics and its Applications}}.
\newblock Cambridge University Press, 2000.

\bibitem{Baik2}
J.~Baik and E.~M. Rains.
\newblock {Algebraic aspects of increasing subsequences}.
\newblock {\em Duke Math. J.}, 109(1):1--65, 2001.

\bibitem{Bowi1}
M.~J. Bowick and E.~Br\'ezin.
\newblock {Universal scaling of the tail of the density of eigenvalues in
  random matrix models}.
\newblock {\em Phys. Lett.}, B268:21--28, 1991.

\bibitem{Chev1}
Claude Chevalley.
\newblock {\em {Theory of Lie Groups}}.
\newblock Princeton University Press, 1946.

\bibitem{Clar1}
P.~A. Clarkson and J.~B. McLeod.
\newblock {A connection formula for the second Painlev\'e transcendent}.
\newblock {\em Arch. Rational Mech. Anal.}, 103(2):97--138, 1988.

\bibitem{Debr1}
N.~G. de~Bruijn.
\newblock {On some multiple integrals involving determinants}.
\newblock {\em J. Indian Math. Soc.}, 19:133--151, 1955.

\bibitem{Deif2}
P.~Deift and X.~Zhou.
\newblock {Asymptotics for the Painlev\'e II equation}.
\newblock {\em Commun. Pure Appl. Math.}, 48:277--337, 1995.

\bibitem{Dumi1}
I.~Dumitriu and A.~Edelman.
\newblock {Matrix models for beta ensembles}.
\newblock {\em J. Math. Phys.}, 43(11):5830--5847, 2002.

\bibitem{Edel1}
A.~Edelman and N.~R. Rao.
\newblock {Random matrix theory}.
\newblock {\em Acta Numerica}, 14:233--297, 2005.

\bibitem{Elka1}
N.~El~Karoui.
\newblock {On the largest eigenvalue of Wishart matrices with identity
  covariance when $n$, $p$ and $p/n$ tend to infinity}.
\newblock ArXiv:math.ST/0309355.

\bibitem{Forr2}
P.~J. Forrester.
\newblock {The spectrum edge of random matrix ensembles}.
\newblock {\em Nucl. Phys.}, B402:709--728, 1993.

\bibitem{Forr1}
P.~J. Forrester and E.~M. Rains.
\newblock {Interrelationships between orthogonal, unitary and symplectic matrix
  ensembles}.
\newblock In P.~Bleher, A.~Its, and S.~Levy, editors, {\em Random Matrix Models
  and their Applications}, volume~40 of {\em Math. Sci. Res. Inst. Publ.},
  pages 171--207. Cambridge Univ. Press, Cambridge, 2001.

\bibitem{Gohb1}
I.~Gohberg, S.~Goldberg, and M.~A. Kaashoek.
\newblock {\em {Classes of Linear Operators, Vol. I}}, volume~49 of {\em
  Operator Theory: Advances and Applications}.
\newblock Birkh{\"a}user, 1990.

\bibitem{Hasti1}
S.~P. Hastings and J.~B. McLeod.
\newblock {A boundary value problem associated with the second Painlev\'e
  transcendent and the Korteweg--de\thinspace Vries equation}.
\newblock {\em Arch. Rational Mech. Anal.}, 73(1):31--51, 1980.

\bibitem{John1}
I.~M. Johnstone.
\newblock On the distribution of the largest eigenvalue in principal component
  analysis.
\newblock {\em Ann. Stats.}, 29(2):295--327, 2001.

\bibitem{Meht1}
M.~L. Mehta.
\newblock {\em Random Matrices, Revised and Enlarged Second Edition}.
\newblock Academic Press, 1991.

\bibitem{Meht3}
M.~L. Mehta and F.~Dyson.
\newblock {Statistical theory of the energy levels of complex systems. V.}
\newblock {\em J. Math. Phys.}, 4:713--719, 1963.

\bibitem{Pers1}
P.~Persson.
\newblock {Numerical methods for random matrices}.
\newblock {Course project for MIT 18.337}, MIT, 2002.

\bibitem{Sosh2}
A.~Soshnikov.
\newblock {A note on universality of the distribution of the largest
  eigenvalues in certain sample covariance matrices}.
\newblock {\em J. Stat. Phys.}, 108(5--6):1033--1056, 2002.

\bibitem{Stan2}
R.~Stanley.
\newblock {\em {Enumerative Combinatorics}}, volume~2.
\newblock Cambridge University Press, 1999.

\bibitem{Trac7}
C.~A. Tracy and H.~Widom.
\newblock {Fredholm determinants, differential equations and matrix models}.
\newblock {\em Commun. Math. Physics}, 163:33--72, 1994.

\bibitem{Trac3}
C.~A. Tracy and H.~Widom.
\newblock {Level--spacing distributions and the Airy kernel}.
\newblock {\em Commun. Math. Physics}, 159:151--174, 1994.

\bibitem{Trac2}
C.~A. Tracy and H.~Widom.
\newblock {On orthogonal and symplectic matrix ensembles}.
\newblock {\em Commun. Math. Physics}, 177:727--754, 1996.

\bibitem{Trac1}
C.~A. Tracy and H.~Widom.
\newblock {Correlation functions, cluster functions, and spacing distributions
  for random matrices}.
\newblock {\em J. Stat. Phys.}, 92(5--6):809--835, 1998.

\bibitem{Trac4}
C.~A. Tracy and H.~Widom.
\newblock {Airy kernel and Painlev\'e}.
\newblock In {\em Isomonodromic deformations and applications in physics},
  volume~31 of {\em {CRM Proceedings \& Lecture Notes}}, pages 85--98. Amer.
  Math. Soc., Providence, RI, 2002.

\bibitem{Trac5}
C.~A. Tracy and H.~Widom.
\newblock {Matrix kernels for the Gaussian orthogonal and symplectic
  ensembles}.
\newblock ArXiv:math-ph/0405035, 2004.

\end{thebibliography}

\end{document}